\documentclass[brochure,english,12pt]{bourbaki}
\usepackage[matrix,arrow]{xy}
\usepackage{amssymb,amsfonts,amsmath,footnote,url}
\usepackage{csquotes}
\usepackage{tikz}
\usetikzlibrary{shapes,arrows}
\usepackage[Algorithm]{algorithm}
\usepackage[noend]{algpseudocode}
\usepackage[russian,english]{babel}
\usepackage[utf8]{inputenc} 
\usepackage[T1, T2A]{fontenc}    
\newcommand{\El}{\mbox{\usefont{T2A}{\rmdefault}{m}{n}\CYRL}}
\DeclareMathOperator{\Sym}{Sym}
\DeclareMathOperator{\Iso}{Iso}
\DeclareMathOperator{\Aut}{Aut}
\DeclareMathOperator{\Alt}{Alt}
\renewcommand{\labelenumi}{(\alph{enumi})}
\renewcommand{\theenumi}{\alph{enumi}}

\newtheorem{exer}[defi]{Exercise}

\DeclareFontFamily{OT1}{rsfs}{}
\DeclareFontShape{OT1}{rsfs}{n}{it}{<-> rsfs10}{}
\DeclareMathAlphabet{\mathscr}{OT1}{rsfs}{n}{it}

\tikzstyle{decision} = [diamond, draw, fill=gray!20,
    text width=4.5em, text badly centered, node distance=3cm, inner sep=0pt]
\tikzstyle{block} = [rectangle, draw, fill=gray!20,
  text width=5em, text centered, rounded corners, minimum height=4em]
\tikzstyle{noblock} = [rectangle, draw, 
  text width=5em, text centered, rounded corners, minimum height=4em]
\tikzstyle{rect} = [rectangle, draw, text centered, minimum height=4em]
\tikzstyle{line} = [draw, very thick, color=black!50, -latex']
\tikzstyle{cloudin} = [draw, ellipse, text width = 5.75em, node distance=3cm,
  minimum height=5em]
\tikzstyle{cloudout} = [draw, ellipse,  text width = 4em, node distance=3cm,
  minimum height=2em]

\date{January 2017}
\bbkannee{$69^{th}$ year, 2016-2017}
\bbknumero{1125}
\title{Graph Isomorphisms in quasi-polynomial time}

\begin{document}
\maketitle
\vspace{-2cm}

\centerline{\bf $\mathrm{by}$ Harald  Andr\'es HELFGOTT\footnote{corresponding author}}
\vspace{0.5cm}
\centerline{\bf Appendix B $\mathrm{by}$ Jitendra BAJPAI and Daniele DONA}

\medskip

\centerline{\bf Translated from the French original $\mathrm{by}$ Jitendra BAJPAI and Daniele DONA}

\vspace{0.75cm}
\smallskip
\noindent{{\underline {Address :}} \hspace{5.75cm}{\underline {Email :}} }\\
Universit\"at G\"ottingen \hspace{3.3cm} helfgott@math.univ-paris-diderot.fr\\
Mathematisches Institut \hspace{2.9cm}  jitendra@math.uni-goettingen.de\\
Bunsenstrasse 3-5 \hspace{4cm} daniele.dona@mathematik.uni-goettingen.de\\
D-37073 G\"ottingen\\
Germany\\

\vspace{0.5cm}
\noindent{\em Abstract:}
 Let us be given two graphs $\Gamma_1$, $\Gamma_2$ of $n$ vertices. Are they isomorphic? If they are, the set of isomorphisms from $\Gamma_1$ to
$\Gamma_2$ can be identified with a coset $H\cdot\pi$ inside the symmetric group on $n$ elements. How do we find $\pi$ and a set of generators of $H$?

The challenge of giving an always efficient algorithm answering these questions remained open for a long time. Babai has recently shown how to solve these problems -- and others linked to them -- in quasi-polynomial time, i.e. in time $\exp\left(O(\log n)^{O(1)}\right)$. His strategy is based in part on the algorithm by Luks (1980/82), who solved the case of graphs of bounded degree.

\section{Introduction}

Let $\mathbf{x}$, $\mathbf{y}$ be two strings of characters, that is,
two applications $\Omega\to \Sigma$, where $\Sigma$ ({\em alphabet}) and $\Omega$
({\em domain}) are finite sets. Every permutation group\footnote{For us, $G<S$ (or $S>G$) means ``$G$ is a subgroup of $S$, not necessarily proper.''}
$G<\Sym(\Omega)$ acts on the set $\Sigma^\Omega$ of strings of domain $\Omega$ on the alphabet $\Sigma$.
For us, {\em to describe a group $G$}, or {\em to be given a group $G$},
will  always  mean \enquote{to give, or be given, a set of generators 
 of  $G$}; {\em to describe a coset}
$H \pi$ will mean \enquote{to give an element $\pi$ of the coset
and a set of generators of $H$}. 

The {\em string isomorphism problem} consists in determining,
given $\mathbf{x}$, $\mathbf{y}$  and $G$,
whether there exists at least one element $\pi$ of $G$ mapping
$\mathbf{x}$ to $\mathbf{y}$, and in describing such elements $\pi$
({\em isomorphisms}) if they exist. It is clear that the set of isomorphisms
$\Iso_G(\mathbf{x},\mathbf{y})$ forms
a coset $\Aut_G(\mathbf{x}) \pi$ of the group
$\Aut_G(\mathbf{x})$ of automorphisms of $\mathbf{x}$ in~$G$,
i.e. of the group consisting of the elements of $G$ mapping $\mathbf{x}$ to itself. 

The challenge is to give an algorithm that solves the problem in
 polynomial time relative to the size $n = |\Omega|$ of $\Omega$,
or at least in reasonable time.
For example, the time employed could be {\em quasi-polynomial} in $n$,
meaning $\exp\left(O(\log n)^{O(1)}\right)$.
Here, as always,
$O(f(n))$ stands for a quantity bounded by $C\cdot f(n)$, for $n$
large enough and $C>0$ a constant, and $O_\epsilon$ indicates that the constant
$C$ depends on $\epsilon$.

One of the main motivations for the string isomorphism problem comes from the fact that the {\em graph isomorphism problem} reduces to it. This problem is to determine whether two finite
graphs $\Gamma_1$ and $\Gamma_2$ are isomorphic and if they are, describe the set of their isomorphisms. An {\em isomorphism} $\pi:\Gamma_1\to \Gamma_2$ is a  bijection $\pi$ from the set of vertices of $\Gamma_1$ to that of $\Gamma_2$ such that $\pi(\Gamma_1) = \Gamma_2$.
A solution would allow, for example, to find a molecule in a given database.

The graph isomorphism problem reduces in polynomial time to the string isomorphism problem, in the following way. Suppose without loss of generality that $\Gamma_1$ and $\Gamma_2$ have the same set of vertices $V$. Then, we can define $\Omega$ as the set of pairs of elements of $V$ (ordered or unordered, depending on whether our graphs are directed or not). The string $\mathbf{x}_i$, $i=1,2$, is defined as follows: for the pair $a=\{v_1,v_2\}$ (or $a=(v_1,v_2)$, if our graphs are directed), the value of $\mathbf{x}_i(a)$ is 1 if there is an edge between $v_1$ and $v_2$ in $\Gamma_1$,  and 0 otherwise. Let G be the image of the homomorphism $\iota:\Sym(V)\to \Sym(\Omega)$ defined by $\sigma^\iota(\{v_1,v_2\}) = \{\sigma(v_1),\sigma(v_2)\}$ where $\sigma^\iota = \iota(\sigma)$. Then $\iota$ induces a bijection between the set of isomorphisms of $\Gamma_1$  to $\Gamma_2$  and the set $\Iso_G(\mathbf{x}_1,\mathbf{x}_2)$.

\begin{theo}[Babai] The string isomorphism problem $\Omega\to \Sigma$ can be solved in quasi-polynomial time in the number of elements of the domain $\Omega$.
\end{theo}

In November 2015, Babai announced a solution in quasi-polynomial time, with an explicit algorithm. The preparation of this article led me to find a non-trivial error in the time analysis, but Babai was able to repair it by simplifying the algorithm. The proof is now correct. 

\begin{coro}[Babai]
The graph isomorphism problem can be solved in quasi-polynomial time in the number of vertices. 
\end{coro}
Our main reference will be~\cite{Ba}; we will also use the short version~\cite{Ba2}. We will try to examine the proof in as much detail as possible within the format of this presentation, partly to help eliminate any doubt that could remain regarding the result in the current form. The best previously known general bound for the time required by the graph isomorphism problem, due to Luks~\cite{BKL}, was $\exp(O(\sqrt{n \log n}))$.

\begin{center}
  * * *
\end{center}

\emph{Canonicity} plays a crucial role in Babai's treatment. As in category theory, or even in common usage, a choice is \emph{canonical} if it is functorial. The typical situation for us will be the following: a group $G<\Sym(\Omega)$ acts on $\Omega$ and so on $\Sigma^\Omega$. In addition, it also acts on another set $S$ and so on the maps $S\to \mathscr{C}$ where $\mathscr{C}$ is a finite set. A map $S\to \mathscr{C}$ is called a \emph{coloring}; the set $\mathscr{C}$ is called the set of \emph{colors}. A \emph{canonical} choice (with respect to $G$) for a coloring of $S$ for each string $\mathbf{x}\in \Sigma^\Omega$ is a map going from $\Sigma^\Omega$  to $\mathscr{C}^S$ commuting with the action of $G$.
 
In particular, a canonical choice can be a tool to detect non-isomorphisms: if the colorings $C(\mathbf{x})$ and $C(\mathbf{y})$ induced canonically by  $\mathbf{x}$ and $\mathbf{y}$ are not isomorphic to each other - for example, if they have a different number of vermilion elements - then $\mathbf{x}$  and $\mathbf{y}$ are not isomorphic to each other. Even when there are isomorphisms in $G$
sending $C(\mathbf{x})$ to $C(\mathbf{y})$, the set $\Iso_G(C(\mathbf{x}),C(\mathbf{y}))$ of such isomorphisms can be used to limit the set of isomorphisms $\Iso_G(\mathbf{x},\mathbf{y})$
from $\mathbf{x}$ to $\mathbf{y}$, because the latter is necessarily a subset of $\Iso_G(C(\mathbf{x}),C(\mathbf{y}))$. 

The proof also assimilates many ideas developed in previous approaches to the problem. The first step of the procedure consists in trying to follow what is essentially Luks's algorithm~\cite{Lu}. If this algorithm stops, this means that it has crashed against a quotient $H_1/H_2$ isomorphic to $\Alt(\Gamma)$ where $H_2\triangleleft H_1 < G$ and $\Gamma$ is rather large. 

Our main task is to study what happens at this moment. The principal strategy will be to try coloring $\Gamma$ in a way that depends canonically on $\mathbf{x}$.  This will limit the possible automorphisms and isomorphisms to consider. For example, if half of $\Gamma$ is colored in red and the other half in black, the group of possible automorphisms reduces to $\Sym(|\Gamma|/2)\times \Sym(|\Gamma|/2)$. A similar coloring induced by $\mathbf{y}$ limits the isomorphisms to the maps that align the two colorings. We will always find colorings to help us, except when certain structures have a very large symmetry. In that last case,
however, it is precisely that very large symmetry that
will allow a descent to considerably smaller $\Omega$. This double recursion - reduction of the group $H_1/H_2$ or descent to considerably shorter strings - will solve the problem.

\section{Foundations and earlier work}

Following the common usage for permutation groups, we will write $r^{g}$ for the element $g(r)$ to which $g\in \Sym(\Omega)$ sends $r\in \Omega$. Given a string $\mathbf{x}:\Omega\to \Sigma$ and an element $g\in \Sym(\Omega)$, we define $\mathbf{x}^g:\Omega\to \Sigma$ by $\mathbf{x}^g(r) = \mathbf{x}\left(r^{g^{-1}}\right)$. 

On the other hand, we write $\Omega^k$ for the set of $\vec{x}=(x_1,\dotsc,x_k)$ with the action on the left given by $(\phi(\vec{x}))_r = \vec{x}_{\phi(r)}$. The idea is that this is defined not only for $\phi$ a permutation, but for all maps $\phi:\{1,\dotsc,k\}\to \{1,\dotsc,k\}$, even if they are not injective. We call the elements of $\Omega^k$ {\em tuples} rather than {\em strings}.
\newpage
\subsection{Basic algorithms}

\subsubsection{Schreier-Sims}\label{subs:Schreier-Sims}
Several essential algorithms are based on an idea by Schreier~\cite{Sch}. He remarked that, for any subgroup $H$ of $G$ and any subset $A\subset G$ generating $G$ and containing a full set of coset representatives of $H$ in $G$,
\[A' = A A A^{-1}\cap H =
\left\{\sigma_1 \sigma_2 \sigma_3^{-1} : \sigma_i\in A\right\} \cap H\]
is a set of generators of $H$.

The next step is that of Sims~\cite{Si1}, \cite{Si2}, who showed the usefulness of working with a permutation group $G<\Sym(\Omega)$, $\Omega = \{x_1,\dotsc, x_n\}$, in terms of a {\em chain of stabilizers}
 
\[G =G_0> G_1 > G_2 > \dotsc >G_{n-1} = \{e\},\]
where $G_k = G_{(x_1,x_2,\dotsc,x_k)} =
\{g \in G: \forall 1\leq i\leq k\;\; x_i^g = x_i \}$ ({\em pointwise stabilizer}).

The Schreier-Sims algorithm (Algorithm 1; description based on~\cite[\S 1.2]{Lu}) builds sets  $C_i$ of representatives of $G_i/G_{i+1}$ such that $\cup_{i\leq j<n-1} C_j$ generates  $G_i$  for all $0\leq i<n-1$. The time taken by the algorithm is $O(n^5 + n^3 |A|)$, where $A$ is the set of generators of G given to us : the function \textsc{Filter} takes time $O(n)$, and every $g$ for which it is called satisfies $g \in A C \cup C A \cup C^2$, where $C$ is the value of  $\cup_i C_i$  at the end of the procedure. Of course, $|C|\leq n (n+1)/2$. 
Thanks to the algorithm itself, we can always suppose that our sets of generators are of size $O(n^2)$. The time taken by the algorithm is then $O(n^5)$.\footnote{We suppose that the initial set of generators, specifying the group $G$ of the problem, is of size $O(n^C)$, $C$ a constant. The time taken by the first use of the algorithm is then $O\left(n^{\max(5,3+C)}\right)$.}

\begin{algorithm}
  \caption{Schreier-Sims: construction of the sets $C_i$}\label{alg:schrsims}
  \begin{algorithmic}[1]
    \Function{SchreierSims}{$A$, $\vec{x}$}
    \Comment{$A$ generates $G<\Sym(\{x_1,\dotsc,x_n\})$}
    \Ensure{$\cup_{i\leq j<n-1} C_j$ generates $G_i$
    and $C_i\mapsto G_i/G_{i+1}$ is injective\; $\forall i\in \{0,1,\dotsc,n-2\}$}
    \State{$C_i\gets \{e\}$ for all $i\in \{0,1,\dotsc,n-2\}$}
    \State{$B\gets A$}
    \While{$B\ne \emptyset$}
    \State{Choose $g\in B$ arbitrary, and remove it from $B$}
    \State{$(i,\gamma) \gets \text{\textsc{Filter}($g$, $(C_i)$, $\vec{x}$)}$}
    \If{$\gamma\ne e$}
    \State{add $\gamma$ to $C_i$}
    \State{$B\gets B \cup \bigcup_{j\leq i} C_j \gamma \cup
          \bigcup_{j\geq i} \gamma C_j$}  
    \EndIf
    \EndWhile
    \State \Return $(C_i)$ 
    \EndFunction
    \Function{Filter}{$g$, $(C_i)$, $\vec{x}$}
    \Comment{returns $(i,\gamma)$ such that $\gamma\in G_i$,
      $g\in C_0 C_1\dotsb C_{i-1} \gamma$
    }
    \Require 
    $C_i\subset G_i$ and
  $C_i \to G_i/G_{i+1}$ injective\; $\forall i\in \{0,1,\dotsc,n-2\}$
  \Ensure 
  $g\notin C_0 C_1 \dotsb C_i G_{i+1}$ unless $(i,\gamma) = (n-1,e)$
  \State{$\gamma \gets g$}
  \For{$i=0$   \textbf{to} $n-2$}
  \If{$\exists h\in C_i$ such that $x_{i+1}^{h} = x_{i+1}^{\gamma}$}
  \State{$\gamma\gets h^{-1} \gamma$} 
  \Else{}
  \State{\Return $(i,\gamma)$}
    \EndIf
    \EndFor 
    \State \Return $(n-1,e)$
    \EndFunction
  \end{algorithmic}
\end{algorithm}

Once the sets $C_i$ have been built, it becomes possible to accomplish several essential tasks quickly.

\begin{exer}\label{ex:fhl} Show how to accomplish the following tasks in polynomial time, given a group $G<\Sym(\Omega)$, $|\Omega|=n$: 

\begin{enumerate}
\item Determine whether an element $g\in \Sym(\Omega)$ is in $G$. 
\item\label{it:richt} \cite{FHL} Let $H<G$ with $\lbrack G:H\rbrack\ll n^{O(1)}$. Given a test that determines in polynomial time whether an element $g\in G$  belongs to $H$, describe H. {\em Hint: work with $G > H > H_1>H_2 >\dotsc$ instead of $G =G_0> G_1 > G_2 > \dotsc$. }
\item\label{it:sshom} Given a homomorphism $\phi: G \to \Sym(\Omega')$, $|\Omega'|\ll |\Omega|^{O(1)}$, and a subgroup $H<\Sym(\Omega')$, describe $\phi^{-1}(H)$.
  \end{enumerate}
\end{exer}

Here, as always, {``describe''}  means {``find a set of generators''}  and a group is  {``given''} if such a set is given.

The Schreier-Sims algorithm describes the pointwise stabilizer $G_{(x_1,\dotsc,x_k)}$ for arbitrary $x_1,\dotsc , x_k \in \Omega$. However, we cannot ask for a set of generators of a {\em setwise stabilizer} $G_{\{x_{1},...,x_{k}\}} = \{g\in G: \{x_1^g,\dotsc ,x_k^g\} =  \{x_1,\dotsc ,x_k\} \}$: obtaining such a set for general $G$ and $x_i$ would be equivalent to solving the isomorphism problem itself (see Exercise~\ref{hidden0}).

 \subsubsection{Orbits and blocks}\label{subs:orbl}
As always, let us be given a permutation group $ G$ acting on a finite set  $\Omega$. 
 The domain $\Omega$ is the disjoint union of the {\em orbits}
 \mbox{$\{x^g: g\in G\}$} of~$G$. These orbits can be determined in polynomial time\footnote{To be precise:
   $O\left(|\Omega|^{O(1)} + |A| |\Omega|\right)$, where $|A|$ is the size of the given set of generators of $G$. We will omit all mentions of such size from now on, since, as we have already said, we can always keep it under control.}
 in $|\Omega|$. This is a simple exercise (see Exercise~\ref{hidden1}). This
 reduces to the simple task of finding the connected components of a graph. 

Suppose that the action of $G$ is transitive. (Then there is only one orbit). A {\em block} of G is a subset $B\subset \Omega$, $B\notin \{\emptyset, \Omega\}$, such that, for any $g,h\in G$, either $B^g=B^h$ or $B^g\cap B^h=\emptyset$. The collection  $\{B^g: g\in G\}$ ({\em  system of blocks}) for a given $B$ partitions $\Omega$. The action of $G$ is {\em primitive} if there is no block of size $>1$; otherwise, it is called {\em imprimitive}. A system of blocks is {\em minimal}\footnote{To paraphrase \cite[\S 1.1]{Lu}: it should be said that such a system could rather be called maximal. The size of the blocks is maximal, their number is minimal.} if the action of $G$ on it is primitive.

Let us see how to determine whether the action of $G$ is primitive, and if it is not, how to find a system of blocks of size $> 1$. By iterating the procedure, we will obtain a system of minimal blocks in polynomial time. We follow \cite{Lu} who quotes \cite{Si1}.

 For distinct $a, b\in \Omega$, let $\Gamma$ be the graph with $\Omega$ as its set of vertices and the orbit  $\{\{a,b\}^g : g\in G\}$ as its set of edges. The connected component containing $a$ and $b$ is the smallest block containing $a$ and $b$. If $\Gamma$ is connected, then the {``block''} is $\Omega$. The action of $G$ is imprimitive iff $\Gamma$ is not connected for
 an arbitrary $a$ and at least one $b$; in this case, we obtain a block containing $a$ and $b$, and thus a whole system of blocks with each block of size $>1$. One last note: if $G<\Sym(\Omega)$, we say that $G$ is {\em transitive}, or {\em primitive}, if its action on $\Omega$ is so.
 
 \subsection{Luks: the case of groups with  bounded factors}\label{subs:luks}

Luks showed how to solve the graph isomorphism problem in polynomial time in the special case of graphs of bounded degree. The {\em degree} or {\em valence} of a vertex in an undirected graph is the number of edges that contain it. He reduces this to the problem of describing the group of string automorphisms in the case of a group $G$ whose composition factors are bounded, where a composition factor is simply a quotient in a principal series (Jordan-H\"older) of $G$. The reduction process, elegant and far from being trivial, does not concern us here (see Exercise~\ref{hidden1b}). Let us rather see how Luks solves this case of the string isomorphism problem.

We follow the notation of \cite{Ba} and the ideas of \cite{Lu}.

\begin{defi}
Let  $K\subset \Sym(\Omega)$ and $\Delta\subset \Omega$ (the {``window''}). The set of {\em partial isomorphisms} $\Iso_K^\Delta$ is
 \[\Iso_K^\Delta(\mathbf{x},\mathbf{y}) = \{\tau\in K: \mathbf{x}(x) = \mathbf{y}(x^\tau)\;\;\; \forall x\in \Delta\}.\]
The set of {\em partial automorphisms} $\Aut_K^\Delta(\mathbf{x})$ is defined to be $\Iso_K^\Delta(\mathbf{x},\mathbf{x})$.
 \end{defi} 

 $\Iso_K^\Delta$ is thus the set of all the permutations  $g\in K$ sending $\mathbf{x}$ to $\mathbf{y}$ - at least from what we can see from the window $\Delta$. We generally work with $K$  of the form $H \pi$, where $H$ leaves $\Delta$ invariant (as a set). 
 
It is clear that for $K,K_1,K_2\subset \Sym(\Omega)$ and $\sigma\in \Sym(\Omega)$,
 \begin{equation}\label{eq:udu1}
 \Iso_{K \sigma}^\Delta(\mathbf{x},\mathbf{y}) =
 \Iso_K^\Delta\left(\mathbf{x},\mathbf{y}^{\sigma^{-1}}\right) \sigma,\end{equation}
 \begin{equation}\label{eq:udu2}\Iso_{K_1\cup K_2}^\Delta(\mathbf{x},\mathbf{y}) =
 \Iso_{K_1}^\Delta(\mathbf{x},\mathbf{y}) \cup
 \Iso_{K_2}^\Delta(\mathbf{x},\mathbf{y}).
 \end{equation}

It is also clear that, if $G$ is a subgroup of $\Sym(\Omega)$ and $\Delta$ is invariant under $G$, then $\Aut_G(\mathbf{x})$ is a subgroup of $G$, and, for all $\sigma\in \Sym(\Omega)$, $\Iso_{G \sigma}(\mathbf{x},\mathbf{y})$ is either empty, or a right coset of the form $\Aut_G(\mathbf{x}) \tau$, $\tau \in \Sym(\Omega)$. Let $\Delta_1, \Delta_2\subset \Omega$, $\Delta_1$ invariant under $G$. For $G' = \Aut_G(\mathbf{x})$ and $\sigma, \tau$ such that  $\Iso_{G \sigma}^{\Delta_1}(\mathbf{x},\mathbf{y}) = G' \tau$,
 \begin{equation}\label{eq:udu3}
  \Iso_{G \sigma}^{\Delta_1\cup \Delta_2}(\mathbf{x},\mathbf{y}) =
 \Iso_{G' \tau}^{\Delta_2}(\mathbf{x},\mathbf{y}) =
 \Iso_{G'}^{\Delta_2}\left(\mathbf{x},\mathbf{y}^{\tau^{-1}}\right) \tau,\end{equation}
where the second equation is an application of (\ref{eq:udu1}). Babai calls (\ref{eq:udu3})  the {\em chain rule}.

The following statement does not use the classification of finite simple groups. 

\begin{theo}[\cite{BCP}\footnote{ As a matter of fact, \cite[Thm 1.1]{BCP} is more general than stated here; for example, arbitrary abelian factors (not bounded) are allowed. This gives a generalization of Theorem \ref{thm:luxor}.}]\label{thm:bcp}
  Let $G<\Sym(\Omega)$ be a primitive group. Let $n = |\Omega|$.
  If all composition factors of $G$ are of order $\leq k$, then
  $|G|\leq n^{O_k(1)}$.
\end{theo}
Here, as usual, $O_k(1)$ designates a quantity that depends only on $k$.

\begin{theo}[Luks \cite{Lu}]\label{thm:luxor}
  Let $\Omega$ be a finite set and $\mathbf{x}, \mathbf{y}: \Omega\to
  \Sigma$ be two strings.
  Let us be given a group $G<\Sym(\Omega)$ such that
  any composition factor of $G$ is of order $\leq k$.
  It is possible to determine $\Iso_G(\mathbf{x},\mathbf{y})$ in polynomial time in $n = |\Omega|$.
\end{theo}
  \noindent{\sc Proof} ---
  {\bf Case 1: $G$ not transitive.}
  Let $\Delta_1\subsetneq \Omega$,
  $\Delta_1\ne \emptyset$ and invariant under the action of $G$.
  Define $\Delta_2 = \Omega\setminus \Delta_1$.
  Then, by (\ref{eq:udu3}), it is sufficient
  to calculate $\Iso_G^{\Delta_1}(\mathbf{x},\mathbf{y})$ (equal
  to a coset that we call $G' \tau$) and $\Iso_{G'}^{\Delta_2}(\mathbf{x},\mathbf{y}')$ for $\mathbf{y}' = \mathbf{y}^{\tau^{-1}}$.
  Now, to determine $\Iso_G^{\Delta_1}(\mathbf{x},\mathbf{y})$, we determine, recursively,
  $\Iso_G\left(\mathbf{x}|_{\Delta_1},\mathbf{y}|_{\Delta_1}\right)$, then, by Schreier-Sims, the pointwise stabilizer $G_{(\Delta_1)}$.
  In the same way, determining $\Iso_{G'}^{\Delta_2}(\mathbf{x},\mathbf{y}')$ for $\mathbf{y}' = \mathbf{y}^{\tau^{-1}}$ reduces to determining the group of isomorphisms (inside a group $G'$)
  between two chains of length $|\Delta_2|$. Since
  $|\Delta_1| + |\Delta_2| = n$ and Schreier-Sims takes time $O(n^5)$,
  everything is fine. (Accounting is left to the reader.)

  {\bf Case 2: $G$ transitive.} Let $N$ be the stabilizer of a minimal system of blocks for $G$;
  then, $G/N$ is primitive. By Theorem
  \ref{thm:bcp}, $|G/N|\leq m^{O_k(1)}$, where $m$ is the number of blocks.
  Now, for $\sigma_1,\dotsc,\sigma_{\ell}$ ($\ell = |G/N|$)
  such that $G = \cup_{1\leq i\leq \ell} N \sigma_i$,
  \begin{equation}\label{eq:rulu}
    \Iso_G(\mathbf{x},\mathbf{y}) = \Iso_{\cup_i N \sigma_i}(\mathbf{x},
  \mathbf{y}) = \bigcup_{1\leq i\leq \ell} \Iso_{N \sigma_i}(\mathbf{x},
  \mathbf{y}) = \bigcup_{1\leq i\leq \ell} \Iso_{N}(\mathbf{x},
  \mathbf{y}^{\sigma_i^{-1}}) \sigma_i
  \end{equation}
  by (\ref{eq:udu1}) and (\ref{eq:udu2}). Since the orbits of $N$ are
  contained inside the blocks, which are of size $n/m$,
  determining $\Iso_N(\mathbf{x},\mathbf{y}_i)$
  ($\mathbf{y}_i = \mathbf{y}^{\sigma_i^{-1}}$)
  reduces -- by rule (\ref{eq:udu3}) --
  to determining the group of isomorphisms of $m$ pairs
  of strings of length $n/m$.
  We have then reduced the problem to the solution of
  $\ell\cdot m = m^{O_k(1)}$
  problems for strings of length $n/m$.
  
  The last step is making the union of the cosets in (\ref{eq:rulu}).
  We have a description of each
  $\Iso_N(\mathbf{x},\mathbf{y}_i)$, either as the empty set,
  or as a right coset $H \tau_i$
  of the group $H=\Aut_N(\mathbf{x})$, of which we have
  found a description, i.e. a set of generators $A$. Then
  \[\begin{aligned}   \Iso_G(\mathbf{x},\mathbf{y}) &=
\bigcup_{1\leq i\leq \ell} \Iso_{N}(\mathbf{x},\mathbf{y}_i) \sigma_i =
\bigcup_{1\leq i\leq \ell} H \tau_i \sigma_i\\
&= \left\langle A \cup \left\{
\tau_i \sigma_i (\tau_1 \sigma_1)^{-1}: 1\leq i\leq \ell\right\}\right\rangle
  \tau_1 \sigma_1.
  \end{aligned}\]
\qed

We could have avoided a few calls to Schreier-Sims by always working with partial isomorphisms, but this has little qualitative importance.

\subsection{Relations, partitions, configurations}\label{subs:secf1}

Let  $\mathscr{C}$ ({``colors''}) be a finite set. We can assume it to be
ordered (say, from red to violet). A {\em $k$-ary relation} on a finite set $\Gamma$ is a subset $R\subset \Gamma^k$. A
{\em $k$-ary (relational) structure} is a pair
$\mathfrak{X} = (\Gamma, (R_i)_{i\in \mathscr{C}})$
where for each $i\in \mathscr{C}$, $R_i$ is a $k$-ary relation on
$\Gamma$. If the $R_i$ are all nonempty and partition $\Gamma^k$,
we say that $\mathfrak{X}$ is a {\em $k$-ary partition structure}.
In this case, we can describe $\mathfrak{X}$ by a function
$c:\Gamma^k\to \mathscr{C}$ assigning to each $\vec{x}\in \Gamma^k$ the index $i$ of the relation $R_i$ to which it belongs. We say that
$c(\vec{x})$ is the color of $\vec{x}$.

An isomorphism between two $k$-ary structures $\mathfrak{X} =
(\Gamma, (R_i)_{i\in \mathscr{C}})$ and $\mathfrak{X}' =
(\Gamma', (R_i')_{i\in \mathscr{C}})$ is a bijection $\Gamma\to \Gamma'$
 sending $R_i$ to $R_i'$ for each $i$.
It is possible to construct a functor $F_1$ sending
each $k$-ary structure $\mathfrak{X}$ on $\Gamma$
to a $k$-ary partition structure
$F_1(\mathfrak{X})$ on $\Gamma$; furthermore, $\Iso(\mathfrak{X},\mathfrak{Y})
=\Iso(F_1(\mathfrak{X}),F_1(\mathfrak{Y}))$.
The procedure is rather trivial; we give the details (Algorithm
\ref{alg:struraff}) to show what {\em index} means. This allows us to avoid using more than $\min\left(|\Gamma|^k,2^{|\mathscr{C}|}\right)$ colors, while keeping their meaning in terms of original colors
$\mathscr{C}$. 
The time needed to calculate $F_1(\mathfrak{X})$ is $O(|\mathscr{C}| |\Gamma|^{O(k)})$.
We will not take care of the implementation details for the tuple collection $\mathscr{I}$,
but it can simply be a list lexicographically ordered; 
in this case, $|\Gamma|^{O(k)}$ is $|\Gamma|^{2k}$. (In reality,
$\mathscr{I}$ would be implemented with some
{\em hash function}, which is just the art of organizing a library.)

\begin{algorithm}
  \caption{Refinement of a relational structure. Indexer.}\label{alg:struraff}
  \begin{algorithmic}[1]
    \Function{$F_1$}{$\Gamma$,$k$,$\mathscr{C}$,$(R_i)_{i\in \mathscr{C}}$}
    \State $\mathscr{I} \gets \emptyset$
    \For{$\vec{x}\in \Gamma^k$}
    \State $\mathit{a} \gets \{i\in \mathscr{C}: \vec{x}\in R_i\}$
    \State $c(\vec{x}) \gets \text{\textsc{Indexer}}(\mathscr{I},\mathit{a})$
    \EndFor
    \State \Return $(\mathscr{I}, c)$ \Comment{returns $c:\Gamma^k\to \mathscr{C}'$}\\
    \Comment{$\mathscr{C}'$ is the set of indices of $\mathscr{I}$;
      $\mathscr{I}$
      explains $\mathscr{C}'$ in terms of $\mathscr{C}$}
    \EndFunction
    \Function{Indexer}{$\mathscr{I}$,$\mathit{a}$}
    \Comment{$\mathscr{I}$ is a modifiable collection}
    \If{$\mathit{a}$ is not in $\mathscr{I}$}
    \State add $\mathit{a}$ to $\mathscr{I}$
    \EndIf
    \State \Return index of $\mathit{a}$ in $\mathscr{I}$
        \EndFunction
  \end{algorithmic}
\end{algorithm}

An element  $\vec{x}\in \Gamma^k$ defines an equivalence relation 
$\rho(\vec{x})$ on $\{1,\dotsc,k\}$: $i\sim j$ iff $x_i = x_j$.
The monoid $\mathfrak{M}(S)$ ($S$  a set) consists of the
maps $S\to S$ with composition as the binary operation.
 \begin{defi}\label{def:ukuru}
   A $k$-ary partition structure $\mathfrak{X} = (\Gamma,c)$ is called
   {\em $k$-ary configuration} if the following conditions are satisfied.
   \begin{enumerate}
   \item\label{it:klara} For all $\vec{x}, \vec{y}\in \Gamma^k$, if
     $c(\vec{x}) = c(\vec{y})$, then $\rho(\vec{x}) = \rho(\vec{y})$.
   \item\label{it:klarb} There is a homomorphism of monoids
     $\eta: \mathfrak{M}(\{1,\dotsc,k\}) \to \mathfrak{M}(\mathscr{C})$
    such that for all $\tau\in \mathfrak{M}(\{1,\dotsc,k\})$,
     $c\left(\tau(\vec{x})\right) = \tau^\eta(c(\vec{x}))$, $\forall \vec{x}\in \Gamma^k$.
     \end{enumerate}
 \end{defi}
   So, for example, for $k=2$, (\ref{it:klara}) means that the color
   of $\vec{x}=(x_1,x_2)$ {``knows''} whether $x_1=x_2$ or not, in the sense that
   if we know $c(\vec{x})$ then we know whether $x_1=x_2$ or not. Similarly, (\ref{it:klarb}) indicates that the color of $\vec{x}$ knows the colors of $(x_2,x_1)$, $(x_1,x_1)$ and $(x_2,x_2)$.

We can define a functor $F_2$ sending
each $k$-ary partition structure
 $\mathfrak{X}$ on $\Gamma$ to a $k$-ary configuration;
 as similar to $F_1$, the fact that $F_2(\mathfrak{X})$ is a refinement 
of $\mathfrak{X}$ implies that
$\Iso(\mathfrak{X},\mathfrak{Y}) = \Iso(F_2(\mathfrak{X}),F_2(\mathfrak{Y}))$.
The procedure for calculating $F_2$ is very similar to that for calculating
$F_1$ (Algorithm \ref{alg:struraff}). Instead of assigning
$\vec{x}$ the color $\{i\in \mathscr{C}: \vec{x}\in R_i\}$, we assign
it the color $\left(\rho(\vec{x}),
(c(\phi(\vec{x})))_{\phi \in \mathfrak{M}(\{1,\dotsc,k\})}\right)$.

It is easy to see that $F_2(\mathfrak{X})$ is the coarsest partition structure
$\mathfrak{X}$ that is a configuration, in the same way that $F_1(\mathfrak{X})$ is the coarsest refinement of a structure $\mathfrak{X}$ that is a partition structure. 

\begin{defi}
  Let $\mathfrak{X} = (\Gamma,c)$, $c:\Gamma^k\to \mathscr{C}$ be
  a $k$-ary partition structure.
  For $1\leq l\leq k$, we define $c^{(l)}:\Gamma^l\to
  \mathscr{C}$ as follows:
  \[c^{(l)}(\vec{x}) = c(x_1,x_2,\dotsc,x_l, x_l, \dotsc x_l).\]
  The $l$-ary partition structure
  $\mathfrak{X}^{(l)} = \left(\Gamma,c^{(l)}\right)$ is called the {\em ($l$)-skeleton} of $\mathfrak{X}$.
\end{defi}
The empty string will be of the color of viridescent seas.

\begin{exer}\label{ex:skeleconfig}
Every skeleton of a configuration is a configuration.
\end{exer}

Here, the fact that the axiom (\ref{it:klarb}) in the definition of {\em configuration}
is valid even for~$\eta$ not injective is crucial.

For $\mathfrak{X} = (\Gamma,c)$ a partition structure 
and $\Gamma'\subset \Gamma$, the {\em induced substructure}
$\mathfrak{X}\lbrack \Gamma'\rbrack$ is the structure
$(\Gamma',c|_{\Gamma'})$ defined by the restriction of $c$ to $\Gamma'$.
It is clear that, if $\mathfrak{X}$ is a configuration, $\mathfrak{X}\lbrack \Gamma'\rbrack$ is so (see Exercise~\ref{hidden1d}).

\begin{center}
  * * *
\end{center}

A partition structure should not be confused with what we will call
a {\em colored partition}.
A colored partition of a set $\Gamma$
is a coloring of $\Gamma$ to which a partition of each color class is added. (A {\em color class}
is the set of vertices of a given color.) A colored partition is said to be
{\em admissible} if each set $B$ in each partition is of size $\geq 2$. For $\alpha<1$,
a {\em colored $\alpha$-partition} is an admissible colored partition such that 
$|B|\leq \alpha |\Gamma|$ for each $B$.

A colored partition is a structure that is finer than the coloring that it refines,
but less fine than the structure that we would obtain if we gave a different color to the
elements of each part of the partition. An automorphism or
isomorphism of a colored partition must preserve its colors, but 
could permute the sets of the same size belonging to
the partition of a color. As sets of different sizes cannot evidently be permuted, it is 
clear that we can assume without loss of generality that each color
is partitioned into sets of the same size.  We add this to the definition of colored $\alpha$-partitions from now on. 

\subsection{$k$-ary coherent configurations}\label{subs:concoh}

For $\vec{x}\in \Gamma^k$, $z\in \Gamma$ and $1\leq i\leq k$, we
define $\vec{x}^i(z) \in \Gamma^k$ as follows:
\[\vec{x}^i(z) = (x_1, x_2, \dotsc,x_{i-1},z,x_{i+1},\dotsc, x_k).\]
\begin{defi}\label{def:coh}
  A {\em $k$-ary coherent configuration} $\mathfrak{X}=(\Gamma,c)$
is a $k$-ary configuration having the following property: there is a
  function $\gamma:\mathscr{C}^k\times \mathscr{C} \to \mathbb{Z}_{\geq 0}$ such 
  that, for arbitrary $\vec{k}\in \mathscr{C}^k$ and $j\in \mathscr{C}$ and 
  any $\vec{x}\in \Gamma^k$ such that $c(\vec{x})=j$,
  \[|\{z\in \Gamma: c(\vec{x}^i(z)) = k_i\; \forall 1\leq i \leq k\}| = \gamma(\vec{k},j).\]
  The values $\gamma(\vec{k},j)$ are called {\em intersection numbers}
  of $\mathfrak{X}$.
\end{defi}
A coherent configuration is said to be {\em classical} if $k=2$.

\begin{rema}
  The classical coherent configurations were introduced 
  by Higman \cite{Hi}. The first examples were
  of {\em Schurian} type: a configuration is
  {\em Schurian} if it is the partition
  of $\Gamma^2$ into its orbits ({``orbitals''}) under the action of a group $G<\Sym(\Gamma)$.
\end{rema}

\begin{defi}
  If a classical coherent configuration has only two colors, one for
  $\{(x,x) : x\in \Gamma\}$ and the other one for its complement, then the configuration
  is said to be a {\em clique} or {\em trivial}.
\end{defi}

\begin{exer}\label{ex:skelecoh}
Any skeleton of a coherent configuration is coherent.
\end{exer}

Once again the axiom (\ref{it:klarb}) of configurations plays a key role.

\begin{exer}\label{ex:rescoh}
  Let $\mathfrak{X} = (\Gamma,c)$ be a coherent configuration and
  $\Gamma'\subset \Gamma$ a color class with respect to the coloring 
  induced by $c$ on $\Gamma$. Then the induced substructure   
   $\mathfrak{X}\lbrack \Gamma'\rbrack$ is a coherent configuration.
\end{exer}

Here, we must use a special case of (\ref{it:klarb}): the color 
$c(x_1,\dotsc,x_n)$ {``knows''} the colors $c(x_1),\dotsc,c(x_n)$,
because $c(x_i)=c(x_i,\dotsc,x_i)$.

Let $0\leq l<k$ and $\vec{x} \in \Gamma^l$. We color
$\Gamma^{k-l}$ as follows: for $\vec{y}\in \Gamma^{k-l}$,
\[c_{\vec{x}}(\vec{y}) = c(\vec{x} \vec{y}).\]
 The result is a $(k-l)$-ary partition structure
  $\mathfrak{X}_{\vec{x}} = \left(\Gamma,c_{\vec{x}}\right)$.

\begin{exer}\label{ex:veccoh}
    Let $\mathfrak{X} = (\Gamma,c)$ be a configuration;
    let $\vec{x}\in \Gamma^l$, $0\leq l <k$. Then
    \begin{enumerate}
      \item
    $c_{\vec{x}}$ is a refinement of the coloring of the skeleton $\mathfrak{X}^{(k-l)}$. \item If $\mathfrak{X}$ is coherent, $\mathfrak{X}_{\vec{x}}$
        is so.
    \end{enumerate}
  \end{exer}
  
It is clear that $\mathfrak{X}_{\vec{x}}$ is
{\em canonical with respect to $\vec{x}$}, meaning that
  $\mathfrak{X}\to \mathfrak{X}_{\vec{x}}$
  commutes with the action on $\Gamma$ of the stabilizer in $\Sym(\Gamma)$
  of the points $x_1,\dotsc,x_l$.

  \begin{defi}\label{def:uniprim}
    A coherent configuration $(\Gamma,c)$
    is said to be {\em homogeneous} if the color
    $c(x,x,\dotsc,x)$ of any vertex $x\in \Gamma$ is the same.
    A classical coherent configuration
    is called {\em primitive} if it is homogeneous
 and the graphs $\mathscr{G}_r = \{(x,y): x,y\in \Gamma, c(x,y)=r\}$
    (for any color $r$ such that $c(x,y)= r$ for at least one pair
    $(x,y)$ with $x\ne y$) are all connected. It is said to be {\em uniprimitive} if it is primitive and non-trivial.
  \end{defi}
We do not need to specify whether these graphs are strongly connected
  (that is, there is a path from any vertex to the other,
  respecting the orientation) or weakly connected (without considering the orientation):
  the fact that $(\Gamma,c)$ is coherent, classical and homogeneous implies that 
  $d^+_r(x)=|\{y\in \Gamma: (x,y)\in \mathscr{G}_r\}|$
  is independent from $x$ (why?),
  and that implies in turn that any weakly connected component of
  $\mathscr{G}_r$ is connected (see Exercise~\ref{hidden2}).

  \begin{exer}\label{ex:noclique}
    Let $\mathfrak{X}=(\Gamma,c)$ be a classical coherent configuration that is not
    a clique.
    There is no set $B\subset \Gamma$, $|B|>|\Gamma|/2$, such that
    the restriction of $\mathfrak{X}$ to $B$ is a clique.
  \end{exer}
  \noindent{\sc Solution} --
  If we suppose the edges of the large clique to be white,
  let {\em black} be another edge color of
  $\mathfrak{X}$, and let $\mathscr{G} = \mathscr{G}_{\text{black}}$.
  Now, for
  a biregular\footnote{See the definition in \S \ref{subs:hypdes}.} nonempty directed graph $\mathscr{G}$ with $\Gamma$ as set of vertices, it is impossible that there be a set $B\subset \Gamma$, $|B|>|\Gamma|/2$, such that the reduction of the graph to $B$ is empty (why?).
  \qed
  
  \begin{exer}\label{ex:samesize}
    Let $(\Gamma,c)$ be a homogeneous classical coherent configuration.
      \begin{enumerate}
      \item\label{it:bibr1} Let $r_0,\dotsc,r_k$ be a sequence of colors.
        Then, if $x_0,x_k\in \Gamma$ are such that $c\left(x_0,x_k\right)=r_0$,
        the number of $x_1,\dotsc,x_{k-1}\in \Gamma$ such that
        $c\left(x_{i-1},x_i\right) = r_i$ for every $1\leq i\leq k$ depends
       only on $r_0,\dotsc,r_k$.
      \item\label{it:bibr2}
        For every color $r$, every connected component of
    $\mathscr{G}_r$ is of the same size.
    \end{enumerate}
  \end{exer}

\noindent{\sc Solution (sketch)} --- In (\ref{it:bibr1}), the case
  $k=2$ holds by the definition of {``coherent''}; prove the cases $k>2$
  by induction. To prove (\ref{it:bibr2}), use (\ref{it:bibr1}).

  \subsection{The $k$-ary canonical refinement using Weisfeiler-Leman}\label{subs:wl}

Let us define a functor $F_3$ sending a configuration
  $\mathfrak{X}=(\Gamma,c)$ to a coherent configuration
  $F_3(\mathfrak{X})= (\Gamma,c')$. As $F_3(\mathfrak{X})$ will be
  a canonical refinement of $\mathfrak{X}$, we will have
$\Iso(\mathfrak{X},\mathfrak{Y}) = \Iso(F_3(\mathfrak{X}),F_3(\mathfrak{Y}))$.

Algorithm \ref{alg:weislem}, which computes $F_3$, is based on an idea 
 by Weisfeiler and Leman\footnote{Also called Lehman, but \cite{Ba} indicates that the second author 
  preferred {\em Leman}. Two natural transformations $L\to \El$, $\El\to L$ do not have to be inverses of each other.}
  \cite{WL}.
  It is a matter of iterating a refinement procedure. If, within an iteration,
  no refinement takes place -- that is, if the equivalence classes of the 
new coloring $\mathscr{C}_i$ are the same as those of the old coloring
$\mathscr{C}_{i-1}$ -- then, (a) no refinement
  will be produced in the future, (b) the coloring $\mathscr{C}_{i-1}$ is
  already coherent.

  If the coloring $\mathscr{C}=\mathscr{C}_0$ at the beginning
  has $r$ different colors, it is clear that it can be refined at most 
  $|\Gamma|^k - r$ times. Then, $|\Gamma|^k-r$ iterations are sufficient 
  to produce a coherent configuration.
  In particular, if indexing takes logarithmic time, and the vector
  in step \ref{sta:thisline} of Algorithm \ref{alg:weislem} is stored as a sparse vector (having at most
  $|\Gamma|$ non-zero entries),
  the time taken is $O(k^2 |\Gamma|^{2k+1} \log |\Gamma|)$.
  (A stronger claim is made in \cite[\S 2.8.3]{Ba}.)

Algorithms of Weisfeiler-Leman type were previously seen as
  a plausible approach to the graph isomorphism problem. Since
  \cite{CFI}, \cite{EvP}, it is clear that they are not sufficient on their own. 
  They are an invaluable tool all the same. The $k$-ary version here is due to Babai-Mathon \cite{Ba3} and Immerman-Lander
  \cite{ImL}.
  
\begin{algorithm}
  \caption{Weisfeiler-Leman for $k$-ary configurations.}\label{alg:weislem}
  \begin{algorithmic}[1]
    \Function{WeisfeilerLeman}{$\Gamma$,\; $k$,\; $c:\Gamma^k\to \mathscr{C}$}
    \State $\mathscr{C}_0\gets \mathscr{C}$;
    $c_0\gets c$; $i_0\gets |\Gamma^k|-|c(\Gamma^k)|$
    \For{$i=1$   \textbf{to} $i_0$}
    \State $\mathscr{I}_i\gets \emptyset$
    \For $\vec{x} \in \Gamma^k$
    \State $\nu\gets
    \left(c_{i-1}(\vec{x}),
    \left(\left|\left\{z\in \Gamma: c_{i-1}(\vec{x}^{j}(z)) = r_j\;\;\;
    \forall\; 1\leq j\leq k\right\}\right|\right)_{\vec{r}\in \mathscr{C}_{i-1}^k}
    \right)$ \label{sta:thisline}
    \State $c_i(\vec{x}) = \text{\textsc{Indexer}}\left(\mathscr{I}_i,
    \nu\right)$\Comment{\textsc{Indexer} is
      like in Algorithm \ref{alg:struraff}}
    \EndFor
    \State $\mathscr{C}_i\gets \text{indices of $\mathscr{I}_{i}$}$    
    \EndFor
    \State \Return $\left(c_{n-i_0}:\Gamma^k\to \mathscr{C}_{n-i_0},
    (\mathscr{I}_i)_{1\leq i
      \leq n-i_0}\right)$ \Comment{$(\mathscr{I}_i)$ gives a meaning to $\mathscr{C}_{n-r}$}
    \EndFunction
  
  \end{algorithmic}
\end{algorithm}

\subsection{Graphs, hypergraphs and block designs}\label{subs:hypdes}
We already know that a {\em graph} is a pair $(V,A)$, where $V$
is a set ({``vertices''}) and $A$ is a collection of
pairs of elements of $V$ (or of subsets of $V$ with two elements,
if the graph is undirected). An undirected graph is said to be {\em regular}
if the degree of each vertex is the same; a directed graph
 is said to be {\em biregular} if the {\em outdegree}
$d^+(v)=|\{w\in V: (v,w)\in A\}|$ and the {\em indegree}
$d^-(v)=|\{w\in V: (w,v)\in A\}|$ are
independent from $v$. (For $V$ finite, they are necessarily the same constant.)

A {\em bipartite graph} is a triplet
$(V_1,V_2;A)$ with $A\subset V_1\times V_2$. A bipartite graph
is {\em semiregular} if the degree\footnote{
We omit the prefixes
{``in-''} and {``out-''}, because it is evident that it is the indegree 
in the case of $v_2$ and the outdegree in the case of $v_1$.}
$d^+(v_1)$ is independent
from $v_1\in V_1$ and the degree $d^-(v_2)$ is independent from $v_2\in V_2$. 

\begin{exer}\label{ex:biparhom}
      Let $\mathfrak{X}=(\Gamma,c)$ be a classical coherent configuration.
  \begin{enumerate}
  \item\label{it:bar1}
    Let $C_1$, $C_2$ be two color classes and let {\em green} be
an edge color in $C_1\times C_2$. Then,
  the bipartite graph $(C_1,C_2;\mathscr{G}_{\text{green}})$ is semiregular.
\item\label{it:bara2}
  Let $y\in \Gamma$, and $L_i(y) = \{x\in \Gamma: c(x,y) = i\}$. Let
  {\em aqua}, {\em beige} and {\em cyan} be three edge colors. Then,
  for $L_1 = L_{\text{aqua}}(y)$ and $L_2 = L_{\text{beige}}(y)$,
  the bipartite graph $(L_1,L_2; R_{\text{cyan}} \cap (L_1\times L_2))$
  is semiregular.
\end{enumerate}
\end{exer}

\begin{exer}\label{ex:biparcomp}
      Let $\mathfrak{X}$ be a classical coherent configuration.
      Let $C_1$, $C_2$ be two color classes. Let {\em green} be an edge color in $C_1\times C_2$ and 
    {\em red} be an edge color in $C_2\times C_2$. Let
    $B_1,\dotsc,B_m$ be the connected components of $\mathscr{G}_{\text{red}}$
    in $C_2$. Define the bipartite graph $Y = (C_1,\{1,\dotsc,m\}; D)$
    as follows: $(x,i)\in D$ if and only if $(x,y)\in \mathscr{G}_{\text{green}}$ for at least one
    $y\in B_i$. Then $Y$ is semiregular.
\end{exer}

\noindent{\sc Solution} --- Notice that, for $y\in B_i$ and $x\in C_1$, $(x,i)\in D$ if and only if there exist $y_{0},y_{1},\dotsc,y_{k}=y$ such that $(x,y_{0})$ is green and $(y_{j},y_{j+1})$ is red for $0\leq j<k$. Conclude by Exercises
\ref{ex:samesize}\ref{it:bibr1} and \ref{ex:biparhom}\ref{it:bar1} that all vertices in $\{1,\dotsc,m\}$ have
the same degree in $Y$.

In an analogous way, show that, for $x\in C_1$ and $y\in B_i$
such that $(x,y)$ is green, the number of $z\in B_i$ such that $(x,z)$
is green does not depend on $x$, $y$ or $i$. Call this number~$q$.
Then, the degree of any $x\in C_1$ is the number of green edges $(x,y)$ in ${\mathfrak X}$, divided by $q$.
By Exercise \ref{ex:biparhom}\ref{it:bar1}, conclude that it does
not depend on $x$. \qed

A bipartite graph is {\em complete} (as a bipartite graph) if
$A = V_1\times V_2$. A bipartite graph that is neither empty nor complete is 
called {\em non-trivial}.

A {\em hypergraph} $\mathscr{H} = (V,\mathscr{A})$ consists of a set
$V$ ({``vertices''}) and a collection $\mathscr{A}$ of subsets of $V$
({``edges''}), possibly with repeated subsets. A hypergraph is said to be {\em $u$-uniform} if
$|A|=u$ for every $A\in \mathscr{A}$. It is said to be {\em regular of degree~$r$}
if every $v\in V$ belongs exactly to $r$ sets $A$ in $\mathscr{A}$.

The complete $u$-uniform hypergraph on $V$ is $(V,\{A\subset V: |A|=u\})$,
where each set $A$ is counted once. An {\em edge coloring} of the complete hypergraph is a map
from $\{A\subset V: |A|=u\}$ to a finite set $\mathscr{C}$.

A {\em balanced block design} (BBD) of parameters $(v,u,\lambda)$
is a hypergraph with $|V|=v$ vertices, 
$u$-uniform and regular of degree $r\geq 1$, such that every pair 
$\{v_1,v_2\}$ of distinct vertices is contained in 
exactly $\lambda\geq 1$ edges ({``blocks''}).
A {\em degenerate} block design has the same definition, but with $\lambda=0$,
and the additional condition of being a regular hypergraph. (The regularity
can be deduced from the definition if $\lambda\geq 1$.) A block design
is {\em incomplete} if $u<v$.
Let us call $b$ the number $|\mathscr{A}|$ of edges of a BBD.
\begin{prop}[Fisher's inequality\footnote{Yes, R. A. Fisher, the
      statistician.
Here {\em design} comes from {\em experimental design}.} \cite{F}]\label{prop:fishy}
  For every balanced incomplete block design, $b\geq v$.
\end{prop}
It is easy to see that this inequality is true even for degenerate designs. 

Block designs admit a generalization. A
{\em $t$-$(v,u,\lambda)$ design} is a hypergraph $(V,\mathscr{A})$
$u$-uniform with $v=|V|$ vertices such that every $T\subset V$ of size $t$
is contained in exactly $\lambda$ edges. Here $t\geq 2$ and
$\lambda\geq 1$.
We always write
$b=|\mathscr{A}|$.
\begin{prop}[\cite{RChW}]\label{prop:RW}
  For every $t$-$(v,u,\lambda)$ design and every $s\leq \min(t/2,v-u)$, we
  have $b\geq \binom{v}{s}$.
\end{prop}

\subsection{Johnson schemes}\label{sec:schjoh}

An {\em association scheme} is a classical coherent configuration
$(\Gamma,c:\Gamma^2\to \mathscr{C})$ such that
$c(x,y)=c(y,x)\;\; \forall x,y\in \Gamma$.
(Thus, the meaning of the word  {\em scheme} has nothing to do with the {\em schemes} of algebraic geometry.)

Let $s\geq 2$ and $r\geq 2 s + 1$. A {\em Johnson scheme}
$\mathscr{J}(r,s) = (\Gamma,c)$ is given by
\[
  \Gamma = \mathscr{S}_s(\Lambda) =
  \{S\subset \Lambda: |S|=s\},\;\;\;\;\;\;\;\;\;\;\;\;
  c(S_1,S_2) = |S_1\setminus (S_1\cap S_2)|,\]
  where $\Lambda$ is a set of $r$ elements. The relation $R_i$
  is of course the set \[R_i = \{(S_1,S_2): c(S_1,S_2) =i\}.\]
  
  Notice that we have implicitly defined a functor from the category
  of sets $\Lambda$ with $|\Lambda|=r$ to the category of 
  Johnson schemes. This is a full functor; that is, the only
  automorphisms of $\mathscr{J}(r,s)$ are those induced by
    $\Sym(\Lambda)$.

\subsection{Identification of groups and schemes}\label{subs:idgroup}

It is one thing to prove that two groups $G$, $H$ are isomorphic,
and another to explicitly build an isomorphism $\phi$ between them.
 This last task involves, at least, giving the images $\phi(g_1),\dotsc,\phi(g_r)$
of generators $g_1,\dotsc,g_r$ of $G$.

Let us see a specific case that will be crucial to us. We will have a permutation group
$G<\Sym(\Gamma)$, and we will know that it is isomorphic to the
abstract group $\Alt_m$. How to build an isomorphism?

If $m$ is not too small with respect to $n = |\Gamma|$, it is known
that $G$ must be isomorphic to a permutation group of the form
$\Alt_m^{(k)}$, which is none other than the group $\Alt_m$ acting
on the set $\mathscr{S}_k(\Lambda_0) =
\{S\subset \Lambda_0: |S|=k\}$ of $\binom{m}{k}$
elements, where $\Lambda_0$ is a set of $m$ elements.\footnote{
  Babai calls the groups $\Alt_m^{(k)}$ {\em Johnson groups}, by
analogy with Johnson schemes. Since $\Alt_m^{(k)}$ is just $\Alt_m$ in disguise, should we not call the latter a {\em Ramerrez group}?} In other terms, 
there exists a bijection $\iota_0:\Gamma\to \mathscr{S}_k(\Lambda_0)$
and an isomorphism $\phi_0:G\to \Alt(\Lambda_0)$ such that
\[\iota_0\left(\omega^g\right) = \iota_0(\omega)^{\phi_0(g)}.\]
The problem consists in building $\iota:\Gamma\to \mathscr{S}_k(\Lambda)$
and $\phi:G\to \Alt(\Lambda)$,
computable in polynomial time, with these properties.

We follow \cite{BLS}. Let $\Upsilon\subset \Gamma\times \Gamma$ be
the smallest orbital of $G$ (outside the diagonal
$(\{\omega,\omega\}: \omega\in \Gamma\}$); let $\Delta\subset
\Gamma\times \Gamma$ be the largest orbital. We will suppose
that $m>(k+1)^{2}-2$, which means that $n$ is not too big with respect to
$m$.\footnote{If $m\leq (k+1)^2-2$, then $n$ is so large that $m! = n^{O(\log n)}$. In that case, we will be able to
  afford to remove $G$ (that is, in the application
  of interest for us, a quotient $G/N$) -- brutally, as in case
  2 of the proof of Theorem \ref{thm:luxor} (Luks). We can also simply do without the assumption $m>(k+1)^2-2$ at the cost of some complications in what follows. In particular,
  $\phi(\Delta)$ would not equal $R_k$ as in (\ref{eq:ludod}), but some other
$R_j$.}
Then,
\begin{equation}\label{eq:ludod}\begin{aligned}
\phi(\Upsilon) &= R_1 = \{(S_1,S_2)\in \mathscr{S}_k(\Lambda_0): |S_1\cap S_2|=k-1\},\\
\phi(\Delta) &= R_k = \{(S_1,S_2)\in \mathscr{S}_k(\Lambda_0): S_1\cap S_2=\emptyset\}.
\end{aligned}\end{equation}

Let us define, for $(x,y)\in \Upsilon$,
\[
B(x,y) = \{z\in \Gamma: (x,z)\notin \Delta, (y,z)\in \Delta\}.\]
This is the set of all $z$ such that $\iota_0(z)$ intersects
$\iota_0(x)$ but not $\iota_0(y)$. Let
\[C(x,y) = \Gamma \setminus  \bigcup_{z\in B(x,y)} \{r: (z,r) \in \Delta\}.\]
Then \[\begin{aligned}
\iota_0&(C(x,y))\\ &= \{S\in \mathscr{S}_k(\Lambda_0): 
S\cap S' \ne \emptyset\;\; \forall S'\in \mathscr{S}_k(\Lambda_0)
\; \text{s.t. $S'\cap \iota_0(x)\ne \emptyset$, $S'\cap \iota_0(y) = \emptyset$}
\}\\ &= \{S\in \mathscr{S}_k(\Lambda_0): i\in S\},\end{aligned}\]
where $i$ is the element of $\iota_0(x)$ that is not in $\iota_0(y)$.

Let $\Lambda$ be the collection $\{C(x,y): (x,y)\in \Upsilon\}$, without
multiplicities. We can compute and compare $C(x,y)$ for given $(x,y)$,
and compute and index the set $\Lambda$, all in polynomial time.
We compute, also in polynomial time, the action of $G$ on
$\Lambda$ induced by the action of $G$ on $\Upsilon$. This defines
$\phi:G\to \Alt(\Lambda)$.

There is a natural bijection 
$j:\Lambda\to \Lambda_0$ commuting with the action of $G$:
it sends $C(x,y)$ to $i$, where $i$ is
the element of $\Lambda_0$ such that $\iota_0(C(x,y)) =
\{S\in \mathscr{S}_k(\Lambda_0): i\in S\}$.
It is clear that, for $\omega\in \Gamma$, $\omega\in C(x,y)$ if and only if
$j(C(x,y))\in \iota_0(\omega)$.
Thus, we obtain the bijection
$\iota:\Gamma\to \mathscr{S}_k(\Lambda)$, given by
\[\iota(\omega) = \{\gamma\in \Lambda: \omega \in \gamma\}.\]
This one satisfies $\iota\left(\omega^g\right) = \iota(\omega)^{\phi(g)}$.

The applications $\phi$, $\iota$ are thus the ones we wanted; we
have built an explicit isomorphism between $G$ and $\Alt(\Lambda)$.
Notice that this very procedure allows us to build an explicit isomorphism 
between, on one side,
an association scheme (\S \ref{sec:schjoh}) that we know to be
isomorphic to a Johnson scheme $\mathscr{J}(m,k)$, and, on the other side, that Johnson scheme.
 
\newpage
\section{The main procedure}\label{sec:procprin}
\begin{center}
\begin{tikzpicture}[scale=2, node distance = 2cm, auto]
  \node [rect] (init) {Function String-Isomorphism};
    \node [cloudin, left of=init, node distance = 5.5cm] (input) {
    Input:\\ $G<\Sym(\Omega)$\\ $\mathbf{x},\mathbf{y}:\Omega\to \Sigma$};
    \node [cloudout, right of=init, node distance = 5cm] (output) {Output: $\Iso_G(\mathbf{x},\mathbf{y})$};
    \node [decision, below of=init] (trans) {$G$ transitive?};
    \node [decision, below of=trans] (idtop) {$G/N \sim \Alt_m$?};
    \node [block, right of=idtop, node distance = 4.1cm] (luks) {\bf Recursion
      $n'\leq n/2$};
        \node [block, right of=trans, node distance = 4.1cm] (petrec) {\bf Recursion
      $n'<n$};
    \node [block, left of = trans, node distance = 3.5cm] (align) {Align};
    \node [decision, below of=idtop, node distance = 3cm] (gampetit) {$m$ small?};
    \node [block, left of=gampetit, node distance = 3.5cm] (struc)
          {$\text{Blocks} \sim \binom{\Gamma}{k}$};
    \node [decision, below of=struc] (primi)
          {$G$ primitive?};
          \node [decision, right of = primi, node distance = 3.5cm] (k1) {$k=1$?};
          \node [block, right of =k1, node distance = 2.95cm] (triv) {\bf Trivial case};
          \node [decision, below of = k1] (sym) {Symmetry $>1/2$?};
                    \node [block, right of =sym, node distance = 4.75cm] (pull) {Pullback};
                    \node [block, below of = sym, node distance = 3cm] (xrel) {$\mathbf{x},
                    \mathbf{y}\to$ $k$-ary relations on $\Gamma$};
                    \node [block, left of =xrel, node distance = 4cm] (wl) {$k$-ary Weisfeiler - Leman};
                    \node [block, below of = primi, node distance = 3.25cm] (cert) {Local certificates};
  \node [decision, left of = cert, node distance = 5cm] (splirel) {Split or relations?};
                    \node [decision , below of = splirel, node distance = 4.2cm] (certsym) {Fullness $>1/2$?};
                    \node [block, left of = idtop, node distance = 3.5cm] (redsqrt) {\bf Reduction of $G/N$ to $\Alt_{m'}$\\ $m'\ll \sqrt{m}$};
                    \node [decision, left of = redsqrt, node distance = 3cm] (spj) {Split or {\textsc Johnson}?};
                    \node [block, above of = spj, node distance = 3cm] (redhalf)
                          {\bf Reduction of $G/N$ to $\Alt_{m'}$\\ $m'\leq |m|/2$};
                          \node [decision, below of = spj, node distance = 4.2cm] (coldom) {Does a color dominate?};
                          \node [block, below of = coldom, node distance = 3.5cm] (design) {Design Lemma};

    \path[line] (input) -- (init);
    \path[line] (init) -- (output);      
    \path [line] (init) -- (trans);
        \path [line] (trans) -- node [near start, color=black] {no} (petrec);
    \path [line] (trans) -- node [,color=black] {yes} (idtop);
    \path [line] (idtop) -- node [, color=black] {no} (luks);
    \path [line] (idtop) -- node [, color=black] {yes:
      $G/N \sim \Alt(\Gamma)$} (gampetit);
    \path [line] (gampetit) -| node [near start, color=black] {yes} (luks);
    \path [line] (gampetit) -- node [near start, color=black] {no}(struc);
    \path [line] (struc) -- (primi);
    \path [line] (align) -- (trans);
    \path [line] (primi) -- node [, color=black] {yes} (k1);
    \path [line] (k1) -- node [, color=black] {yes} (triv);
    \path [line] (k1) -- node [, color = black] {no} (sym);
    \path [line] (sym) -- node [, color = black] {yes} (pull);
    \path [line] (pull) -- (luks);
    \path [line] (sym) -- node [near start, color = black] {no} (xrel);
    \path [line] (xrel) -- (wl);
    \path [line] (primi) -- node [, color = black] {no} (cert);
    \path [line] (cert) -- (certsym);
    \path [line] (certsym) -- node [near start, color=black] {no} (splirel);
    \path [line] (splirel) -- node [near start, color = black] {
      rels.} (wl);
    \path [line] (certsym) -| node [very near start, color = black] {
      yes} (pull);
    \path [line] (wl) -- (design);
    \path [line] (design) -- (coldom);
    \path [line] (coldom) -| node [near start, color=black] {no\;\;\;\;\;\;\;\;\;\;\;\;\;\;\;\;\;\;\;\;\;\;\;\;} (luks);
    \path [line] (coldom) -- node [near start, color=black] {yes} (spj);
    \path [line] (redhalf) -- (align);
    \path [line] (redsqrt) -- (align);
    \path [line] (spj) -- node [near start, color=black] {split} (redhalf);
    \path [line] (spj) -- node [, color= black] {$\mathbb J$} (redsqrt);
    \path [line] (splirel) |- node [near start, color=black]  {split} (redhalf);
\end{tikzpicture}
\end{center}
\subsection{First step: recursion in the style of Luks}

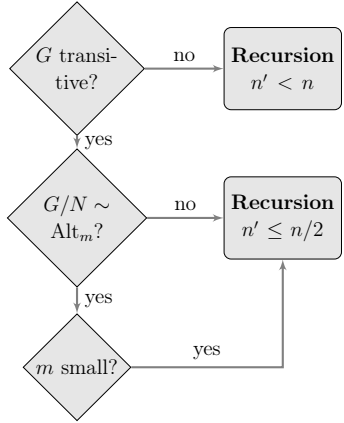
\begin{center}
  \scalebox{0.66}{
\begin{tikzpicture}[scale=2, node distance = 2cm, auto]
    \node [decision, below of=init] (trans) {$G$ transitive?};
    \node [decision, below of=trans, node distance = 3cm] (idtop) {$G/N \sim \Alt_m$?};
           \node [block, right of=trans, node distance = 4.1cm] (petrec) {\bf Recursion
      $n'<n$};
    \node [block, right of=idtop, node distance = 4.1cm] (luks) {\bf Recursion
      $n'\leq n/2$};
    \node [decision, below of=idtop, node distance = 3cm] (gampetit) {$m$ small?};
    \path [line] (trans) -- node [, color=black] {no} (petrec);
    \path [line] (trans) -- node [,color=black] {yes} (idtop);
        \path [line] (idtop) -- node [,color=black] {yes} (gampetit);
    \path [line] (idtop) -- node [, color=black] {no} (luks);
    \path [line] (gampetit) -| node [near start, color=black] {yes} (luks);
\end{tikzpicture}
}
\end{center}

The first steps of the procedure are those of the proof of Theorem 
\ref{thm:luxor} (Luks). 
In particular, if $G<\Sym(\Omega)$ is not transitive, we proceed exactly
as in the intransitive case of the proof of Theorem \ref{thm:luxor}.
Although it is possible that $n=|\Omega|$ decreases only slightly,
the recursion works, because its cost is also slight: we 
just have to subdivide the problem among the orbits of $G$.

Suppose that $G$ is transitive. We know that we can quickly find
a minimal system of blocks $R = \{B_i: 1\leq i\leq r\}$,
$B_i\subset \Omega$ (\S \ref{subs:orbl}). By Schreier-Sims, we also find, in polynomial time, 
the subgroup $N\triangleleft G$ of elements
$g\in G$ such that $B_i^g = B_i$ for every $i$. The group $H = G/N$ acts on $R$.

Instead of Theorem \ref{thm:bcp}
\cite{BCP}, we will use a consequence of the Classification
of Finite Simple Groups (CFSG). It was first derived by
Cameron, then refined by Mar\'oti. 

\begin{theo}[\cite{Cam}, \cite{Ma}]\label{thm:cam}
  Let $H< \Sym(R)$ be a primitive group, where $|R| = r$ is larger
  than an absolute constant. Then, either\footnote{For us, $\log_2$
    designates the logarithm in base $2$, and not the iterated logarithm
    $\log \log$.} 
  \begin{enumerate}
   \item\label{it:rila} $|H| < r^{1 + \log_2 r}$, or
   \item\label{it:rilb}
    there is a $M\triangleleft H$
    such that $R$ is divided
    \footnote{The statement in \cite{Cam}, \cite{Ma} is stronger: it describes the entire action of $H$ on $R$. To tell the truth, $M$ is isomorphic, as a permutation group, to $(\Alt_m^{(k)})^s$, $s\geq 1$. We have $r =\binom{m}{k}^s$.}
    into a system of $\binom{m}{k}$ blocks on which $M$ acts as a group
    $\Alt_m^{(k)}$, $m\geq 5$. Furthermore,
    $\lbrack H: M\rbrack \leq r$.
  \end{enumerate}
\end{theo}
The bound $\lbrack H:M\rbrack\leq r$ is deduced from $m>2$, $|H| \geq r^{1+\log_2 r}$, $|H|\leq m!^s s!$, $m^s\leq r$ and
      $\lbrack H:M\rbrack\leq 2^s s!$, where $s\geq 1$ is a parameter in Cameron-Maroti.
      
It is possible \cite{BLS} to find in polynomial time 
the normal subgroup $M$ and the blocks of the action of $M$.
We have already seen in \S \ref{subs:idgroup} how to explicitly identify
the action of $M$ with that of $\Alt_m^{(k)}$.

Besides, the Schreier-Sims algorithm allows us to calculate
$|H|$ in polynomial time, and thus also tells us whether we are in case
 (\ref{it:rila}). If this is the case, we proceed
as in the transitive case of the proof of Theorem \ref{thm:luxor}.
We thus reduce the problem to $r^{2+\log_2 r}$
instances  of the problem for strings of length $\leq n/r$.

If we are in case (\ref{it:rilb}) we always start by
reducing the problem to $\lbrack H:M\rbrack$ instances of the problem with
$M$ instead of $H$: by equation (\ref{eq:udu2})
and as in equation (\ref{eq:rulu}),
\[\Iso_H(\mathbf{x},\mathbf{y}) = \bigcup_{\sigma \in S} \Iso_M\left(
\mathbf{x}, \mathbf{y}^{\sigma^{-1}}\right) \sigma,\]
where $S$ is a system of representatives of the cosets of $M$ in $H$.

If $m\leq C \log n$, where $C$ is a constant,
\[|M| = \frac{m!}{2} < m^m \leq m^{C \log n}\leq (m')^{C \log n},\]
where $m' = \binom{m}{k}$.
Thus, here as in case (\ref{it:rila}), we can
proceed as in the transitive case of the proof of Theorem \ref{thm:luxor}.
We obtain a reduction to $\leq r\cdot (m')^{C \log n} = (m')^{O(\log n)}$
instances of the problem for strings of length $n/m'$.
This is completely consistent with the goal of having a solution in 
quasi-polynomial time in $n$ (or even in time $n^{O(\log n)}$).

It is left to know what to do in the following case:
there is an isomorphism $\phi:G/N \to \Alt(\Gamma)$, $|\Gamma|> C \log n$,
$C$ a constant.
(Here we have already
(i) replaced $G$ by the preimage of $M$ in the reduction $G \to G/N$ and, after that, (ii) replaced $N$ by the stabilizer of the blocks in part (\ref{it:rilb})
of Theorem~\ref{thm:cam}.) This case will occupy us for the rest of the article.

\begin{center}
  * * *
\end{center}

Babai indicates how to remove the dependence on CFSG at this point.
Let $G$ and $N$ be as before, with $G$ transitive. Then $G/N$ is
a primitive group acting on the set of blocks $R$.

If a permutation group on a set $R$ is such that its action on
the set of pairs of distinct elements of $R$ is transitive, the
group is said to be {\em doubly transitive}. Now, a result of Pyber
that does not depend on CFSG \cite{Py} tells us that such a group
is either $\Alt(R)$, or $\Sym(R)$, or of order
$\leq |R|^{O(\log^2 |R|)}$.

If $G/N$ is $\Alt(R)$ or $\Sym(R)$, we are in the case that we will discuss
from here until the end. If $G/N$ is doubly
transitive, but is equal neither to $\Alt(R)$ nor to $\Sym(R)$,
we can proceed as in the transitive case of the proof
of Theorem \ref{thm:luxor}, because
$|G/N|\leq r^{O(\log^2 r)}$, $r = |R|\leq n$.
(Babai proposes also an alternative treatment,
even more efficient and elementary.)

Let us thus suppose that $G/N$ is not doubly transitive. Then
the Schurian coherent configuration (\S \ref{subs:concoh}) that it induces is 
not a clique. Therefore, we can give
this configuration to the procedure
\textsc{Split-or-Johnson} (\S \ref{sec:coujoh}), and resume the argument at that point. 

\section{The structure of the action of $\Alt$}

\subsection{Stabilizers, orbits and alternating quotients}\label{subs:stab}

We will need several results on the epimorphisms $G\to \Alt_k$.
They will play a crucial role in the method of local certificates (\S \ref{sec:certloc}).
In the original version \cite{Ba}, they also played a similar role as \cite{BLS} in this article.

\begin{lemm}\label{lem:mustaf}
  Let $G<\Sym(\Omega)$ be primitive. Let $\phi:G\to \Alt_k$ be an epimorphism 
  with $k>\max(8,2+\log_2 |\Omega|)$. Then $\phi$ is an isomorphism.
\end{lemm}
Proving this lemma is more or less an exercise in finite group theory; 
 \cite[Prop.~1.22]{BaPS} must be used for the case of abelian socle and
 Schreier's conjecture for the case of non-abelian socle (see Exercise~\ref{hidden3}).
Schreier's conjecture is a theorem, but a theorem whose proof depends, in turn, on CFSG. 

On the other hand, Pyber \cite{Py2} gave a proof of Lemma \ref{lem:mustaf}
that does not use CFSG, with a stricter condition:
$k> \max(C,(\log |\Omega|)^5)$, $C$ constant.
The dependence on CFSG has thus been completely removed  from the proof of the 
main theorem.

\begin{defi}
  Let $G<\Sym(\Omega)$. Let $\phi:G\to \Sym_k$ be a homomorphism whose
  image contains $\Alt_k$. Then $x\in \Omega$ is said to be {\em affected}
    if
  $\phi(G_x)$ does not contain $\Alt_k$.
\end{defi}

\begin{lemm}\label{lem:preaffect} Let $G<\Sym(\Omega)$.
  Let $\phi:G\to \Alt_k$ be an epimorphism with
  $k>\max(8,2+\log_2 n_0)$, where 
  $n_0$ is the size of the largest orbit of $G$.
  \begin{enumerate}
  \item\label{it:transitive}
    If $G$ is transitive, every $x\in \Omega$ is affected.
  \item\label{it:atr}
    At least one $x\in \Omega$ is affected.
  \end{enumerate}
\end{lemm}
\noindent{\sc Proof (sketch)} ---
(\ref{it:transitive})  This results directly from Lemma \ref{lem:mustaf} if $G$ is primitive,
or if $K<\ker(\phi)$ for $K$ the stabilizer of a minimal system of blocks. 
It remains the case of $\phi:K\to \Alt_k$ surjective. In general:
\begin{quote} {\sc Lemma.}---
For $K_i$ arbitrary, $K<K_1\times\dots \times K_s$ and an epimorphism
\mbox{$\phi:K\to S$,} $S$ simple, there must be an $i$ such that $\phi$ factors
as follows: \mbox{$K\to K_i \stackrel{\psi}{\to} S$,} $\psi$ an epimorphism.\end{quote}
Using this lemma for the restrictions $K_i$ of $K$ to the orbits of $K$, we move to an orbit $K_i$,
and proceed by induction.

(\ref{it:atr}) Let $\Omega_1,\dotsc,\Omega_m$ be the orbits of $G$, and
let $G_i = G|_{\Omega_i}$ be the restriction of $G$ to $\Omega_i$. By the lemma in (\ref{it:transitive}),
there must be an $i$ such that $\phi$ factors as
$G\to G_i \stackrel{\psi}{\to} \Alt_k$, $\psi$ an epimorphism.
Then, by (\ref{it:transitive}), $(G_x)^\psi = ((G_i)_x)^\psi \ne \Alt_k$
for every $x\in \Omega_i$.
\qed

The following proposition will play a crucial role in \S \ref{sec:casimp}. 
\begin{prop}\label{prop:atinl}
  Let $G<\Sym(\Omega)$ be transitive
  and $\phi:G\to \Alt_k$ be an epimorphism. Let
  $U\subset \Omega$ be the set of unaffected elements. 
  \begin{enumerate}
  \item\label{it:unaffstab} Suppose that $k\geq \max(8,2+\log_2 n_0)$,
    where $n_0$ is the size of the largest orbit of $G$.
    Then $(G_{(U)})^\phi = \Alt_k$.
  \item\label{it:afforb} Suppose that $k\geq 5$. If $\Delta$ is an orbit
    of $G$ containing some unaffected elements, then each orbit of
    $\ker(\phi)$ contained in $\Delta$ is of length $\leq |\Delta|/k$.
  \end{enumerate}
\end{prop}
Remember that $G_{(U)} = \{g\in G: x^g = x \; \forall x\in U\}$
(pointwise stabilizer).

\noindent{\sc Proof} ---
(\ref{it:unaffstab}) It is easy to see that
$G$ fixes $U$ as a set. Then,
$G_{(U)}\triangleleft G$, and thus $(G_{(U)})^\phi \triangleleft G^\phi$.
Now, $G^\phi = \Alt_k$.
Suppose that $(G_{(U)})^\phi = \{e\}$. Then $\phi$ factors as follows:
\[G\to G|_U \stackrel{\psi}{\to} \Alt_k,\]
since $G_{(U)}$ is the kernel of $G\to G|_U$. 
Here $\psi$ is an epimorphism, and thus, by Lemma
\ref{lem:preaffect} (\ref{it:afforb}), there exists an $x\in U$ such that
$((G|_U)_x)^\psi \ne \Alt_k$. Now $((G|_U)_x)^\psi = (G_x)^\phi = \Alt_k$,
because $x$ is in $U$, i.e. unaffected. Contradiction.

(\ref{it:afforb}) As $\Delta$ contains some affected elements and is an orbit of $G$,
every element of $\Delta$ is affected.
Let $N=\ker(\phi)$, $x\in \Delta$. 
The length of the orbit $x^N$ is
\[\begin{aligned}
\left|x^N\right| &= \lbrack N:N_x\rbrack = \lbrack N:(N\cap G_x)\rbrack
= \lbrack N G_x : G_x\rbrack = \frac{\lbrack G:G_x\rbrack}{\lbrack G : N
  G_x\rbrack}\\ &= \frac{|\Delta|}{\lbrack G^\phi: (G_x)^\phi\rbrack}
= \frac{|\Delta|}{\lbrack \Alt_k : (G_x)^\phi\rbrack}.\end{aligned}\]
Now, every proper subgroup of $\Alt_k$ is of index $\geq k$. Thus
$\left|x^N\right|\leq |\Delta|/k$.
\qed

\subsection{The case of large symmetry}\label{sec:grasym}

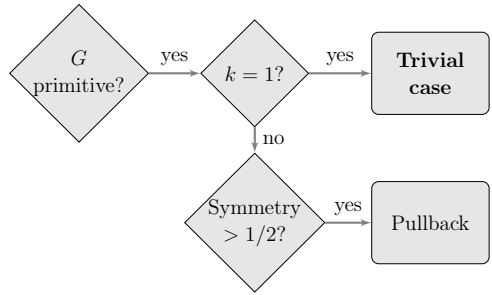
\begin{center}
    \scalebox{0.66}{
    \begin{tikzpicture}[scale=2.5, node distance = 2cm, auto]
         \node [decision] (primi)
          {$G$ primitive?};
      \node [decision, right of = primi, node distance = 3.5cm] (k1) {$k=1$?};
      \node [block, right of =k1, node distance = 3.5cm] (triv) {\bf Trivial case};
      \node [decision, below of = k1] (sym) {Symmetry $>1/2$?};
      \node [block, right of =sym, node distance = 3.5cm] (pull) {Pullback};
          \path [line] (primi) -- node [, color=black] {yes} (k1);
      \path [line] (k1) -- node [, color=black] {yes} (triv);
      \path [line] (k1) -- node [, color = black] {no} (sym);
      \path [line] (sym) -- node [, color = black] {yes} (pull);
    \end{tikzpicture}
  }
\end{center}

Consider the case of $G$ primitive. We can suppose that
$G$ is isomorphic as a permutation group to $\Alt_m^{(k)}$,
since we have already eliminated the other cases in \S \ref{sec:procprin}
(possibly by passing to an imprimitive group $M$; the imprimitive 
case will be treated in \S \ref{sec:casimp}).
As we have seen in
\S \ref{subs:idgroup}, we can build a bijection $\iota$
between $\Omega$ and the set $\mathscr{S}_k(\Gamma)$ of subsets with $k$ elements of a set $\Gamma$. This bijection induces 
an isomorphism $\phi:G\to \Alt(\Gamma)$.

If $k=1$, then $\Omega$ is in bijection with $\Gamma$, and $G\sim \Alt_n =
\Alt_m$.
We are thus in the trivial case, the group $\Aut_G(\mathbf{x})$
consists of the elements of $\Alt_n$ permuting the letters of
$\mathbf{x}$ of the same color, and $\Iso_G(\mathbf{x},\mathbf{y})$
is nonempty if and only if $\mathbf{x}$ and $\mathbf{y}$ have exactly the same
number of letters of the same color -- where, if no letter is repeated either in 
$\mathbf{x}$ or in $\mathbf{y}$, we add the condition that the
permutation of $\{1,\dotsc,n\}$ inducing $\mathbf{x}\mapsto \mathbf{y}$
is in $\Alt_n$.

Thus, let $G$ primitive, $k>1$. 

Two elements $\gamma_1,\gamma_2\in \Gamma$ are {\em twins}
with respect to an object if the transposition $(\gamma_1 \gamma_2)$ leaves it unchanged. 
It is clear that the twins form equivalence classes, and that, 
for every such equivalence class $C$, every $\Sym(C)$ leaves the object
unchanged. Our object will be the chain $\mathbf{x}$ (or $\mathbf{y}$):
$\gamma_1$, $\gamma_2$ are twins with respect to $\mathbf{x}$ if, for every $i\in \Omega$,
$\mathbf{x}(i) = \mathbf{x}(\tau^{\phi^{-1}}(i))$, where
$\tau = (\gamma_1 \gamma_2)$.

We can thus easily determine (and in polynomial time)
the equivalence classes in
$\Gamma$ (called {\em twin classes}),
and verify whether there is an equivalence class $C$
of size $>|\Gamma|/2$.
Let us examine this possibility as we will have to exclude it later on. 

The class $C$ of size $>|\Gamma|/2$ is obviously unique and thus canonical. If
$\mathbf{x}$ has such a class and $\mathbf{y}$ does not have it, or if
both have such classes, but of different sizes, then 
$\mathbf{x}$ and $\mathbf{y}$ are not isomorphic.

If $\mathbf{x}$, $\mathbf{y}$ have twin classes
$C_{\mathbf{x}}$, $C_{\mathbf{y}}$ of the same size $>|\Gamma|/2$,
we choose $\sigma\in \Alt(\Gamma)$ such that
$C_{\mathbf{x}} = \left(C_{\mathbf{y}}\right)^\sigma$.
(We suppose $m>1$.)
By replacing $\mathbf{y}$
by $\mathbf{y}^{\sigma'}$, where $\sigma' = \phi^{-1}\left(\sigma^{-1}\right)$, we reduce our problem to the case 
$C_{\mathbf{x}} = C_{\mathbf{y}}$. (Here is the simplest example of
what Babai calls {\em to align}; we have {\em aligned}
$C_{\mathbf{x}}$ and $C_{\mathbf{y}}$.)

Then, let $C = C_\mathbf{x} = C_{\mathbf{y}}$. 
The partition
$\{C,\Gamma\setminus C\}$ of $\Gamma$ induces a partition
$\{\Omega_j\}_{0\leq j\leq k}$ of $\Omega$: $\omega\in \Omega_j$
  if and only if $\psi(\omega)$ contains $k-j$ elements of $C$ and $j$ elements
  of $\Gamma\setminus C$.
It is easy to show that $\alpha^{k-j} (1-\alpha)^j \binom{k}{j}
  <1/2$ for $\alpha\in (1/2,1\rbrack$, $1\leq j\leq k$; thus,
  $|\Omega_j| < n/2$ for $1\leq j\leq k$.

We have reduced our problem to that of determining 
  $\Iso_H(\mathbf{x},\mathbf{y})$, where
  $H = \phi^{-1}\left(\Alt(\Gamma)_C\right)$. Here the need to take 
  a setwise stabilizer (namely,
  $\Alt(\Gamma)_C$) does not create any problems: we generate $H$
  by taking preimages $\phi^{-1}(h_1),\dotsc,\phi^{-1}(h_5)$
  of two generators $h_1$, $h_2$ of $\Alt(C)<\Alt(\Gamma)$, two generators
  $h_3$, $h_4$ of $\Alt(\Gamma\setminus C)<\Alt(\Gamma)$ and
  an element $h_5\in \Alt(\Gamma)$ of the form $(\gamma_1 \gamma_2)
  (\gamma_3 \gamma_4)$, where $\gamma_1,\gamma_2\in C$, $\gamma_3,\gamma_4\in
  \Gamma\setminus C$.
  (If $|\Gamma|<8$, the number of generators is smaller, and the
  discussion gets simpler.) 
 Our problem reduces to that of determining
  $\Iso_{H'}(\mathbf{x},\mathbf{y}')$ for $\mathbf{y}' =
  \mathbf{y}$ and $\mathbf{y}' = \mathbf{y}^{h_5}$,
  where $H' = \phi^{-1}(\Alt(C)\times
  \Alt(\Gamma\setminus C)) = \phi^{-1}(\langle h_1,\dotsc,h_4\rangle)$.
  
  As $C$ is a twin class for $\mathbf{x}$, every element of
  $\phi^{-1}(\Alt(C))$ leaves $\mathbf{x}$ unchanged.
  If $\mathbf{x}|_{\Omega_0} \ne \mathbf{y}|_{\Omega_0}$,
  then $\Iso_{H'}(\mathbf{x},\mathbf{y}) = \emptyset$.

  Let thus $\mathbf{x}|_{\Omega_0} = \mathbf{y}|_{\Omega_0}$.
  We have reduced our problem to that of determining
  $\Iso_{H'|_{\Omega'}}(\mathbf{x}|_{\Omega'},\mathbf{y}|_{\Omega'})$,
  where $\Omega' = \Omega  \setminus \Omega_0$.
  Remember that $H'|_{\Omega'}$ acts on $\Omega'$ with orbits of
  length $|\Omega_i|<n/2$. We thus proceed as in the intransitive case of Luks's method (Proof of Thm.~\ref{thm:luxor}).

\section{Descent from a relational structure}\label{subs:chasche}

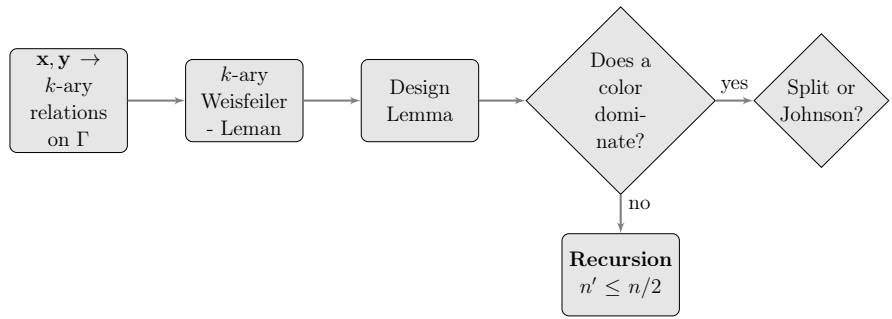
\begin{center}
      \scalebox{0.66}{
\begin{tikzpicture}[scale=2, node distance = 2cm, auto]
    \node [block, node distance = 3cm] (xrel) {$\mathbf{x},
    \mathbf{y}\to$ $k$-ary relations on $\Gamma$};
    \node [block, right of =xrel, node distance = 3.5cm] (wl) {$k$-ary Weisfeiler - Leman};
    \node [block, right of = wl, node distance = 3.5cm] (design) {Design Lemma};
    \node [decision, right of = design, node distance = 4cm] (coldom) {Does a color dominate?};
    \node [decision, right of = coldom, node distance = 4cm] (spj) {Split or {\textsc Johnson}?};
    \node [block, below of = coldom, node distance = 3.5cm] (luks) {\bf Recursion $n'\leq n/2$};
    \path [line] (xrel) -- (wl);
    \path [line] (wl) -- (design);
    \path [line] (design) -- (coldom);
    \path [line] (coldom) -- node [, color = black] {yes} (spj);
        \path [line] (coldom) -- node [near start, color = black] {no} (luks);
\end{tikzpicture}
}
\end{center}

Let us discuss now the case of $G$ primitive and, more precisely, $G$ isomorphic to $\Alt_m^{(k)}$, $k\geq 2$. Now we can suppose that our strings $\mathbf{x}, \mathbf{y}$ do not have twin classes of size $>m/2$. The main tools that we will develop (Design Lemma, Split-or-Johnson) will be useful to us, even essential, also in the case of $G$ imprimitive.

We have a bijection between the elements of $\Omega$
and $\{S\subset \Gamma: |S|=k\}$. For a given $\mathbf{x}:\Omega\to \Sigma$, we have then a $k$-ary relational structure $\mathfrak{X}=(\Gamma, (R_i)_{i\in \Sigma})$ on $\Gamma$: $(x_1,\dotsc,x_k)\in R_i$ if $x_1,\dotsc,x_k$ are all different and $\mathbf{x}(\omega) = i$, where $\omega$ is the element of $\Omega$ that corresponds to $\{x_1,\dotsc,x_k\}$.

We apply to $\mathfrak{X}$ the functor $F_1$ (\S \ref{subs:secf1}), that makes it a partition structure, then the functor $F_2$ (\S \ref{subs:secf1}), that gives us a $k$-ary configuration, and, finally, the functor $F_3$ defined by $k$-ary Weisfeiler-Leman (\S \ref{subs:wl}). We obtain thus a refinement $F_3(F_2(F_1(\mathfrak{X}))) =
(\Gamma,c_\mathbf{x}:\Gamma^k\to \mathscr{C})$ that is a $k$-ary coherent configuration.

Since $F_1$, $F_2$, $F_3$ are functors, the assignment of $c_\mathbf{x}$
to $\mathbf{x}$ is canonical. It will be then useful to us: if $c_{\mathbf{x}}$ and $c_{\mathbf{y}}$ are not isomorphic under the action of $\Alt_m$, then $\mathbf{x}$ and $\mathbf{y}$ are not isomorphic under the action of $\Alt_m^{(k)}$ either.

We will obtain canonically a classical coherent configuration from $c_\mathbf{x}$ (Design Lemma). Either this new configuration is non-trivial or we obtain a canonical coloring without a dominant color. Such a coloring will immediately allow us to reduce the problem to a certain number of problems for shorter strings, as in Luks's algorithm.

Suppose then that we have a non-trivial classical coherent configuration canonically assigned to $\mathbf{x}$. The \textsc{Split-or-Johnson} procedure will give us one of these two results: either a canonical {\em colored partition} of $\Gamma$, or a {\em Johnson scheme} canonically embedded in $\Gamma$. In both cases, having such a canonical structure restricts greatly the set of possible isomorphisms and automorphisms. We will be able to reduce $G$ to a subgroup $\sim \Alt_{m'}$, with $m'\leq m/2$ in the case of a colored partition or $m'\ll \sqrt{m}$ in the case of a Johnson scheme. The bound $m'\leq m/2$ is already sufficient for a successful recursion.

\subsection{Design Lemma}\label{sec:design}

Given a configuration $\mathfrak{X}=(\Gamma,c:\Gamma^k\to\mathscr{C})$ and a parameter $1/2\leq \alpha<1$, a color $i$ is called {\em $\alpha$-dominant} if $c(\gamma,\dotsc,\gamma) = i$ for $\geq \alpha |\Gamma|$ values of $\gamma\in \Gamma$. The color class $\{\gamma\in \Gamma: c(\gamma,\dotsc,\gamma) = i\}$ is also called {\em dominant}. If on the other hand for every color $i$ the class $\{\gamma\in \Gamma: c(\gamma,\dotsc,\gamma) = i\}$ is of size $<\alpha|\Gamma|$, we say that our coloring is an {\em $\alpha$-coloring}.

As before, two elements $\gamma_1,\gamma_2\in \Gamma$ are {\em twins} with respect to a structure~$\mathfrak{X}$ (here, a coherent configuration on $\Gamma$) if $(\gamma_1 \gamma_2) \in \Aut(\mathfrak{X})$.

\begin{prop}[Design Lemma]\label{prop:design}
  Let $\mathfrak{X}=(\Gamma, c:\Gamma^k\to \mathscr{C})$
  a $k$-ary coherent configuration, where
  $2\leq k\leq |\Gamma|/2$. Let $1/2\leq \alpha < 1$.
  Suppose that there is no twin class
  in $\Gamma$ with $>\alpha |\Gamma|$ elements.

  Then, at least one of the following options is true:
  \begin{enumerate}
  \item\label{it:loro1}
    there exist $x_1,\dotsc,x_\ell\in \Gamma$,
    $0\leq \ell<k$, such that
    $\mathfrak{X}_{\vec{x}}^{(1)}$ has no $\alpha$-dominant color;
  \item\label{it:loro2} 
    there exist $x_1,\dotsc,x_\ell\in \Gamma$,
    $0\leq \ell<k-1$, such that
    $\mathfrak{X}_{\vec{x}}^{(1)}$ has an $\alpha$-dominant color $C$
    and $(\mathfrak{X}_{\vec{x}})^{(2)}\lbrack C\rbrack$ is not a clique.
  \end{enumerate}
  \end{prop}
The notation has been defined in \S\S \ref{subs:secf1} -- \ref{subs:concoh}. In particular, the $1$-skeleton $\mathfrak{X}_{\vec{x}}^{(1)}$ is simply a coloring of $\Gamma$. 

\begin{lemm}[Large clique lemma]
  Let $\mathfrak{X} = (\Gamma,c)$ be a classical coherent configuration.
  Let $C\subset \Gamma$ be a color class with
  $|C|\geq |\Gamma|/2$. If $\mathfrak{X}\lbrack C\rbrack$ is a clique, then
  $C$ is a twin class.
\end{lemm}

\noindent{\sc Proof} --- Suppose that $C$ is not a twin class. Then there is an $x\in \Gamma$ and a color (say {\em azure}) such that $c(x,y)$ is equal to this color for at least one $y\in C$ but not for all of them. Since $\mathfrak{X}\lbrack C\rbrack$ is a clique, $x\notin C$. Call the color of $C$ {\em crimson}, and that of $x$ {\em bronze}. Let $B\subset \Gamma$ be the set of bronze elements.

We want to construct a balanced block design (\S \ref{subs:hypdes}) that contradicts Fisher's inequality (Prop.~\ref{prop:fishy}). Define $A_b = \{y\in \Gamma : c(b, y) = \text{azure}\}$ for $b\in B$. Since $x\in B$ and $c(x, y) = \text{azure}$ for at least one $y\in C$, and $c(x,y)$ knows the color of $y$, all elements of $A_b$ are crimson.

By the coherence of $\mathfrak{X}$ and the definition of structure constants (Def.~\ref{def:coh}), \[|A_b| = \gamma(\text{azure}^{-1},\text{azure},\text{bronze}),\] so $|A_b|$ does not depend on $b$. As we said at the beginning, $1\leq |A_x|< |C|$; therefore, $1\leq |A_b|<|C|$ for all $b\in B$.

Show similarly that, for $v\in C$, the size of $\{b\in B: v\in A_b\} = \{b\in B: c(b, v) = \text{azure}\}$ does not depend on $b$. Since $\mathfrak{X}\lbrack C\rbrack$ is a clique,  $c(v,v')$ is of the same color for all $v,v'\in C$, $v \ne v'$; call this color {\em daffodil}. Show that \[ |\{b\in B: v,v'\in A_b\}| = \gamma(\text{azure},\text{azure}^{-1},\text{daffodil}).\] Then $(C,\{A_b\}_{b\in B})$ is an incomplete balanced block design.

Therefore, by Fisher's inequality, $|B|\geq |C|$. Now, we know that $|C|>|\Gamma|/2$, $B,C\subset \Gamma$ and $B\cap C = \emptyset$. Contradiction.
\qed

\noindent{\sc Proof of the Design Lemma} (Prop. \ref{prop:design}) ---
Suppose that
for every $\vec{x} \in \Omega^\ell$,
$0\leq \ell < k$,
$C_{\vec{x}}$ has an $\alpha$-dominant color $C(\vec{x})$, and moreover,
if $\ell<k-1$, $(\mathfrak{X}_{\vec{x}})^{(2)}\lbrack C(\vec{x})\rbrack$
is a clique. We will come to a contradiction.

Let $C = C(\text{empty})$. Since $|C|>\alpha |\Gamma|$, $C$ is too big to be a twin set. Then there exist $u,v\in C$, $u\ne v$, such that $\tau = (u v) \notin \Aut(\mathfrak{X})$. Let $\vec{y}$ be of minimal length $r$ among the strings satisfying $c(\vec{y}^{\,\tau}) \ne c(\vec{y})$. By this minimality and the rules in Definition \ref{def:ukuru}, $y_1,\dotsc y_r$ are all distinct. Permuting them, we can ensure that $u,v\notin \{y_1,y_2,\dotsc,y_{r-2}\}$, and, without loss of generality, that either (i) $y_{r-1}\ne u,v$, $y_r=u$, or (ii) $y_{r-1} = u$, $y_r=v$. In case (i), we choose $\vec{x} = y_1,\dotsc,y_{r-1}$, $\ell=r-1$, and we see that $c_{\vec{x}}(u) \ne c_{\vec{x}}(v)$; in case (ii), we choose $\vec{x} = y_1,\dotsc,y_{r-2}$, $\ell=r-2$, and we obtain $c_{\vec{x}}(u,v) \ne c_{\vec{x}}(v,u)$. We will then have a contradiction to our hypothesis once we will have proved that $u,v\in C(\vec{x})$.

The fact that $u,v\in C(\vec{x})$ will immediately follow from the equality $C(\vec{x}) = C\setminus \{x_1,\dotsc,x_\ell\}$; this equality in turn can be deduced from the fact that, for $\vec{y}$ of length $\leq k-2$ and $\vec{x} = \vec{y} z$, $z\in \Omega$,
\begin{equation}\label{eq:kurut}
  C(\vec{x}) = C(\vec{y})\setminus \{z\}.\end{equation}
Why is (\ref{eq:kurut}) true? We are supposing that $\mathfrak{X}_{\vec{y}}^{(2)}\lbrack C(\vec{y})\rbrack$ is a clique, and that $|C(\vec{y})|>\alpha |\Gamma|\geq |\Gamma|/2$. Therefore, by the large clique lemma, all the elements of $C(\vec{y})$ are twins in $\mathfrak{X}_{\vec{y}}^{(2)}$. In particular, for $u\in C(\vec{y})\setminus \{z\}$, $c_{\vec{x}}(u) = c_{\vec{y}}(z u)$ does not depend on $u$. Since the coloring of vertices in $C_{\vec{x}}$ is a refinement of that in $C_{\vec{y}}$ (by the second rule of Definition~\ref{def:ukuru}), it follows that, either $C(\vec{x}) = C(\vec{y})\setminus \{z\}$, or $C(\vec{x})\subset \Gamma \setminus C(\vec{y})$, or $C(\vec{x}) = \{z\}$. Since $|C(\vec{x})|, |C(\vec{y})| > \alpha |\Gamma| \geq |\Gamma|/2$, the two last possibilities are excluded.
\qed

We apply the Design Lemma (with $\alpha=1/2$) to the $k$-ary coherent configuration $\mathfrak{X}' = F_3(F_2(F_1(\mathfrak{X})))$, where $\mathfrak{X}$ is given by $\mathbf{x}$ in the fashion described at the beginning of the section. We run through all the possible tuples $\vec{x} = (x_1,\dotsc,x_{\ell}) \in \Gamma^\ell$, $0\leq \ell<k$, until we find a tuple for which the first or second conclusion of the Design Lemma is true.

If the first conclusion is true, we define $c_{\mathfrak{X}} = \mathfrak{X}_{\vec{x}}^{(1)}$ and we skip to \S\ref{subs:cacoudo}. If the second conclusion is true, we go to \S \ref{sec:coujoh}, having defined $\mathfrak{X}'' = \mathfrak{X}_{\vec{x}}^{(2)}\lbrack C\rbrack$, where $C$ is the $\alpha$-dominant color of $\mathfrak{X}_{\vec{x}}^{(1)}$.

\subsection{Split-or-Johnson}\label{sec:coujoh}

We have a non-trivial homogeneous classical coherent configuration $\mathfrak{X}'' = (\Gamma,c)$. (We remind that it is a coloring $c$ of the complete graph on $\Gamma$ such that (a) the vertices have their own color ({``diagonal color''}), (b) the edges $(x,y)$, $x\ne y$, are not all of the same color, (c) the color $c(x,y)$
of the edge $(x,y)$ determines $c(y,x)$, and (d) the axiom of coherence (\ref{def:coh}) is verified.) We want to find structures that depend canonically on $\mathfrak{X}''$ and that reduce its group of automorphisms.

It is reasonable to expect that these structures exist: by Theorem \ref{thm:cam}, if the group of automorphisms is transitive, either it is imprimitive (so that it leaves invariant a partition), or it is close to being $\Alt_m^{(k)}$, $k\geq 2$ (which leaves a Johnson scheme invariant), or it is small (so that the stabilizer of some point will have small orbits, and thus will give us a coloring without a dominant color). The challenge is to find these structures, and to do it canonically.

If $\mathfrak{X}''$ is not primitive (Def.~\ref{def:uniprim}), the task is rather easy: let $r$ be the reddest non-diagonal color such that the graph $\mathscr{G}_r = \{(x,y): x,y\in \Gamma,\; c(x,y)=r\}$ is not connected; by Exercise \ref{ex:samesize}, this gives a partition of $\Gamma$ into sets of the same size  $\leq |\Gamma|/2$.

\begin{theo}[Split-or-Johnson]\label{thm:coupjohn}
  Let $\mathfrak{X} = (\Gamma,c)$ be a uniprimitive classical coherent configuration.
  Let $2/3\leq \alpha< 1$.
  In time $|\Gamma|^{O(1)}$, we can find
  \begin{itemize}
  \item either a colored $\alpha$-partition of $\Gamma$,
  \item or a Johnson scheme embedded in
     $\Gamma_0\subset \Gamma$, $|\Gamma_0|\geq \alpha |\Gamma|$,
  \end{itemize}
  and a subgroup $H<\Sym(\Gamma)$ with \[\lbrack \Sym(\Gamma):H\rbrack
  = |\Gamma|^{O(\log |\Gamma|)}\]
  such that the colored partition, or the scheme, is canonical with respect to $H$.
  \end{theo}

We will define the group $H$ as a pointwise stabilizer in $\Gamma$.
The value $2/3$ in the statement is rather arbitrary; actually,
any value greater than $1/2$ would do, though values close to $1/2$ would
affect the implied constants.

\noindent{\sc Proof } ---
Choose an arbitrary $x\in \Gamma$. Give each $y\in \Gamma$ the color $c(x,y)$. This coloring is canonical with respect to $G_x$. If there is no color class $C$ of size $>\alpha |\Gamma|$, the trivial partition (non-partition) of each class will give us a colored $\alpha$-partition of $\Gamma$, and we are done.

Suppose on the contrary that there is a color class -- say $C_{\text{aqua}}$ -- of size $> \alpha |\Gamma|$. Since $\alpha n>n/2$, the relation $R_{\text{aqua}}$ of this color is non-oriented ($c(y,z) = \text{aqua}$ iff $c(z,y) = \text{aqua}$). The complement of $R_{\text{aqua}}$ (or of any other relation) is of diameter $2$ (see Exercise~\ref{hidden4}). Let $x,z\in \Gamma$ be such that $c(x,z)=\text{aqua}$, and let $y\in \Gamma$ be such that $c(x,y), c(z,y)\ne \text{aqua}$. Let us call $c(x,y)$ {\em beige} and $c(z,y)$ {\em cyan}.

Consider the bipartite graph $(V_1,V_2;A)$ with vertices $V_1 = C_{\text{aqua}}$, $V_2 = C_{\text{beige}}$ and edges $R_{\text{cyan}}\cap (V_1\times V_2)$. The graph is nonempty by definition and semiregular by Exercise \ref{ex:biparhom}\ref{it:bara2}. By homogeneity and coherence, the number of $y$ such that $c(y,w)$ is of a given color $c_0$ is independent from $w$. Therefore, it is always $\leq (1-\alpha) n < n/2$ for $c_0\ne \text{aqua}$. Applying this to $c_0 = \text{cyan}$ and $V_2$, we see that the degree $|\{v_1\in V_1: (v_1,v_2)\in A\}|$ is $<n/2$, so, since $|V_1|>n/2$, the graph is not complete.

We apply then Proposition  \ref{prop:bicoup} to $(V_1,V_2;A)$ with $\beta = \alpha |\Gamma|/|V_1|$. Notice that \mbox{$|V_2|\leq \beta |V_1|$.}
\qed

We will work then with a bipartite graph $(V_1,V_2;A)$. The strategy will be to try either to make $V_2$ smaller (by at least a constant factor), or to find structures in it. Either these structures will allow us to reduce $V_2$ all the same, or they will help us to partition $V_1$, or to find a large enough Johnson scheme on~$V_1$.

First, we will have to limit the symmetry in $V_1$, i.e. reduce, even eliminate the twins. There are two reasons for that.
\begin{itemize}
\item Even if we discover a rich enough structure in $V_2$, this will imply nothing in $V_1$ if enough elements of $V_1$ are connected to $V_2$ in the same way.
\item If $V_2$ is small, we will color each vertex of $V_1$ according to its neighbors in $V_2$. This will give us a canonical coloring with respect to $G_{(V_2)}$. Now, in this coloring, two vertices in $V_1$ will have the same color iff they are twins; then, if no twin class in $V_1$ has $>\alpha |V_1|$ elements, we will have an $\alpha$-coloring.
\end{itemize}

\begin{exer}\label{ex:gita}
  Let $(V_1,V_2;A)$ be a non-trivial semiregular bipartite graph. Then, no twin class in
  $V_1$ has more than $|V_1|/2$ elements.
\end{exer}

\noindent{\sc Solution} ---
  We will ensure that $|A|\leq |V_1| |V_2|/2$ by taking the complement if necessary.
  Let $d_2$ be the degree of vertices in $V_2$ and $S$ a twin class in $V_1$.
  Show that $d_2\geq |S|$, and then $|A|\geq |S| |V_2|$.
\qed

  \begin{exer}\label{ex:reduc2}
    Let $(V_1,V_2;A)$ be a bipartite graph without twins in $V_1$.
    Let $V_2 = C_1\cup C_2$, $C_1\cap C_2=\emptyset$.
    Show that, for at least one $i=1,2$,
    there is no class of $\geq |V_1|/2+1$ twins in $V_1$
    in the graph $(V_1,C_i;A\cap (V_1\times C_i))$.
  \end{exer}

  \begin{exer}\label{ex:twindet}
    Let $\mathfrak{X}=(\Gamma,c)$ be a classical coherent configuration. Let
    $C_1$, $C_2$ be two color classes in $\Gamma$. Let {\em brown} be
    an edge color in $C_1 \times C_2$. Then, for $x,y\in C_1$,
    the color $c(x,y)$ determines whether $x$ and $y$ are twins in the bipartite graph
    $(C_1,C_2;\mathscr{G}_{\text{brown}})$.
\end{exer}

\begin{prop}[Bipartite Split-or-Johnson, or {``Una partita a poker''}]\label{prop:bicoup}
  Let $X = (V_1,V_2;A)$ be a bipartite graph
  with $|V_2|<\beta |V_1|$, where $2/3\leq \beta<1$, and such that no twin class
  in $V_1$ has more than $2|V_1|/3$ elements.
  Then, we can find, in time  $|V_1|^{O(1)}$,
  \begin{itemize}
  \item either a colored $\beta$-partition of $V_1$,
  \item or a Johnson scheme embedded in 
    $V_0\subset V_1$, $|V_0|\geq \beta |V_1|$,
  \end{itemize}
  and a subgroup $H<G$, $G = \Sym(V_1)\times \Sym(V_2)$, with
  \[\lbrack G:H\rbrack = |V_1|^{O(\log |V_1|)}\]
  such that the colored partition, or the scheme,
  is canonical with respect to $H$.
\end{prop}
The graph obtained at the end of the proof of Thm.~\ref{thm:coupjohn} satisfies the condition on the twin classes given here (even for $1/2$ instead of $2/3$), thanks to Exercise~\ref{ex:gita}.

With regards to the time of the procedure, we will make explicit some details that could seem non-evident. What is most delicate is the index $\lbrack G:H\rbrack$. We will define the group $H$ as a pointwise stabilizer; we will have to control the number of points that we stabilize.

Let us sketch now the basic strategy of the proof. What we want is a reduction
to Prop.~\ref{prop:cohcoup}, ``Coherent Split-or-Johnson''.
We can construct a classical coherent configuration $\mathfrak{X}$ on
$V_1 \cup V_2$ starting from the graph $X$,
 simply by Weisfeiler-Leman. The tricky part will be to guarantee
 that the restriction $\mathfrak{X}\lbrack C_2\rbrack$ to the dominant color
 class $C_2\subset V_2$ (if there is one) is non-trivial.

 In order to obtain a non-trivial configuration on $C_2$, we will note that
 the graph $X$ itself induces a $d$-ary relation on $C_2$, where
 $d$ is at most the degree of the majority of elements of $V_1$ (should there be
 such a thing; otherwise degrees give us a partition of $V_1$). If this
 relation is trivial, in the sense of
 containing all $d$-tuples of distinct elements in $C_2$, we
 obtain a Johnson scheme. If it is non-trivial but has many twins,
 it gives us a way to descend recursively to a smaller $C_2$. If it does
 not have many twins, we can use the Design Lemma (supplemented by
 a standard lemma on designs) to obtain the non-trivial classical coherent
 configuration on $C_2$ that we want.

 \noindent{\sc Proof } ---
If $|V_1|\leq c$, where $c$ is a constant, we color each $v\in V_1$ by itself. This coloring is canonical with respect to $H=\{e\}$; in other words, it is not canonical at all. It does not matter: trivially, $|G|\leq (c!)^2 \leq |V_1|^{O(\log |V_1|)}$. We then can suppose that $|V_1|>c$.

If $|V_2|\leq (6 \log |V_1|)^{3/2}$ (say), then, by the reasoning above, we obtain a $(2/3)$-coloring of $V_1$ (and then a colored $(2/3)$-partition of $V_1$). This coloring is canonical with respect to an $H$ of index
\[|V_2|! \leq |V_2|^{|V_2|} \leq (6 \log |V_1|)^{\frac{3}{2}
   (6 \log |V_1|)^{\frac{3}{2}}} \ll |V_1|^{\log |V_1|}.\]
We then can suppose that $|V_2|> (6 \log |V_1|)^{3/2}$.

Our first task is to eliminate the twins. We divide $V_1$ into its twin classes and we color each $v\in V_1$ by its number of twins and by its degree in the graph $(V_1,V_2;A)$. We obtain a colored $\beta$-partition of $V_1$, unless there is an integer $d_1$ such that the set $V_1'$ of vertices $v$ without twins and of degree $d_1$ is of size     $|V_1'|> \beta |V_1|$. Suppose from now on that this is the case. Since $|V_1'| > |V_2|$ and there are no twins in $V_1'$, we see that $1<d_1<|V_2|-1$; we can suppose that $d_1\leq |V_2|/2$ by replacing $A$ with its complement, if necessary.

Let $\mathscr{H} = (V_2,\mathscr{A})$ be the hypergraph whose edges are the neighborhoods in $(V_1,V_2;A)$ of the vertices in $V_1'$. (They are all contained in $V_2$.) The hypergraph is $d_1$-uniform. Since there are no twins in $V_1'$, there are no identical edges. If $\mathscr{H}$ is the $d_1$-uniform complete hypergraph, then $V_1'$ can be identified with the Johnson scheme $\mathscr{S}_{d_1}(V_2)$. {\em
  (Scoppia in un pianto angoscioso e abbraccia la
  testa di Johnson.)}

Suppose then that $\mathscr{H}$ is not complete. We want to have a canonical coloring on $V_2^d$ for a $d\ll l$, $l=(\log |V_1'|)/\log |V_2|$, such that the elements of $V_2$ are not all twins. If $d_1\leq 6 \lceil l\rceil$, we define $d=d_1$ and color $\{(v_1,\dotsc,v_d)\in V_2^d: \{v_1,\dotsc,v_d\}\in \mathscr{H}\}$ in scarlet, and the rest in gray.

Suppose, on the contrary, that $d_1> 6 \lceil l\rceil$. Let $d = 6 \lceil l \rceil$. We color $\vec{v}=(v_1,\dotsc,v_d)$ in gray if the $v_i$ are not all distinct; otherwise, we give to $\vec{v}$ the color
\begin{equation}\label{eq:juj}
  |\{H\in \mathscr{H}: \{v_1,\dotsc,v_d\} \subset H\}|.\end{equation}
This coloring operation can be done in time of order
\[|V_1|\cdot \binom{d_1}{d} \leq |V_1|\cdot |V_2|^d =
|V_1| \cdot |V_2|^{6 \left\lceil \frac{\log |V_1'|}{\log |V_2|}
  \right\rceil} = |V_1|\cdot |V_1'|^{O(1)} = |V_1|^{O(1)}.\]

If the tuples with distinct $v_1,\dotsc,v_d$ had all the same color $\lambda$, we would have a $d-(|V_2|,d_1,\lambda)$ design with $|V_1'|$ edges. Then, by Proposition~\ref{prop:RW}, $|V_1'|\geq \binom{|V_2|}{s}$ for $s =
3\lceil l\rceil$. Since $|V_2|\geq (6 \log |V_1'|)^{3/2}$ and we can suppose $|V_1'|$ bigger than a constant,
\[\binom{|V_2|}{s} \geq 
\left(\frac{|V_2|}{s}\right)^s > \left(\frac{|V_2|}{6 l}\right)^s > 
\left(|V_2|^{1/3}\right)^{3 l} = |V_2|^{\frac{\log |V_1'|}{\log |V_2|}} = |V_1'|,\]
which gives a contradiction.

Therefore, for $d_1$ arbitrary, the tuples with distinct $v_1,\dotsc,v_d$ do
not have all the same color; in other words, the elements of $V_2$ are not all
twins with respect to our new $d$-ary structure. If there is a twin class $S$
of size $> |V_2|/2$, then, by Exercise~\ref{ex:reduc2}, at least one of the two graphs $(V_1',S;A\cap (V_1'\times S))$, $(V_1',V_2\setminus S; A\cap (V_1'\times (V_2\setminus S)))$ has no class of $> |V_1'|/2+1$ twins in $V_1$.
Since $2 |V_1'|/3\geq |V_1'|/2+1$ for $|V_1|>6$,
we apply Proposition~\ref{prop:bicoup} itself to one of these two graphs (say, the one on $V_1'\times S$ if the two are acceptable), and we are done, by induction. (It could happen that the size of $V_2$ has descended only to $|V_2|-1$, but all our choices were canonical -- non-choices, if you will -- therefore cost-free. We have lost nothing but time; to be precise, time $|V_1|^{O(1)}$, which is acceptable.)

Then, we have a coloring of $V_2^d$ in relation to which there is no twin class in  $V_2$ of size $> |V_2|/2$. We apply the functors $F_1$, $F_2$, $F_3$ (Weisfeiler-Leman) to this coloring. Then we use the Design Lemma (Prop.~\ref{prop:design}) with $\alpha=2/3$.
We find the elements $x_1,\dotsc,x_\ell\in V_2$ ($\ell=d-1$ or $\ell=d-2$) in the statement of Proposition~\ref{prop:design} by brute force, in time proportional to $|V_2|^d = |V_1|^{O(1)}$. We fix them, and we impose that $H$ fixes $x_1,\dotsc,x_\ell$, which has cost $|V_1|^{O(1)}$, in the sense that  \[\lbrack G: G_{(x_1,\dotsc,x_\ell)}\rbrack \leq |V_2|^d = |V_1|^{O(1)}.\]

If we are in the first case of the Design Lemma (no dominant color), we pick
color classes, starting from the reddest (interpret the quantity in (\ref{eq:juj}) as a wavelength), until we have a union of classes $S\subset V_2$ with $|V_2|/3 <|S|\leq (2/3) |V_2|$. (This works if there is no class of size $> |V_2|/3$;
if such classes exist, we simply let $S$ be the redder one.)
We apply Exercise~\ref{ex:reduc2}, and we obtain a graph $(V_1',V_2',A\cap (V_1'\times V_2'))$ satisfying the conditions of Proposition~\ref{prop:bicoup}
with $V_2' = S$ or $V_2' = V_2\setminus S$, and so $|V_2'|\leq \alpha |V_2|$.
Then, we apply Proposition~\ref{prop:bicoup} to this graph; the recursion works. (It is important here that $|V_2'|\leq \alpha |V_2|$, since we have already incurred a considerable cost ($|V_1|^{O(1)}$) in the index.)

Let us examine then the second case of the Design Lemma: we have a coloring of $V_2$ with a color class $C\subset V_2$ such that $|C|\geq (2/3) |V_2|$, and a homogeneous classical coherent configuration $\mathfrak{Y}$ non-trivial on $C$.

We define a graph with vertices $V_1'\cup V_2$, where $V_1'$ has, say,
the color of mother-of-pearl, and we have colored $V_2$ using the Design Lemma;
the edges will be not only those in $A\cap (V_1'\times V_2)$ (colored in black) but also the edges between elements of $V_2$, in the colors given by $\mathfrak{Y}$. We apply the refinements $F_1$, $F_2$ and $F_3$ (Weisfeiler-Leman) to this graph, and we obtain a coherent configuration $\mathfrak{X}$.


The configuration $\mathfrak{X}\lbrack V_2\rbrack$ is a refinement of $\mathfrak{Y}$. If there is no $\alpha$-dominant color, we reduce our problem to that for $(V_1',V_2',A\cap (V_1'\cap V_2'))$, $|V_2'|\leq \alpha |V_2|$, as before;
we can apply Proposition~\ref{prop:bicoup} to such a graph without
changing $\beta$ because \[|V_2'|\leq \frac{2}{3} |V_2|\leq \beta |V_2| < \beta |V_1'|.\]
Recursion works here also because of $|V_2'|\leq 2 |V_2|/3$: it is important
for $V_2'$ to be smaller than $V_2$ by a constant factor,
since we have already incurred a considerable cost ($|V_1|^{O(1)}$) in the index.

Suppose, on the other hand,
that there is a class of $(2/3)$-dominant colors $C_2$ in
$\mathfrak{X}\lbrack V_2\rbrack$. It must be a subset of $C$ because
$2/3+2/3>1$. We know that $\mathfrak{X}\lbrack C_2\rbrack$
is not a clique: if it were, the restriction of $\mathfrak{Y}$ to $C_2$ would
also have been, and that is impossible by Exercise~\ref{ex:noclique}.

We can suppose that there exists a color class $C_1\subset V_1'$ in $\mathfrak{X}_1$ satisfying $|C_1|> \beta |V_1|$; if not, we have a $\beta$-coloring of $V_1$, and we are done. Note that $|C_1| > \beta |V_1|$ implies
$|C_1|> |V_2|\geq |C_2|$.

We can also assume that not all edges of $\mathfrak{X}$ in
$C_1\times C_2$ are of the same color. If it were,
there would be a class of $\geq |C_1|> \beta |V_1| > |V_1|/2+1$ twins in $V_1$ in the graph $(V_1,C_2;A\cap (V_1\times C_2))$, of which $\mathfrak{X}\lbrack V_1\times C_2\rbrack$ is a refinement. In this case, by Exercise~\ref{ex:reduc2}, we would have a reduction to $(V_1,V_2\setminus C_2;A\cap (V_1\times (V_2\setminus C_2)))$, and we could finish by using Prop.~\ref{prop:bicoup} recursively.

We have thus reduced matters to Proposition \ref{prop:cohcoup}: we apply it
to $\mathfrak{X}\lbrack C_1\cup C_2\rbrack$. We obtain either
a colored $(1/2)$-partition of $C_1$, or a bipartite graph
$(W_1,W_2;A')$, $W_1\subset C_1$, $W_2\subset C_2$, with
$|W_1|\geq |C_1|/2$,
$|W_2|\leq |C_2|/2\leq |V_2|/2$, and such that no twin class in $W_1$
has more than $|W_1|/2$ elements.
We can assume that $|W_1|> \beta |V_1|$, as otherwise we have obtained
a colored $\beta$-partition of $V_1$. Hence $|W_2|< |W_1|/2$.
We can thus recur: we apply Proposition \ref{prop:bicoup} with
$(W_1,W_2;A')$ instead of $(V_1,V_2,A)$.

Recursion works in $O(\log |V_2|)$ steps
because $|W_2|\leq |V_2|/2$. As above -- if, in some iteration, the size of
$W_1$ (or $V_1$) goes under $\beta |V_1|$ (for the original value of
$|V_1|$), then we have obtained a colored $\beta$-partition of $V_1$.
\qed
\\

As we have seen, Bipartite Split-or-Johnson uses Coherent Split-or-Johnson,
which, in turn, will reduce matters to Bipartite Split-or-Johnson
for a bipartite graph $(V_1,V_2;A)$ with $V_2$ at most half of the size
of the original $V_2$.

\begin{prop}[Coherent Split-or-Johnson]\label{prop:cohcoup}
  Let $\mathfrak{X}=(C_1\cup C_2,c)$ be a coherent configuration with
  vertex color classes $C_1$, $C_2$, where $|C_1|>|C_2|$.
  Assume that neither $c|_{C_1\times C_2}$ nor $c|_{C_2\times C_2}$ is a constant
  function.
      
  Then, we can find, in time  $|C_1|^{O(1)}$, either
  \begin{itemize}
  \item a colored $(1/2)$-partition of $C_1$, or
  \item a bipartite graph $(V_1,V_2;A)$, $V_i\subset C_i$,
$|V_1|\geq |C_1|/2$, $|V_2|\leq |C_2|/2$, such that no twin class in
    $V_1$ has more than $|V_1|/2$ elements,
  \end{itemize}
  and an element $y\in C_2$
  such that the colored partition, or the bipartite graph,
  is canonical with respect to $G_y$, where $G=\Sym(C_1)\times \Sym(C_2)$.
\end{prop}
Of course, saying that $c|_{C_2\times C_2}$ is constant is
the same as saying that 
$\mathfrak{X}\lbrack C_2\rbrack$ is a clique.\\

\noindent{\sc Proof } ---
If the restriction $\mathfrak{X}\lbrack C_1\rbrack$ were a clique then, by coherence, for any color in $C_1\times C_2$ -- purple, say -- the neighborhoods in $(C_1,C_2;\mathscr{G}_{\text{purple}})$ of the vertices in $C_2$ would give us a balanced block design on $C_1$ (possibly degenerate). The block
design is incomplete because $c$ is not monochromatic on $C_1\times C_2$.
Fisher's inequality gives us $|C_2|\geq |C_1|$, contradicting our assumptions.
Therefore, $\mathfrak{X}\lbrack C_1\rbrack$ is not a clique.

If $\mathfrak{X}\lbrack C_1\rbrack$ is not primitive, the reddest of its non-connected relations gives us a canonical colored
$(1/2)$-partition of $V_1$, by Exercise~\ref{ex:samesize}. Hence,
we may assume that $\mathfrak{X}\lbrack C_1\rbrack$ is primitive.

We have two cases to consider: $\mathfrak{X}\lbrack C_2\rbrack$
primitive and $\mathfrak{X}\lbrack C_2\rbrack$ imprimitive.

Suppose first that $\mathfrak{X}\lbrack C_2\rbrack$ is imprimitive. The reddest
non-connected relation in $\mathfrak{X}\lbrack C_2\rbrack$ gives us a partition of $C_2$ into sets $B_1,\dotsc,B_m$, $m\geq 2$, all of the same size $\geq 2$. We will use
this partition, either to induce a partition in $C_1$, or to reduce $|C_2|$ by a constant factor. The first step consists in showing there are no twins in
$C_1$.

Since our configuration is coherent, the color of an edge in $C_1$ knows whether its vertices are twins with respect to $C_2$ (Ex.~\ref{ex:twindet}). Thus,
if there were twins in $C_1$ with respect to $C_2$, either one of the edge
colors in $C_1$ would give a non-connected relation -- and that contradicts the fact that
$\mathfrak{X}\lbrack C_1\rbrack$ is uniprimitive -- or all the elements of $C_1$ would be twins with
respect to $C_2$. In the last case, by Exercise~\ref{ex:gita},
$c|_{C_1\times C_2}$ would be monochromatic, and that is not the case. In conclusion,
there are no twins in $C_1$ with respect to $C_2$.

Our intention is to apply Exercise~\ref{ex:biparcomp} to obtain a contracted bipartite graph $C_1\times \{1,2,\dotsc,m\}$, with $m\leq |C_2|/2$.
We must just be careful that this graph does not end up being trivial.

Let $d_k$ be the degree of all $w\in C_2$ in the bipartite graph
$(C_1,C_2;\mathscr{G}_k)$ for a given color $k$, where $\mathscr{G}_k$ consists of the edges of color $k$.
(By Ex.~\ref{ex:biparhom}\ref{it:bar1}, the degree $d_k$
does not depend on $w$.)
If $d_k\leq |C_1|/2$ for all $k$, we fix a $w\in C_2$ (non-canonically) and
obtain a $(1/2)$-coloring of $C_1$: assign the color $c(x,w)$ to the vertex
$x\in C_1$.
Suppose then that there is a color -- which we will call {\em violet} -- such that $d_{\text{violet}}> |C_1|/2$. If there is a $1\leq i\leq m$ such that there is no class of more than $|C_1|/2$ twins in $C_1$ with respect to $B_i$,
we fix some element $y\in B_i$ of some such $i$ (non-canonically),
thus fixing that $i$. We thus obtain a reduction to a bipartite graph
$(C_1,B_i;\mathscr{G}_{\text{violet}}\cap (C_1\times B_i))$.

Suppose that this is not the case.
Then, for each $i$, there exists a twin class  $T_i\subset C_1$ with respect to $B_i$ such that $|T_i|>|C_1|/2$. For each $w\in B_i$, the edges from $w$ to each $v\in T_i$ are of the same color; then, they must be violet. That is,
the edges from every $v\in T_i$ to every $w\in B_i$ are violet.
Let {\em green} be an edge color in $C_1\times C_2$ that is not violet. Then,
the graph $X=(C_1,\{1,\dotsc,m\};D)$ in Exercise~\ref{ex:biparcomp} is
nonempty; since $(v,i)$ is violet for all $v\in T_i$, $X$ is not complete
either. Since $X$ is semiregular, there is no twin class in $C_1$ with respect
to $\{1,\dotsc,m\}$ with $> |C_1|/2$ elements (Ex. \ref{ex:gita}).
We have thus reduced matters to a bipartite graph $X$ of the kind we wanted.

Consider now the case of $\mathfrak{X}\lbrack C_2\rbrack$
primitive\footnote{The problem in the original proof of Babai was at this precise point. What follows is an alternative argument proposed by him ({\em col rumore sordo di un galoppo}) when this article was in the process of being edited. It is more concise and elegant than the original argument, in addition to being correct. Before, the proof used a forked reduction -- that is, the
  statement we want to prove reduced to two or more instances of the
  same statement. This made
 the index $\lbrack G:H\rbrack$ grow catastrophically.}.
Fix an arbitrary $y\in C_2$ (non canonically). We can suppose that there is a color -- say, violet -- such that $d_{\text{violet}}> |C_1|/2$,
since, otherwise, the colors of the edges that connect the elements of $C_1$
with $y$ will give us a $(1/2)$-coloring of $C_1$.
Write $V_1 = L_{\text{violet}}(y) = \{x\in C_1: c(x,y) = \text{violet}\}$.
Then $|V_1|>|C_1|/2$. Let {\em blue} be an edge color in
$\mathfrak{X}\lbrack C_2\rbrack$
such that the degree of $\mathscr{G}_{\text{blue}}$ is (positive and) $< |C_2|/2$; such a color exists because $\mathfrak{X}\lbrack C_2\rbrack$ is not a clique.
(If there are several colors like this, we choose the bluest one of them.)
Then, $V_2 = L_{\text{blue}}(y)\subset C_2$ satisfies $1\leq |V_2|< |C_2|/2$.

The bipartite graph $(V_1,V_2;\mathscr{G}_{\text{violet}}\cap (V_1\times V_2))$ is
semiregular by Exercise~\ref{ex:biparhom}\ref{it:bara2}. It is nonempty
because, for all $u\in V_2$,
$|L_{\text{violet}}(u)|> |C_1|/2$, and so
\mbox{$L_{\text{violet}}(u)\cap V_1 \ne \emptyset$.}
If it were complete, we would have $V_1\subset L_\text{violet}(u)$ for all
$u\in V_2$; since $|V_1| = |L_{\text{violet}}(y)| = |L_{\text{violet}}(u)|$, this would imply that $V_1=L_{\text{violet}}(u)$. Now, that would mean that $y$ and $u$ are twins in the graph $(C_1,C_2;\mathscr{G}_{\text{violet}})$. By the same argument as before (based on Exercise~\ref{ex:twindet}), the primitivity of $\mathfrak{X}\lbrack C_2\rbrack$
and the fact that $c|_{C_1\times C_2}$ is non-constant
imply that there are no twins in $C_2$ with respect to the graph
$(C_1,C_2;\mathscr{G}_{\text{violet}})$. Therefore,
$(V_1,V_2;\mathscr{G}_{\text{violet}}\cap (V_1\times V_2))$ is not complete.
By Exercise~\ref{ex:gita}, we obtain that no twin class in
$(V_1,V_2;\mathscr{G}_{\text{violet}}\cap (V_1\times V_2))$ has more than $|V_2|/2$
elements. We are thus done.
\qed

\subsection{Recursion and reduction}\label{subs:recred}

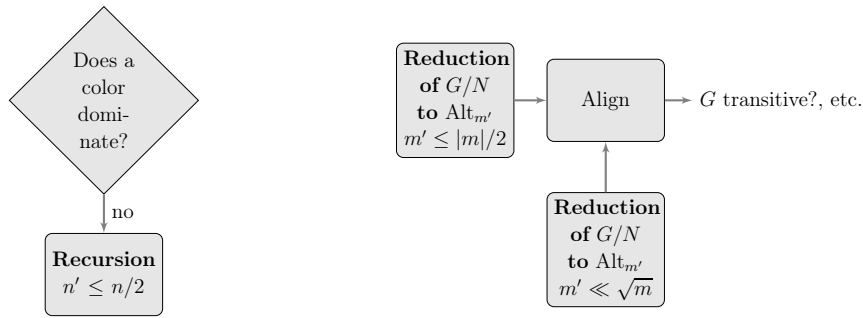
\begin{center}
      \scalebox{0.66}{
\begin{tikzpicture}[scale=2, node distance = 2cm, auto]
    \node [below of=init] (trans) {$G$ transitive?, etc.};
    \node [block, left of = trans, node distance = 3.5cm] (align) {Align};
    \node [block, left of = align, node distance = 3cm] (redhalf)
                {\bf Reduction of $G/N$ to $\Alt_{m'}$\\ $m'\leq |m|/2$};
    \node [block, below of = align, node distance = 3cm] (redsqrt)
          {\bf Reduction of $G/N$ to $\Alt_{m'}$\\ $m'\ll \sqrt{m}$};
    \node [decision, left of = redhalf, node distance = 7cm] (coldom) {Does a color dominate?};
    \node [block, below of = coldom, node distance = 3.5cm] (luks) {\bf Recursion $n'\leq n/2$};       
      \path [line] (align) -- (trans);
      \path [line] (redhalf) -- (align);
      \path [line] (redsqrt) -- (align);
      \path [line] (coldom) -- node [, color = black] {no} (luks);
  \end{tikzpicture}
    }
\end{center}

\subsubsection{The case without dominant colors}\label{subs:cacoudo}

We are in the case where a coloring $c_{\mathfrak{X}}:\Gamma\to \mathscr{C}$ has no dominant color. Here $c_{\mathfrak{X}}$ is the image of a structure $\mathfrak{X}$ under a functor $F$ that commutes with the action of $H_{(x_1,\dotsc,x_\ell)}$, where $H=\Alt(\Gamma)$, $x_i\in \Gamma$. The fact that $c_\mathfrak{X}$ has no dominant color will help us find or rule out its possible isomorphisms in $H_{\vec{x}}
= H_{(x_1,\dotsc,x_\ell)}$. In order to find or rule out isomorphisms in the whole $H=\Alt(\Gamma)$, we only have to work with a set of representatives $\{\sigma_1,\dots,\sigma_s\}$, $s\leq |\Gamma|^\ell = m^\ell$, of cosets of $H_{\vec{x}}$ in $H$, and to take the union of $\Iso_{H_{\vec{x}}}\left(c_{\mathfrak{X}}, c_{\mathfrak{Y}_i}\right)$ for $\mathfrak{Y}_i = \mathfrak{Y}^{\sigma_i^{-1}}$:
\begin{equation}\label{eq:kozron}
  \Iso_{H}(\mathfrak{X},\mathfrak{Y}) =
  \bigcup_{1\leq i\leq s} \Iso_{H_{\vec{x}}}\left(
  \mathfrak{X},\mathfrak{Y}_i\right)
\sigma_i,\;\;\;\;\;\;\;\;\;\;\;\;
\Iso_{H_{\vec{x}}}(\mathfrak{X},\mathfrak{Y}_i) \subset
\Iso_{H_{\vec{x}}}(c_{\mathfrak{X}},c_{\mathfrak{Y}_i}). 
\end{equation}
This is similar to equation (\ref{eq:rulu}), in \S \ref{subs:luks}. The cost of the procedure is multiplied by $s \leq m^\ell$.

If the coloring $c_{\mathfrak{X}}$ is not a permutation (in $H_{\vec{x}}$) of the coloring $c_{\mathfrak{Y}_i}$, then $\Iso_{H_{\vec{x}}}(c_{\mathfrak{X}},c_{\mathfrak{Y}_i}) = \emptyset$. Suppose, on the contrary, that there is at least one $\tau_i\in H_{\vec{x}}$ such that $c_{\mathfrak{X}} = c_{\mathfrak{Y}_i}^{\tau_i}$. (We say that $\tau_i$ {\em aligns} $c_{\mathfrak{X}}$ and $c_{\mathfrak{Y}_i}$.) It is trivial to find $\tau_i$. Now,
    \[
    \Iso_{H_{\vec{x}}}\left(\mathfrak{X},\mathfrak{Y}_i\right) =
    \Iso_{H_{\vec{x}}}\left(\mathfrak{X},\mathfrak{Y}_i^{\tau_i}\right) \tau_i^{-1}
    \subset \Aut_{H_{\vec{x}}}(c_{\mathfrak{X}}) \tau_i^{-1}. 
    \]
    Since $c_{\mathfrak{X}}$ has no dominant color, this is rather constraining, and
    that is what we want.

Let us apply this general procedure to the case of $G$ primitive, which we are about to discuss. There is a bijection $\iota:\Omega\to \{S\subset \Gamma: |S|=k\}$; therefore, $c_{\mathfrak{X}}$ induces a coloring $c':\Omega \to \{(k_i)_{i\in \mathscr{C}}: k_i\geq 0, \sum_i k_i =k\}$. We are in a situation similar to that at the end of \S \ref{sec:grasym}, but better: it is easy to show that, since no color class of $c_{\mathfrak{X}}$ has more than $\alpha |\Gamma|$ elements, no color class of $c'$ has more than $\alpha |\Omega|$ elements.

We proceed then as in the intransitive case in the proof of Luks (Thm.~\ref{thm:luxor}). In this way, the problem reduces to $\leq n$ problems of string isomorphism for strings of length $\leq \alpha n$ and total length $\leq n$. The last step ({\em lifting}) consists in finding elements of $G$ inducing $\tau_i$. Given a bijection $\iota$, this is trivial.

\subsubsection{The case of a colored partition}

Consider now a colored $\alpha$-partition (end of \S \ref{subs:secf1}) of a set of vertices $\Gamma$. This colored partition will be given canonically, i.e., as the image of a structure $\mathfrak{X}$ under a functor, like the coloring in \S \ref{subs:cacoudo}. We can suppose that the colored partition has a dominant color class $C$ ($|C| > \alpha |\Gamma|$, $\alpha>1/2$), since otherwise we can pass to \S \ref{subs:cacoudo}.

We want to know which elements of $\Alt(\Gamma)$ respect the colored $\alpha$-partition; this will help reduce the isomorphisms of $\mathfrak{X}$, like in (\ref{eq:kozron}). By definition of colored $\alpha$-partition, $C$ is partitioned into $\ell\geq 2$ sets of the same size $\geq 2$. The only allowed permutations in $\Alt_{m_0}$, $m_0 = |C|$, are those that respect the partition. The group that respects the partition is isomorphic to $\Alt_{m_0/\ell}$.

We have then reduced our problem to a problem with $m'=m_0/\ell\leq m/2$. After having solved this problem, we work -- as in  \S \ref{subs:cacoudo} -- on the other color classes.

Given two colored $\alpha$-partitions, we verify whether they have the same number of sets of the same size for each color, then we align the two partitions, and we proceed exactly as in the automorphism problem.

\subsubsection{The case of a Johnson scheme}
Let us be given a Johnson scheme on a set of vertices $\Gamma$, or rather two Johnson schemes $\mathscr{J}(m_i,k_i)$, $2\leq k_i\leq m_i/2$, on two sets of vertices $\Gamma_1$, $\Gamma_2$ of the same size. We have seen in \S \ref{subs:idgroup} now to identify $\Gamma_i$ (there, $\Omega$) explicitly with the subsets of size $k_i$ of a set $\Lambda_i$ (there, $\Gamma$) of size $m_i$. If $k_1\ne k_2$ and $m_1\ne m_2$, our structures are not isomorphic. If $k_1=k_2$ and $m_1=m_2$, we establish a bijection between $\Lambda_1$ and $\Lambda_2$ and we align the two structures. We have reduced our problem to a problem with $m'\ll \sqrt{m}$ instead of $m$.

The situation is then more favorable than in the case of a colored partition. Again, we leave the calculations to the reader.

\begin{center}
  * * *
\end{center}

A little confession: the case of $G$ primitive, which we have finished treating, could be treated exactly as in the case of $G$ imprimitive, which we will now examine. The motivation for the separate treatment for $G$ primitive is pedagogical. No effort has been in vain, since all the techniques that we have studied will be essential in the imprimitive case.

\section{The imprimitive case}\label{sec:casimp}

We have an explicit surjective function \[\phi:G\to \Alt(\Gamma),\] where $G<\Sym(\Omega)$ is a permutation group, $|\Gamma|=m$, $|\Omega|=n$. We can suppose that $|\Gamma|\geq C \log n$, $C$ arbitrary. The function $\phi$ factors as follows
\[G\to G/N\to \Alt(\Gamma),\]
where $N$ is the stabilizer of a block system, and $G/N\to \Alt(\Gamma)$ is an isomorphism.

We must determine $\Iso_G(\mathbf{x},\mathbf{y})$, where $\mathbf{x}$, $\mathbf{y}$ are strings. We have already solved the case $N = \{e\}$.

We attack the problem locally: for $T\subset \Gamma$, we arrive to obtain a certificate, either of the fact that $\phi(\Aut_{G_T}(\mathbf{x}))|_T$ contains $\Alt(T)$ ({``certificate of fullness''}), or of the contrary. (Here $G_T$ is the group $\{g\in G: T^{\phi(g)} = T\}$.) We calculate all our certificates for $T$ of moderate size $k$. If the number of certificates of fullness is very large, we will have proved that $\phi(\Aut_G(\mathbf{x}))$ contains a large alternating group; what is left to do is a version of the procedure of \S \ref{sec:grasym} ({``pull-back''}).

In the other case, the certificates will form a $k$-ary structure whose symmetry is bounded. We will then be able to apply the Design Lemma, followed by Split-or-Johnson, as before. There are also other particular cases, but they will bring us to colored $\alpha$-partitions, $\alpha<1$, which is also good.

\subsection{The local certificates}\label{sec:certloc}
\subsubsection{Certificates of automorphisms}\label{subs:certaut}
A local certificate\footnote{Or {``local-global''}, in Babai's nomenclature. {``Global''} refers to $\Aut_{G_T}(\mathbf{x}) <\Sym(\Omega)$.} for $T\subset \Gamma$ is
\begin{itemize}
\item either a triplet $(\text{{``not full''}},W,M(T))$, where $W\subset \Omega$, $M(T)<\Sym(T)$, $M(T)\neq \Alt(T)$ (therefore {``not full''}) and $\phi\left(\Aut_{G_T}^W(\mathbf{x})\right)|_T < M(T)$,
\item or a pair $(\text{{``full''}},K(T))$, where $K(T)<\Aut_{G_T}(\mathbf{x})$, and $\phi(K(T))|_T = \Alt(T)$.
\end{itemize}
The local certificate depends on $\mathbf{x}$ canonically. It is clear that a certificate of fullness (or non-fullness) guarantees that $\phi(\Aut_{G_T}(\mathbf{x}))|_T$ is, or is not, equal to $\Alt(T)$.

If $T$ is given as an ordered tuple, its certificate depends on the order of $T$ only in the sense that it does not depend on it: the same group $\{(2 3), e\} < \Sym(\{1,2,3\})$ (say) has a different appearance if we see it from the point of view of the order  $(1,2,3)$ or of the order $(2,1,3)$.

We construct the certificate by iteration. At the beginning of each step, $W\subset \Omega$ and $A(W)$ is the group $\Aut_{G_T}^W(\mathbf{x})$; the window $W$ will be invariant under $A(W)$. At the very beginning of the procedure, $W=\emptyset$ and $A(W) = G_T$. (We can calculate $G_T$ as in Exercise~\ref{ex:fhl}\ref{it:richt} in time $|\Omega|^{O(k)}$, where $k = |T|$.) At each step, we add to $W$ all the elements {\em affected} by $A(W)$ (see \S \ref{subs:stab}), and
then we update $A(W)$ according to the new $W$. We stop if $\phi(A(W))|_T \ne \Alt(T)$ (non-fullness) or if $W$ does not grow anymore; in the latter case,
no element of $\Omega\setminus W$ is affected by $A(W)$.

It is clear that there will be $\leq |\Omega|$ iterations. At the end, in the case of non-fullness, we return $(\text{{``not full''}},W,\phi(A(W)))$; in the case of fullness, we return $\left(\text{{``full''}},A(W)_{(\Omega\setminus W)}\right)$. It is clear that the pointwise stabilizer $A(W)_{(\Omega\setminus W)}$ is contained not only in $\Aut_{G_T}^W(\mathbf{x})$, but also in $\Aut_{G_T}(\mathbf{x})$, since it fixes all the points of $\Omega\setminus W$. We know that $\phi\left(A(W)_{(\Omega\setminus W)}\right) = \Alt(T)$ by Proposition \ref{prop:atinl}\ref{it:unaffstab}, under the condition that $|T|\geq \max(8,2+\log_2 |\Omega|)$.

Verifying whether $\phi(A(W))|_T = \Alt(T)$ is easy: we only have to verify, using Schreier-Sims, whether two arbitrary generators of $\Alt(T)$ are in $\phi(A(W))|_T$. Similarly, it is easy to determine which elements are affected by $A(W)$: we calculate $A(W)_x$ for each $x\in \Omega$ (by Schreier-Sims) and, always by Schreier-Sims, we verify whether $\phi(A(W)_x)|_T = \Alt(T)$. This takes polynomial time in $|\Omega|$.

It remains to see how to update $A(W)$, given $A\left(W^-\right)$, where we write $W^-$ for the old value of $W$. All elements of $A(W)$ are inside $A(W^-)$, and so $A\left(W\right) = \Aut_{A\left(W^-\right)}^{W}(\mathbf{x})$. As in equation (\ref{eq:rulu}),
  \begin{equation}\label{eq:ruyur}
    \Aut_{A\left(W^-\right)}^{W}(\mathbf{x}) =
  \bigcup_\sigma \Aut^{W}_{N\sigma}(\mathbf{x}) =
  \bigcup_\sigma \Iso_N^{W}\left(\mathbf{x},\mathbf{x}^{\sigma^{-1}}\right),
  \end{equation}
where $N$ is the kernel of $\phi|_{A(W^-)}$ and $\sigma$ runs through the representatives of the $k!/2$ cosets of $N$ in $A(W)$. We can rapidly find a $\sigma\in A(W^-) \cap \phi^{-1}(\{\tau\})$ for all $\tau\in \Sym(\Gamma)$, by Schreier-Sims.

Proposition \ref{prop:atinl}\ref{it:afforb} gives us that all orbits of $N$ contained in $W$ (the set of elements affected by $A(W^-)$) are of length $\leq |W|/k \leq |\Omega|/k$. Consequently, by rule (\ref{eq:udu3}), updating $A(W)$ reduces to  $k\cdot (k!/2)$ problems  of determining $\Iso$ for strings of length $\leq |\Omega|/k$.

Since the number of iterations is $\leq |\Omega|$, the procedure calls String-Isomorphism $\leq |\Omega| k \cdot (k!/2)$ times for strings of length $\leq |\Omega|/k$. This is acceptable for $k\ll (\log |\Omega|)^\kappa$.
For that matter, the
routine that took time $|\Omega|^{O(k)}$ is also acceptable under the same
condition. We will choose $\kappa=1$.

\subsubsection{Comparing certificates}
A little modification of the procedure above allows us to clarify the relation between two local certificates for two strings. Let $\mathbf{x}, \mathbf{x}':\Omega\to \Sigma$, $T,T'\subset\Gamma$, $|T|=|T'|=k$. For $S\supset T$, let $\mathbf{x}^S$ be the string
  \[\mathbf{x}^S(i) = \begin{cases} \mathbf{x}(i) & \text{if $i\in S$,}\\
    \text{glaucous} & \text{if $i\notin S$}\end{cases}\]
  where $\text{glaucous} \notin \Sigma$. We want to calculate \begin{equation}\label{eq:headtail}
    \Iso_{G_T \cdot \tau_{(T,T')}}\left(\mathbf{x}^W,\mathbf{x}^{W'}\right),\end{equation}
    where $G_T \cdot \tau_{(T,T')}$ is the set of elements of $G$ sending the set $T$ to $T'$, and $W'$ is the value of $W$ returned when the input is $T'$ instead of $T$.

To determine (\ref{eq:headtail}), we follow the procedure (\S \ref{subs:certaut}), modified as follows: we update, in each iteration, not only $A(W)$, but also the set $A(W) \tau$ of isomorphisms in $G_T \cdot \tau_{(T,T')}$ from $\mathbf{x}^W$ to $(\mathbf{x}')^{W'}$. Here is how to do it, in an analogous fashion as in (\ref{eq:ruyur}):
    \begin{equation}\label{eq:rurik}
    \bigcup_\sigma \Iso_{N\sigma}\left(\mathbf{x}^{W},
    \left(\mathbf{x}'\right)^{W'}\right) =
    \bigcup_\sigma \Iso_N\left(\mathbf{x}^{W},
    \left(\left(\mathbf{x}'\right)^{W'}\right)^{\sigma^{-1}}\right),\end{equation}
where $N$ is the kernel of $\phi|_{A(W^-)}$ and $\sigma$ runs through the representatives of the $k!/2$ cosets of $N$ contained in $A(W^-) \tau^-$. Since $W$ is stabilized by $A(W^-)$ (and then by $N$), the fact that $\sigma$ sends $W$ to $W'$ or not depends only on the coset of $N$ to which $\sigma$ belongs. (The set $\Iso$ in the last expression of (\ref{eq:rurik}) is empty if $W^\sigma \ne W'$.)

As before, all orbits of $N$ contained in $W$ are of length $\leq |W|/k$, and the problem reduces to $k\cdot (k!/2)$ calls by iteration to String-Isomorphism for strings of length $\leq |W|/k\leq |\Omega|/k$.

Furthermore, if $T$ and $T'$ are given as ordered tuples $(T)$, $(T')$, it is easy to determine
      \begin{equation}\label{eq:isoma}I(\mathbf{x},\mathbf{x}',T,T') =
        \Iso_{G_{(T)} \cdot \tau_{((T),(T'))}}\left(\mathbf{x}^W,
      \left(\mathbf{x}'\right)^{W'}\right),\end{equation}
      where $G_{(T)} \cdot \tau_{((T),(T'))}$ is the set of elements of $G$ sending
      the ordered tuple $(T)$ to $(T')$. In fact, we only have to determine (\ref{eq:headtail}), then use Schreier-Sims to detect the elements of (\ref{eq:headtail}) sending $(T)$ to $(T')$ in the right order.

\subsection{The aggregation of certificates}\label{subs:agrecert}

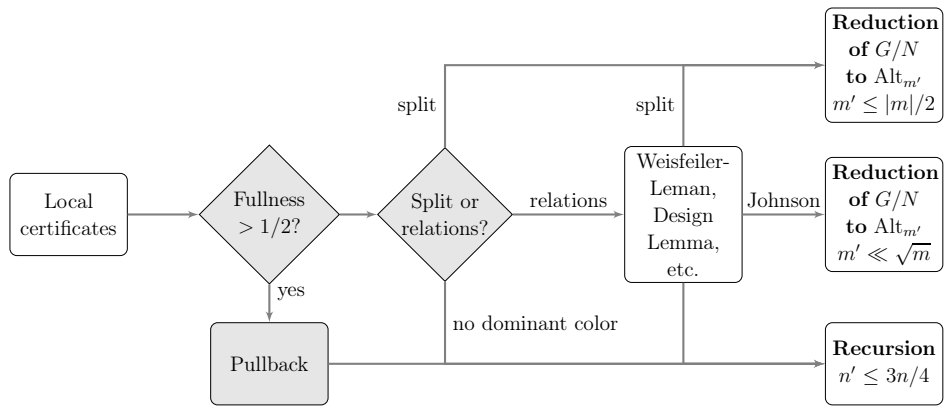
\begin{center}
\scalebox{0.66}{
  \begin{tikzpicture}[scale=2, node distance = 2cm, auto]
    \node [noblock] (cert) {Local certificates};
    \node [decision , right of = cert, node distance = 4cm] (certsym)
       {Fullness $>1/2$?};
    \node [block, below of = certsym, node distance = 3cm] (pull)
             {Pullback};
    \node [decision, right of = certsym, node distance = 3.5cm] (couprel)
           {Split or relations?};
    \node [noblock, right of = couprel, node distance = 4.75cm] (designs)
          {Weisfeiler-Leman, Design Lemma, etc.};
    \node [noblock, right of = designs, node distance = 4cm] (redsqrt) {\bf Reduction of $G/N$ to $\Alt_{m'}$\\ $m'\ll \sqrt{m}$};
              \node [noblock, below of = redsqrt, node distance = 3cm] (luks) {\bf Recursion
      $n'\leq 3n/4$};
   \node [noblock, above of = redsqrt, node distance = 3cm] (redhalf)
           {\bf Reduction of $G/N$ to $\Alt_{m'}$\\ $m'\leq |m|/2$};
    \path[line] (couprel) -- node [color=black] {relations} (designs);
    \path[line] (cert) -- (certsym);
    \path [line] (certsym) -- node [near start, color=black] {yes} (pull);
    \path[line] (pull) -- (luks);
    \path [line] (certsym) -- (couprel);
    \path [line] (couprel) |- node [near start, color=black] {split} (redhalf);
   \path [line] (designs) |- node [near start, color=black] {split} (redhalf);
   \path [line] (couprel) |- node [near start, color=black]
         {no dominant color} (luks);
         \path [line] (designs) -- node [color=black] {Johnson} (redsqrt);
            \path [line] (designs) |- (luks);

\end{tikzpicture}
  }
\end{center}

Following the procedure of \S \ref{subs:certaut} for a string $\mathbf{x}$, we find local certificates for each $T\subset \Gamma$ of size $k$, where $k$ is a constant $\sim C \log |\Omega|$ ($C>1/\log 2$) and $k<|\Gamma|/10$. Let $F<\Aut_G(\mathbf{x})$ be the group generated by the certificates of fullness $K(T)$. Let $S\subset \Gamma$ be the support of $\phi(F)$, i.e. the set of elements of $\Gamma$ that are not fixed by all the elements of $\phi(F)$.

Our objective is to determine the isomorphisms $\Iso_G(\mathbf{x}, \mathbf{x}')$ from $\mathbf{x}$ to another string~$\mathbf{x}'$. Since the certificates are canonical, the assignation of $F$ and $S$ to a string is so. Therefore, if we arrive to two different cases below while following the procedure for $\mathbf{x}$ and $\mathbf{x}'$, the two string are not isomorphic.

{\noindent Case 1: $|S|\geq |\Gamma|/2$, but no orbit of $\phi(F)$ is of length $>|\Gamma|/2$.}

Then, we color each element of $\Gamma$ according to the length of the orbit containing it. This is a canonical coloring. Either no color class is of size $>|\Gamma|/2$, or a color class of size $>|\Gamma|/2$ is cut into $\geq 2$ sets of the same size $\geq 2$. In one case or the other, we pass to a reduction/recursion.

{\noindent Case 2: $|S|\geq |\Gamma|/2$ and an orbit $\Phi$ of $\phi(F)$ is of length $>|\Gamma|/2$.}

{\em Case 2a: $\Alt(\Phi) < \phi(F)|_\Phi$.} We are in the case of large symmetry. We proceed as in \S \ref{sec:grasym}, up to the point where we must determine $\Iso_H(\mathbf{x},\mathbf{y})$ (where $\mathbf{y}$ is $\left(\mathbf{x}'\right)^{\sigma'}$, $\sigma'\in G$, and $H = \phi^{-1}\left(\Alt(\Gamma)_\Phi\right)$). Define $K = \phi^{-1}\left(\Alt(\Gamma)_{(\Phi)}\right)$ and let $\sigma_1,\sigma_2\in G$ be (arbitrary) preimages under $\phi$ of two generators of $\Alt(\Phi)<\Alt(G)$, found through Schreier-Sims. We know that the sets $\Aut_{K \sigma_i}(\mathbf{x})$, $i=1,2$, are nonempty, since $\Alt(\Phi)<\phi(F)|_\Phi$. Since $K$ has no orbit of length $>|\Omega|/2$, we can determine these two sets by calling String-Isomorphism for strings of length $\leq |\Omega|/2$ and total length $\leq 2 |\Omega|$. These generate $\Aut_H(\mathbf{x})$.

Again by the fact that $\Alt(\Phi)<\phi(F)|_\Phi$, the set $\Iso_H(\mathbf{x},\mathbf{y})$ is nonempty iff $\Iso_K(\mathbf{x},\mathbf{y})$ is nonempty. We can determine the latter set by calling String-Isomorphism as above, since $K$ has no orbits of length $>|\Omega|/2$. If it is not empty, we obtain the answer
\[\Iso_H(\mathbf{x},\mathbf{y}) = 
\Aut_H(\mathbf{x}) \Iso_K(\mathbf{x},\mathbf{y}).\]

{\em Case 2b: $\Alt(\Phi) \nless \phi(F)|_\Phi$.} Let $d\geq 1$ be the maximal integer with the property that $\phi(F)|_\Phi$ is $d$-transitive, i.e., $\phi(F)|_\Phi$ acts transitively on the set of $d$-tuples of distinct elements of $\Phi$. By CFSG, $d\leq 5$; if we do not want to use CFSG, we have the classical bound $d\ll \log |\Gamma|$.

Choose $x_1,\dotsc,x_{d-1}\in \Phi$ arbitrarily. The rest of our treatment of this case will then be canonical only with respect to \[G_{(x_1,\dotsc,x_{d-1})} = \{g\in G: x_i^{\phi(g)} = x_i \; \forall 1\leq i\leq d-1\}.\] This, as we know, is not a problem; see the beginning of \S \ref{subs:cacoudo}.

The restriction of the group $\phi(F)_{(x_1,\dotsc,x_{d-1})}$ to $\Phi' =
\Phi\setminus \{x_1,\dotsc,x_{d-1}\}$ is transitive on $\Phi'$, but it is not doubly transitive. Therefore, the Schurian coherent configuration that corresponds to it is not a clique. We give this configuration to Split-or-Johnson (\S \ref{sec:coujoh}), as at the end of \S \ref{sec:coujoh}.

To compare the configurations corresponding to two strings $\mathbf{x}$, $\mathbf{x}'$, we first align their orbits $\Phi$. (If they are not of the same size, or if a string falls into case 2a and the other does not, the strings are not isomorphic.) The isomorphisms will then be contained in the stabilizer $H<G$ of the set $\Phi$ (easy to determine, as towards the end of \S \ref{sec:grasym}, since $\phi$ is surjective). We can replace $\phi$ with the map $g\mapsto \phi(g)|_\Phi$ from $H$ to $\Alt(\Phi)$. Then we construct configurations as above, and we compare what Split-or-Johnson gives us.

Finally, we take care of the complement of $C$. As always, it involves calling String-Isomorphism for strings of length $\leq |\Omega|/2$ and total length $<|\Omega|$.

{\noindent Case 3: $|S|< |\Gamma|/2$.}
We start by aligning the supports $S$ for the strings $\mathbf{x}$, $\mathbf{x}'$, and by replacing $\phi$ with $g\mapsto \phi(g)|_{\Gamma \setminus S}$, as in case 2b.

We are going to define a $k$-ary relation with very few twins, in order to give it afterwards to the Design Lemma.

Let us look at the category of all strings $\Omega \to \Sigma$, where $\Omega$ and $\Sigma$ are fixed, an action of $G$ on $\Omega$ is given, and $\phi:G\to \mbox{Alt}(\Gamma)$ is also given. We have already been looking at it for quite a while, since we must compare colors in configurations induced by different strings to decide whether these strings are isomorphic.

This time, we define colors according to equivalence classes: two pairs $(\mathbf{x},(T))$, $(\mathbf{x}',(T'))$ ($T,T'\subset \Gamma\setminus S$, $|T|=|T'|=k$) are equivalent if the set of isomorphisms in (\ref{eq:isoma}) is nonempty. We color $(T)$ -- in the coloring of $(\Gamma\setminus S)^k$ corresponding to $\mathbf{x}$ -- according to the equivalence class of $(\mathbf{x},(T))$. Here, $(T)$ is an ordered $k$-tuple without repetitions; if $(T)$ has repetitions, it is colored in gray.

For a given $\mathbf{x}$, no twin class in $\Gamma$ can have $\geq k$ elements: if such a set existed with $\geq k$ elements, it would contain a subset $T$ with $k$ elements, and all orderings $(T)$ of $T$ would have the same color. This would mean that the set of isomorphisms in (\ref{eq:isoma}) would be nonempty for any two orderings $(T)$, $(T')$ of $T$. Consequently, $\Aut_{G_T}(\mathbf{x}^W)$ would contain elements giving all possible permutations of $T$. This would give us a contradiction, since $T$, being contained in $\Gamma\setminus S$, is not full.

Then, provided that $k\leq |\Gamma|/4$, we have a coloring of $(\Gamma\setminus S)^k$ with no twin class with $\geq |\Gamma\setminus S|$ elements. We will then be able to apply the Design Lemma, after applying the usual refinements $F_2$, $F_3$ (Weisfeiler-Leman).

But -- can we calculate these colorings? The equivalence classes are enormous. On the other hand, there is no need to calculate them. All we need to compare structures coming from strings $\mathbf{x}$, $\mathbf{y}$ is being able to compare two tuples $(T)$ (on the configuration given by $\mathbf{x}$ or $\mathbf{y}$) and $(T')$ (on the configuration given by $\mathbf{x}'=\mathbf{x}$ or $\mathbf{x}'=\mathbf{y}$) and say whether they are of the same color. In other terms, we will have to calculate -- at the beginning of the procedure, for all pairs $((T),(T'))$, $|T|=|T'|=k$, and for the pairs of strings $(\mathbf{x},\mathbf{x})$, $(\mathbf{x},\mathbf{y})$, $(\mathbf{y},\mathbf{y})$ -- the set of isomorphisms in (\ref{eq:isoma}). This we already know how to do. The colors then are, in practice, entries in an index that we enrich and to which we refer in the course of
our procedures.

We invoke then the Design Lemma, followed by Split-or-Johnson, and the rest of the procedure.

\begin{center}
  {--- \textsc{Fine dell'opera} ---}\\ \vskip 7pt
\end{center}
The reader can verify that what we have said so far (time taken by the procedures, type of recursion) is enough to give a bound of type $\exp((\log |\Omega|)^c)$ for the time of the algorithm that solves the string isomorphism problem. This gives a bound $\exp((\log n)^c)$ for the graph isomorphism problem with $n$ vertices. With a little bit of work it becomes clear that in both cases $c=3$. We give the details in the appendix. The exponent $c=3$ is smaller than the original one; it has become possible thanks to some improvements and simplifications that I have been able to contribute.

Acknowledgments .--- I warmly thank L. Babai, J. Bajpai, L. Bartholdi, D. Dona,
E. Kowalski,
W.~Kantor, G. Puccini, L. Pyber, A. Rimbaud, C. Roney-Dougal and W. Wordsworth\footnote{English version only.} for corrections and suggestions. In particular, L. Babai answered many of my questions, and gave me corrected or improved versions of several sections of \cite{Ba}. In particular, \S\S \ref{subs:secf1}--\ref{subs:concoh} and \S \ref{sec:design} are based on these new versions. I want also to thank V.~Ladret and V. Le Dret for many typographical and linguistic corrections.

\appendix

\section{Analysis of execution time}

\subsection{Some clarifications about the main procedure}

At every given moment, we work with a transitive group $G<\Sym(\Omega)$ acting on a system of disjoint blocks $\mathscr{B} = \{B_i\}$, $\Omega = \bigcup_i B_i$; we call $N$ the kernel of the action on $\mathscr{B}$. As a matter of fact, we will have a whole tower of systems of blocks $\mathscr{B}_1,\mathscr{B}_2,\dots, \mathscr{B}_k$, where $B_i$ is a refinement of $B_{i+1}$; $\mathscr{B}$ will mean $\mathscr{B}_k$, the coarsest system. At the beginning, there is but one system, $\mathscr{B}_1$, whose blocks $B_i$ are all of size $1$, and whose kernel $N$ is trivial.

We will want the action of $G$ on $\mathscr{B}$ to be primitive. So, if it is not, we add to the tower a minimal system $\mathscr{B}_{k+1}$ such that $\mathscr{B}_k$ is a refinement of $\mathscr{B}_{k+1}$. We redefine $\mathscr{B} = \mathscr{B}_{k+1}$; $N$ will be the kernel of the new $\mathscr{B}$.

If $G/N$ is small ($\leq b^{O(\log b)}$, where $b=|\mathscr{B}|$; case (a) of Theorem~\ref{thm:cam} (Cameron)), we reduce our problem to several instances of the problem with $N$ instead of $G$. Each of these instances decomposes into several instances -- one for each orbit of $N$. Each orbit $\Omega'$ of $N$ is contained in a block of $\mathscr{B}$. The intersections of $\Omega'$ with the blocks of $\mathscr{B}_1,\mathscr{B}_2,\dotsc$ give us a tower of systems of blocks for $N|_{\Omega'}$.

If we are in case (b) of Theorem~\ref{thm:cam}, we pass to $\leq b$ instances of the problem with $M\triangleleft G$ (where $\lbrack G:M\rbrack \leq b$) instead of $G$. We pass to a new system\footnote{This system could be equal to $\mathscr{B}$ only if $M=G$; see the second footnote in the statement of Theorem~3.1. In that case, the passage from $G$ to $M$ is obviously cost-free.} $\mathscr{B}'$ of $m' = \binom{m}{k} \leq b$ blocks, and we add it to the tower as its new last level. We call $N$ the kernel of the action of $M$ on $\mathscr{B}'$. Then, $M/N' = \Alt_m^{(k)}$. We replace $G$ with $M$ and redefine $\mathscr{B} = \mathscr{B}'$, $N=N'$.

Therefore, we have an isomorphism from $G/N$ to $\Alt_m$. We are in the main case that Babai attacks. His methods bring us to a reduction of $\Alt_m$, either to an intransitive group without large orbits, or to a product  $\Alt_{s_1}\wr \Alt_{s_2}$, $s_1,s_2>1$, $s_1 s_2\leq m$, or to a group $\Alt_{m'}$, $m'\ll \sqrt{m}$. (We are simplifying a little. We could have, say, a product $\Alt_{s_1}\wr \Alt_{s_2}$, acting on an orbit of large size $s_1 s_2\leq m$, and other groups on small orbits or several products acting on small orbits.)

In the intransitive case without large orbits, we proceed as in the proof by Luks. (The procedure would have been more costly than in Luks, but thanks to the lack of large orbits, the gain in the recursion is actually bigger.) In the case of $\Alt_{m'}$, $m'\ll \sqrt{m}$, we iterate the procedure. In the case of $\Alt_{s_1}\wr \Alt_{s_2}$ -- which corresponds to a partition into sets of size $r$ of the same color -- we have a primitive action of $\Alt_s$ on a system of $s$ blocks of size $r$. We pass then to this action and these blocks, without forgetting the blocks $\mathscr{B}'$, to which we will return later after having finished working on $\Alt_s$.

It is clear that this type of procedure reduces completely $\Alt_k$ in a number of iterations that is not bigger than $\log_2 m$.

\subsection{Recursion and time}

Let us examine the total execution time of the algorithm that finds the isomorphisms between two strings. The steps of the algorithm are not individually
onerous; none of them takes time more than $n^{O(\log n)}$.
Our attention must focus on the recursion.

In the general procedure, a recursion is always a descent either to shorter strings or to a group that is smaller or at least cut into finer pieces by a tower of systems of blocks having more levels. In the first type of descent the group stays the same, or rather is replaced by a restriction of itself. In the second case the length of the strings remains the same. (We can have also a mix of the two cases -- even better: the group becomes smaller and the strings also shorten.)

The least expensive and the least beneficial descent is that of the intransitive case of the procedure of Luks. It could happen that $G$ has two orbits on $\Omega$ ($|\Omega|=n$), one of length $n-1$ and one of length $1$. This is compatible even with a polynomial bound on the time, provided that the time taken before the descent is also polynomial: $n^{c+1} \geq (n-1)^{c+1} + 1^{c+1} + n^c$ for $c\geq 0$.

Other types of descent are more expensive, but also more beneficial: we descend to strings of length $\leq n/2$ (or $\leq 2 n/3$) or from $\Alt_m$ to $\Alt_{s_1}\wr \Alt_{s_2}$, $s_1 s_2 \leq m$, $s_1,s_2\leq m/2$, for example. It is clear that it is impossible to descend more than a logarithmic number of times in this way.

It is crucial to never forget that a (considerable) cost could be hidden in a loss of canonicity. If our choices are not canonical but with respect to a subgroup $H$ of our group $G$, the cost of their application will be multiplied by $\lbrack G:H\rbrack$. (Cf. \S \ref{subs:cacoudo}.)

\begin{center}
  * * *
\end{center}

Consider then the cost of each procedure. The intransitive case of Luks is, as we have already seen, compatible even with a polynomial bound. Let us concentrate then on the case where $G$ acts primitively on a system of blocks; let $N$ be the kernel.

If we are in case (\ref{it:rila}) of Theorem~\ref{thm:cam} or in case (\ref{it:rilb}) but with $m\leq C \log n$, we call $(m')^{O(\log n)}$ instances of the main procedure for strings of length $n/m'$ (where $m'\geq m$). This is consistent with a polynomial bound of the type $\exp(O((\log n)^c))$, $c\geq 2$. We can then focus on the case where there exists an isomorphism $\phi:G/N\to \Alt(\Gamma)$, $|\Gamma| = m > C\log n$. (The procedure of \S \ref{subs:idgroup} makes this isomorphism explicit.)

The first step to consider is the creation of local certificates, with the aim of creating a $k$-ary relation on $\Gamma$. (If $G$ is primitive, creating such a relation is trivial; see the beginning of \S \ref{subs:chasche}.) There are $n^k$ local certificates where $k = 2 \log n$ (say); we must calculate them and also compare each pair of certificates. Already the first step of the calculation of a certificate, i.e. the calculation of $G_T$, takes time $n^{O(k)}$ (more precisely $O((n/k)^{O(k)})$). Other calculations take less time. The use of recursion on the other hand is relatively burdensome: we call the main procedure $\leq n^2 \cdot k!$ times for strings of length $\leq n/k$. This happens for each set $T$ of size $k$, i.e. $\leq n^k/k!$ times. The procedure to compare pairs of certificates is analogous.

We call then the main procedure $O(n^{2 k+1})$ times for strings of length $\leq n/k$. In each of these calls, our tower of stabilizers is inherited: our group is transitive, equal to the restriction of $N^-$ to one of its orbits, where $N^-$ (called $N$ in \S \ref{sec:casimp}) is a subgroup of a subgroup $A(W^-)$ of $G$.

For two consecutive systems of blocks $\mathscr{B}_i$, $\mathscr{B}_{i+1}$, call $r_i$ the number of blocks of $\mathscr{B}_i$ inside each block of $\mathscr{B}_{i+1}$. It is clear that this number does not grow when we pass to the restriction of a subgroup of $G$ (for example, $N^-$) to one of its orbits.

Examine now the aggregation of local certificates (\S \ref{subs:agrecert}). There are three cases. In the first one the additional calculation time is almost trivial and we obtain a reduction, either to an intransitive group without large orbits, or to a product $\Alt_{s_1}\wr \Alt_{s_2}$ on a large orbit and possibly other groups on smaller orbits.

Here the analysis already becomes delicate. We must take into consideration not only the size of the domain but also the group that acts on it. More precisely, we must bound the number of times that our tower $\mathscr{B}_1, \mathscr{B}_2,\dotsc,\mathscr{B}_k$ could still be refined or shortened. This will be measured by
\[\rho = \sum_{1\leq i\leq k-1} (2 \lfloor \log_2 r_i\rfloor - 1),\]
where we suppose that we have removed repeated systems from the tower (therefore $r_i>1$).

Call $F(n,r)$ the execution time of the main procedure for strings of length $n$ and a tower of system of blocks for $G$ such that the parameter $\rho$ is $\leq r$. A reduction of $G/N$ makes $R$ decrease by at least $1$; a coloring without large color classes ensures a descent to strings of length $\leq n/2$. We must also include a factor of $\log n^{2k}$, taking into consideration the time required to access our comparisons of pairs of local certificates\footnote{Doing this type of comparison in advance helps us but not doing it in advance will not change asymptotically the order of the time we use.}. So, in the case we are examining, $F(n,r)$ is bounded by
\[n^{O(k)} + \left(n^{2 k +1} F(n/k,r) + F(n_1,r-1) + \sum_{i\geq 2} F(n_i,r)\right)\cdot O(k \log n),\]
where $\sum n_i = n$ and $n_i\leq n/2$ for $i\geq 2$, or
\[n^{O(k)} + \left(n^{2 k +1} F(n/k,r) + \sum_i F(n_i,r)\right) \cdot O(k \log n),
\]
where $\sum n_i = n$ and $n_i\leq n/2$ for $i\geq 1$. Since $k\ll \log n$ this is consistent with $F(n,r) = \exp\left(O\left(r+\log n\right)^{c}\right)$ for $c\geq 3$ or even with \[F(n,r) = \exp\left(O\left((\log n)^{c_1} + (\log r)^{c_2}\right)\right)\] for $c_1\geq 3$ and $c_2\geq 1$, for example.

Case 2a has a very similar cost, with a constant factor of difference. Cases 2b and 3 are different. In the two cases, we come to build a $d$-ary relation, with $d\leq 5$ in case 2b and $d=k\ll \log n$ in case 3. Then we call Weisfeiler-Leman followed by the Design Lemma for $d$-ary configurations and finally Split-or-Johnson.

Weisfeiler-Leman takes time $|\Gamma|^{O(d)} = m^{O(d)}$. The Design Lemma guarantees the existence of a tuple $(x_1,\dotsc,x_\ell)\in \Gamma$, $\ell\leq d-1$, with certain properties. We search for such tuple by brute force; this takes time $O(m^d)$. What is more important is that this choice is not canonical. Therefore, the execution time of all that remains is multiplied by $m^d = m^{O(\log n)}$.

Split-or-Johnson takes time $O(m^d)$. Here again we make choices that are not completely canonical; they impose a factor of $m^{O(\log m)}$ on all that follows. The result of Split-or-Johnson is either a colored $\beta$-partition, which gives
us a reduction to a product of the type $\Alt_{s_1}\wr \Alt_{s_2}$ and/or to shorter strings, or to a Johnson scheme, which gives us a reduction to $\Alt_{m'}$, $m'\ll \sqrt{m}$. Then, either
\begin{equation}\label{eq:urg1}
  F(n,r)\leq n^{O(k)} + O\left(k n^{2k+2} F(n/k,r)\right) +
  m^{O(\log n)} \left(1+ F(n_1,r-1) + \sum_{i\geq 2} F(n_i,r)\right),
  \end{equation}
where $\sum n_i = n$, and $n_i\leq n/2$ for $i\geq 2$, or
\begin{equation}\label{eq:urg2}
  F(n,r)\leq n^{O(k)} + O\left(k n^{2k+2} F(n/k,r)\right) +
m^{O(\log n)} \left(1 + \sum_i F(n_i,r)\right),\end{equation}
where $\sum n_i = n$ and $n_i\leq n/2$ for $i\geq 1$.

Here $m\leq n$. (We could work with a more accurate bound but this will be of little use.) So, the inequalities (\ref{eq:urg1}) and (\ref{eq:urg2}) are consistent with $F(n,r) = \exp\left(O\left(r+\log n\right)^{c}\right)$ for $c\geq 3$.

Since $r\leq 2 \log_2 n$, we conclude that the total execution time of the procedure to determine the isomorphisms between two strings of length $n$ is
\[F(n,r) \leq e^{O\left(\log n\right)^3}.\]

\section{Exercises and solutions}

\noindent{\sc Solution to Ex.~\ref{ex:fhl}} --- (a) The $C_{i}$ built by the Schreier-Sims algorithm are such that $|C_{i}|=|G_{i}/G_{i+1}|$, so that $|G|=\prod_{i=0}^{n-2}|C_{i}|$. Let $g\in \Sym(\Omega)$: we use the algorithm on $A$ (a generating set of $G$) and on $A\cup\{g\}$ (which generates $\langle G,g\rangle$) and we observe that $g\in G$ if and only if $|G|=|\langle G,g\rangle|$. The
point is of course that we can now calculate both $|G|$ and
$|\langle G,\rangle|$. Since Schreier-Sims takes time $O(n^{5})$, the same holds for deciding whether $g\in G$.

A more practical solution is to simply call the function \textsc{Filter} used in Algorithm 1 for this element $g$. The group $G$ contains $g$ if and only if $g$ can be written as a product $\gamma_{0}\gamma_{1}\ldots\gamma_{n-2}$ with $\gamma_{i}\in C_{i}$ (by our construction of the sets $C_{i}$): therefore, using \textsc{Filter} on $g$ we notice that $g\in G$ if and only if the \textbf{if} condition inside the function is always satisfied for the various intermediate elements $\gamma_{i}^{-1}\ldots\gamma_{1}^{-1}\gamma_{0}^{-1}g$; in other words, if \textsc{Filter} returns $(n-1,e)$ then $g\in G$, while if it returns $(i,\gamma)$ with $i<n-1$ then $g\not\in G$.

(\ref{it:richt}) We modify Algorithm 1: instead of using $G=G_0>G_1>G_2>\dotsc>G_{n-1}=\{e\}$ we use the sequence $G>H>H_1>H_2>\dotsc>H_{n-1}=\{e\}$, obtaining sets $C_{i}$ of representatives of cosets of $H_{i+1}$ inside $H_{i}$ for $0\leq i<n-1$ (and a set $C_{-1}$ of representatives of cosets of $H$ inside $G$); $\bigcup_{0\leq i<n-1}C_{i}$ describes $H$. In order to make the algorithm work, the test for $i=-1$ inside the function \textsc{Filter} should be $h^{-1}\gamma\in H$, which corresponds to the test $x_{i+1}^{h}=x_{i+1}^{\gamma}$ required for $0\leq i<n-1$ to decide whether $h^{-1}\gamma\in H_{i+1}$: by hypothesis, the test to decide whether $h^{-1}\gamma\in H$ takes polynomial time, and again by hypothesis every time we call \textsc{Filter} we have to perform the test at most $n^{O(1)}$ times since $|C_{-1}|\leq [G:H]$; therefore the modified algorithm is still polynomial-time.

(\ref{it:sshom}) The homomorphism $\phi$ is ``given'' in the sense that we are provided with pairs $(g,\phi(g))\in G\times\phi(G)$ where $g$ runs through a set of generators of $G$. First, we use Schreier-Sims on $G$ to obtain in polynomial time a set of generators of $G$ consisting of full sets of representatives $C_{i}$: these new generators $\gamma$ are built as products of the old ones, so we still have pairs $(\gamma,\phi(\gamma))$; therefore we can assume from now on that $\phi$ is described by the images of the generators $\gamma$ computed by Schreier-Sims.

As in (\ref{it:richt}), we now work with a modified sequence:
\begin{equation*}
G=G'_{0}>G'_{1}>G'_{2}>\dotsc>G'_{n'-1}=\ker\phi=K>K_{1}>\dotsc>K_{n-1}=\{e\}
\end{equation*}
Here $G'_{i}$ is defined as the preimage of the pointwise stabilizer $\Sym(\Omega')_{(x'_{1},\dotsc,x'_{i})}$, where $\Omega'=\{x'_{1},\dotsc,x'_{n'}\}$, and $K_{i}$ is the usual stabilizer $K_{(x_{1},\dotsc,x_{i})}$ of points in $\Omega$. Algorithm 1 with this sequence still takes polynomial time, following the same reasoning as in (\ref{it:richt}): each $\lbrack G'_{i}:G'_{i+1}\rbrack$ is of size $\ll n^{O(1)}$ thanks to our hypothesis $n'\ll n^{O(1)}$; moreover, we can write any element $g\in G$ as a product of generators $\gamma$ as above in polynomial time, so that we can quickly obtain $\phi(g)$ in order to test whether $g$ belongs to any given $G'_{i}$.

Now that we have an algorithm that uses this sequence, the task is easy to accomplish. For each $h$ among the generators of $H$, we can find an element inside $\phi^{-1}(h)$ by building a suitable product out of the sets of representatives of the $G'_{i}/G'_{i+1}$ so that the image of such a product sends each $x'_{i}$ to $h(x'_{i})$; on the other hand, the union of all the sets of representatives of the $K_{i}/K_{i+1}$ generates $K=\ker\phi$: the union of these and of the products built before for each $\phi^{-1}(h)$ generates $\phi^{-1}(H)$.
\bigskip

\noindent{\sc Solution to Ex.~\ref{ex:skeleconfig}} --- Let $\mathfrak{X} = (\Gamma,c)$ be a $k$-ary configuration and let $1\leq l\leq k$; for $\vec{x}=(x_{1},x_{2},...,x_{l})$ we define a $l$-ary partition structure $\mathfrak{X}^{(l)} = (\Gamma,c^{(l)})$ with $c^{(l)}(\vec{x})=c(x_{1},x_{2},...,x_{l},x_{l},...,x_{l})$.

If $c^{(l)}(\vec{x})=c^{(l)}(\vec{y})$, then $c(x_{1},x_{2},...,x_{l},x_{l},...,x_{l})=c(y_{1},y_{2},...,y_{l},y_{l},...,y_{l})$, which means that $\rho(x_{1},x_{2},...,x_{l},x_{l},...,x_{l})=\rho(y_{1},y_{2},...,y_{l},y_{l},...,y_{l})$ and then $\rho^{(l)}(\vec{x})=\rho^{(l)}(x_{1},x_{2},...,x_{l})=\rho^{(l)}(y_{1},y_{2},...,y_{l})=\rho^{(l)}(\vec{y})$, where $\rho^{(l)}$ is the equivalence relation map for $l$.

Now, let $\tau$ be a function from $\{1,...,l\}$ to itself, and define a function $\tau'$ as $\tau'(i)=\tau(\min\{i,l\})$ for all $1\leq i\leq k$; then:
\begin{eqnarray*}
c^{(l)}(\tau(\vec{x})) & = & c^{(l)}(x_{\tau(1)},x_{\tau(2)},...,x_{\tau(l)}) \\
& = & c(x_{\tau(1)},x_{\tau(2)},...,x_{\tau(l)},x_{\tau(l)},...,x_{\tau(l)})  \\
 & = & c(x_{\tau'(1)},x_{\tau'(2)},...,x_{\tau'(l)},x_{\tau'(l+1)},...,x_{\tau'(k)})  \\
 & = & \tau'^{\eta}(c(x_{1},x_{2},...,x_{l},x_{l+1},...,x_{k}))
\end{eqnarray*}
The map $\eta$ comes from the definition of configuration in the case of $\mathfrak{X}$. The elements $x_{l+1},...,x_{k}$ could be anything, since in the definition of $\tau'$ we have not used their indices, therefore we can choose $x_{l}=x_{l+1}=...=x_{k}$; thus we have $c^{(l)}(\tau(\vec{x}))=\tau'^{\eta}(c^{(l)}(\vec{x}))$, and there is a map $\eta^{(l)}$ defined as $\tau^{\eta^{(l)}}=\tau'^{\eta}$ such that $\mathfrak{X}^{(l)}$ satisfies the second condition for being a configuration.
\bigskip

\noindent{\sc Solution to Ex.~\ref{ex:skelecoh}} --- We already know that the skeleton is a configuration, so we have only to prove that it is coherent. By definition $c^{(l)}(\vec{x}^{i}(z))=c(\vec{x}^{i}(z),x_{l},...,x_{l})$ (we will write the vector $(x_{l},...,x_{l})$ of length $k-l$ as $\vec{x_{l}}$ for simplicity): since we are working with the $l$-ary configuration $\mathfrak{X}^{(l)}$ we care only about the colors obtained when substituting the first $l$ components, while any color will be valid for the remaining $k-l$; therefore for any $\vec{k}\in\mathscr{C}^{l}$ and any $j\in\mathscr{C}$ we have:
 \[\begin{aligned}
  &\gamma^{(l)}(\vec{k},j) = |\{z\in\Gamma: c^{(l)}(\vec{x}^{i}(z)) = k_i\; \forall 1\leq i\leq l\}|  \\
   &=  |\{z\in\Gamma: c(\vec{x}^{i}(z),\vec{x_{l}}) = k_i\; \forall 1\leq i\leq l\}|  \\
   &=  \sum_{\vec{h}} |\{z\in\Gamma: c(\vec{x}^{i}(z),\vec{x_{l}}) = k_i\; \forall 1\leq i\leq l,\ c(\vec{x},\vec{x_{l}}^{i}(z)) = h_i\; \forall 1\leq i\leq k-l\}|  \\
   &=  \sum_{\vec{h}} \gamma(\vec{k}\vec{h},j)
 \end{aligned}\]
where the sum is over all $\vec{h}\in\mathscr{C}^{k-l}$. The numbers $\gamma(\vec{k}\vec{h},j)$ are independent from $\vec{x}$, thus proving the independence of $\gamma^{(l)}(\vec{k},j)$ and the coherence of $\mathfrak{X}^{(l)}$.
\bigskip

\noindent{\sc Solution to Ex.~\ref{ex:rescoh}} --- We prove the following: if $\mathscr{C}'\subset\mathscr{C}$ is the subset of colors that appear in $\mathfrak{X}[\Gamma']$, then for every $\vec{k}\in\mathscr{C}'^{k}$ and every $j\in\mathscr{C}'$ we have an intersection number $\gamma'(\vec{k},j)$ for $\mathfrak{X}[\Gamma']$ that is equal to the intersection number $\gamma(\vec{k},j)$ for $\mathfrak{X}$.
 
 Let us be given a color $r$; since $\mathfrak{X}$ is a configuration, if $c(\vec{x})=r$ then $r$ knows the colors $c(\tau(\vec{x}))$ for all maps $\tau:\{1,...,k\}\rightarrow\{1,...,k\}$: in fact, if $\eta$ is the map in the definition of configuration then $c(\tau(\vec{x}))=\tau^{\eta}(r)$. In particular, if $\tau$ is one of the $k$ constant maps $\forall j(\tau(j)=i)$ we have that $r$ knows the colors $c(x_{i},x_{i},...,x_{i})$ for each $1\leq i\leq k$: this means that $r$ knows the color class to which each of the vertices that constitute a $k$-tuple of color $r$ belongs. But then, if we restrict ourselves to the color class $\Gamma'$ and we consider the substructure $\mathfrak{X}[\Gamma']$, all the colors $r$ that have vertices outside $\Gamma'$ will not appear at all in $\mathfrak{X}[\Gamma']$. Therefore the set $\mathscr{C}'$ of colors that appear in $\mathfrak{X}[\Gamma']$ is really the set of colors with vertices all in $\Gamma'$, so for any choice of $\vec{k}\in\mathscr{C}'^{k},j\in\mathscr{C}'$ we still have $\gamma'(\vec{k},j)=\gamma(\vec{k},j)$ since none of the possible $z\in\Gamma$ that contribute to $\gamma'$ lies outside $\Gamma'$.
\bigskip

\noindent{\sc Solution to Ex.~\ref{ex:veccoh}} --- (a) Let us define $\tau:\{1,...,k\}\rightarrow\{1,...,k\}$ as $\tau(i)=\min\{i+l,k\}$. Then for $\vec{y}\in\Gamma^{k-l}$, using the map $\eta$ in the definition of configuration, we have:
\begin{eqnarray*}
 c^{(k-l)}(\vec{y}) & = & c(y_{1},...,y_{k-l},y_{k-l},...,y_{k-l}) \ = \ c(\tau(x_{1},...,x_{l},y_{1},...,y_{k-l})) \ = \\
 & = & \tau^{\eta}(c(x_{1},...,x_{l},y_{1},...,y_{k-l})) \ = \ \tau^{\eta}(c_{\vec{x}}(\vec{y}))
\end{eqnarray*}
So the color $c_{\vec{x}}(\vec{y})$ knows $c^{(k-l)}(\vec{y})$, and in particular $c_{\vec{x}}(\vec{y})=c_{\vec{x}}(\vec{z}) \Rightarrow c^{(k-l)}(\vec{y})=c^{(k-l)}(\vec{z})$, which means that $c_{\vec{x}}$ is a refinement of $c^{(k-l)}$.

(b) It is easy to see that if $\mathfrak{X}$ is a configuration then $\mathfrak{X}_{\vec{x}}$ is also a configuration: $c_{\vec{x}}(\vec{y})=c_{\vec{x}}(\vec{z})$ means $c(\vec{x}\vec{y})=c(\vec{x}\vec{z})$, from which we have $\rho(\vec{x}\vec{y})=\rho(\vec{x}\vec{z})$ and restricting to the last $k-l$ components $\rho_{\vec{x}}(\vec{y})=\rho_{\vec{x}}(\vec{z})$; moreover, from $\tau:\{1,...,k-l\}\rightarrow\{1,...,k-l\}$ we define $\tau':\{1,...,k\}\rightarrow\{1,...,k\}$ by imposing $\tau'(i)=i$ for $i\leq l$ and $\tau'(i)=\tau(i-l)$ for $i>l$, so that:
\begin{equation*}
 c_{\vec{x}}(\tau(\vec{y}))=c(\vec{x}\tau(\vec{y}))=c(\tau'(\vec{x}\vec{y}))=\tau'^{\eta}(c(\vec{x}\vec{y}))=\tau'^{\eta}(c_{\vec{x}}(\vec{y}))
\end{equation*}
and $\mathfrak{X}_{\vec{x}}$ is a configuration once we define $\eta_{\vec{x}}$ by $\tau^{\eta_{\vec{x}}}=\tau'^{\eta}$.

Suppose now that $\mathfrak{X}$ is coherent. For any $\vec{k}\in\mathscr{C}^{k-l}$ and any $j\in\mathscr{C}$ we have:
 \begin{eqnarray*}
  \gamma_{\vec{x}}(\vec{k},j) & = & |\{z\in\Gamma: c_{\vec{x}}(\vec{y}^{\; i}(z)) = k_i\; \forall 1\leq i\leq k-l\}| \\
  & = & |\{z\in\Gamma: c(\vec{x}\vec{y}^{\; i}(z)) = k_i\; \forall 1\leq i\leq k-l\}| \\
  & = & \sum_{\vec{h}} |\{z\in\Gamma: c(\vec{x}^{\; i}(z)\vec{y}) = h_i\; \forall 1\leq i\leq l,\ c(\vec{x}\vec{y}^{\; i}(z)) = k_i\; \forall 1\leq i\leq k-l\}| \\
  & = & \sum_{\vec{h}} \gamma(\vec{h}\vec{k},j)
 \end{eqnarray*}
where the sum is over all $\vec{h}\in\mathscr{C}^{l}$. The numbers $\gamma(\vec{h}\vec{k},j)$ are independent from $\vec{y}$, thus proving the independence of $\gamma_{\vec{x}}(\vec{k},j)$ and the coherence of $\mathfrak{X}_{\vec{x}}$.
\bigskip

\noindent{\sc Solution to Ex.~\ref{ex:noclique}} --- Suppose that there exists a set $B$ as in the statement, and call white the color of the edges of the large clique; let black be another edge color of $\mathfrak{X}$ and let $\mathscr{G}=\mathscr{G}_{\text{black}}$: $\mathscr{G}$ is a biregular nonempty directed graph, so call $d$ the degree  of each vertex of $\mathscr{G}$ ($d$ is both their outdegree and their indegree). The induced subgraph of $\mathscr{G}$ on the subset $B$ is empty since all edges in $B$ were white in $\mathfrak{X}$, therefore for any $b\in B$ all edges $(b,x)$ of $\mathscr{G}$ have $x\in\Gamma\setminus B$; thanks to this fact we have:
\begin{eqnarray*}
 |B|d & = & \sum_{b\in B}|\{x\in \Gamma: (b,x)\mbox{ edge in }\mathscr{G}\}| \\ 
 &= & \sum_{b\in B}|\{x\in \Gamma\setminus B: (b,x)\mbox{ edge in }\mathscr{G}\}| \\ 
 & = & \sum_{x\in\Gamma\setminus B}|\{b\in B:(b,x)\mbox{ edge in }\mathscr{G}\}| \\ 
 &\leq& \sum_{x\in\Gamma\setminus B}|\{y\in\Gamma:(y,x)\mbox{ edge in }\mathscr{G}\}| \\
 & = & |\Gamma\setminus B|d
\end{eqnarray*}
which contradicts $|B|>|\Gamma|/2>|\Gamma\setminus B|$.
\bigskip

\noindent{\sc Solution to Ex.~\ref{ex:samesize}} --- (\ref{it:bibr1}) We follow \cite{Do}. We proceed by induction on $k$: for $k=2$ the statement is already satisfied by coherence. Suppose now that the statement holds for $k$. We are given $x_{0},x_{k+1}$ with $c(x_{0},x_{k+1})=r_{0}$ and we want to find the number of walks of colors $r_{1},r_{2},...,r_{k+1}$ from $x_{0}$ to $x_{k+1}$: such a walk would be just the composition of two walks, one of colors $r_{1},r_{2},...,r_{k-1}$ (from $x_{0}$ to $x_{k-1}$) and the other of colors $r_{k},r_{k+1}$ (from $x_{k-1}$ to $x_{k+1}$); therefore we can just consider any walk $r_{1},r_{2},...,r_{k-1},r'$ of length $k$ from $x_{0}$ to $x_{k+1}$ (in this case $c(x_{k-1},x_{k+1})=r'$) and for each of these walks we count all the possible $x_{k}$ with $c(x_{k-1},x_{k})=r_{k}$ and $c(x_{k},x_{k+1})=r_{k+1}$ given the $r'$-colored edge $(x_{k-1},x_{k+1})$.

By inductive hypothesis we have already constants $\gamma^{(k)}$ for walks of length $k$, and we have constants $\gamma^{(2)}$ for walks of length $2$ by coherence; thus, defining $\gamma^{(2)}(r_{0},r_{1},r_{2})=\gamma(r_{2},r_{1},r_{0})$, we have constants $\gamma^{(k+1)}$ as follows:
\begin{equation*}
\gamma^{(k+1)}(r_{0},r_{1},...,r_{k},r_{k+1})=\sum_{r'\in {\mathscr C}}\gamma^{(k)}(r_{0},r_{1},...,r_{k-1},r')\gamma^{(2)}(r',r_{k},r_{k+1})
\end{equation*}

(\ref{it:bibr2}) We call $b$ the unique vertex color of $(\Gamma,c)$; we recall also that by the definition of configuration we have $c(x,y)=c(y,x)^{-1}$ for any pair of vertices $x,y\in\Gamma$.

For any two edge colors $r,r'$, by (\ref{it:bibr1}) $r'$ knows whether there exists a finite walk entirely colored in $r$ between any two vertices $x,y$ with $c(x,y)=r'$. Therefore, having fixed a color $r$, we have that the pairs of points in the graph $\mathcal{G}_{r}$ that belong to the same connected component are exactly the edges of the configuration that have a color $r'$ for which such walks exist; calling $\mathscr{C}'\subset\mathscr{C}$ the subset of such colors $r'$, for any $x\in\Gamma$ the connected component $B(x)\subset\Gamma$ containing $x$ satisfies:
\begin{eqnarray*}
 |B(x)| & = & |\{y\in\Gamma: c(x,y)\in\mathscr{C}'\}| \\ 
 &= & \sum_{r'\in\mathscr{C}'}|\{y\in\Gamma: c(x,y)=r'\}| \\
 & = & \sum_{r'\in\mathscr{C}'}|\{y\in\Gamma: c(y,x)=r'^{-1},c(x,y)=r'\}| \\ 
 &= & \sum_{r'\in\mathscr{C}'}\gamma(r'^{-1},r',b)
\end{eqnarray*}
The constants $\gamma$ are independent from $x$, so the same is true for $|B(x)|$.
\bigskip

\noindent{\sc Solution to Ex.~\ref{ex:biparhom}} --- (a) If green is an edge color in $C_{1}\times C_{2}$, then all green edges are in $C_{1}\times C_{2}$; this comes from the fact that $\mathfrak{X}$ is a configuration, which means in particular that the color green knows both $c_{1}$ and $c_{2}$ (the vertex colors defining $C_{1}$ and $C_{2}$): if $c(x,y)=\text{green},c(x,x)=c_{1},c(y,y)=c_{2}$ then $c_{1}=\tau_{1}^{\eta}(\text{green}),c_{2}=\tau_{2}^{\eta}(\text{green})$, where $\eta$ is the map in the definition of configuration and $\tau_{1},\tau_{2}$ are the two constant maps $\tau_{1}(1)=\tau_{1}(2)=1$ and $\tau_{2}(1)=\tau_{2}(2)=2$. Another consequence of the definition of configuration is that $c(x,y)=\text{green}$ if and only if $c(y,x)=\text{green}^{-1}$.
 
 Thanks to these two properties we deduce:
 \begin{eqnarray*}
 d^{+}(v_{1}) & = & |\{v_{2}\in C_{2}: c(v_{1},v_{2})=\text{green}\}| \\
 & = & |\{v\in \Gamma: c(v_{1},v)=\text{green}\}|  \\
 & = & |\{v\in \Gamma: c(v,v_{1})=\text{green}^{-1},c(v_{1},v)=\text{green}\}| \\ 
 & = & \gamma(\text{green}^{-1},\text{green},c_{1})
 \end{eqnarray*}
 which shows that $d^{+}(v_{1})$ is the same for all $v_{1}$ with $c(v_{1},v_{1})=c_{1}$, i.e. for all $v_{1}\in C_{1}$.
 
 Analogously:
 \begin{eqnarray*}
 d^{-}(v_{2}) & = & |\{v_{1}\in C_{1}: c(v_{1},v_{2})=\text{green}\}| \\
 & = & |\{v\in \Gamma: c(v,v_{2})=\text{green}\}| \\
 & = & |\{v\in \Gamma: c(v,v_{2})=\text{green},c(v_{2},v)=\text{green}^{-1}\}| \\ 
 &= & \gamma(\text{green},\text{green}^{-1},c_{2})
 \end{eqnarray*}
 and again this is the same for all $v_{2}\in C_{2}$.
 
 (b) Let $x$ be any vertex in $L_{1}$, i.e. $c(x,y)=\text{aqua}$; then:
 \begin{eqnarray*}
 d^{+}(x) & = & |\{z\in L_{2}: c(x,z)=\text{cyan}\}| \\ 
 &= & |\{z\in \Gamma: c(z,y)=\text{beige},c(x,z)=\text{cyan}\}| \\ 
 & = & \gamma(\text{beige},\text{cyan},\text{aqua})
 \end{eqnarray*}
 
Let $z$ be any vertex in $L_{2}$, i.e. $c(z,y)=\text{beige}$; then:
 \begin{eqnarray*}
 d^{-}(z) & = & |\{x\in L_{1}: c(x,z)=\text{cyan}\}| \\ 
 & = & |\{x\in \Gamma: c(x,y)=\text{aqua},c(x,z)=\text{cyan}\}| \\
 & = & |\{x\in \Gamma: c(x,y)=\text{aqua},c(z,x)=\text{cyan}^{-1}\}| \\ 
 & = & \gamma(\text{aqua},\text{cyan}^{-1},\text{beige})
 \end{eqnarray*}
\bigskip

\noindent{\sc Solution to Ex.~\ref{ex:biparcomp}} --- Notice that, for $y\in B_i$ and $x\in C_1$, $(x,i)\in D$ if and only if there exist $y_{0},y_{1},\dotsc,y_{k}=y$ such that $(x,y_{0})$ is green and $(y_{j},y_{j+1})$ is red for $0\leq j<k$. By Exercise~\ref{ex:samesize}\ref{it:bibr1}, an edge color $r$ would know whether such a walk exists between vertices $x,y$ with $c(x,y)=r$; calling then $\mathcal{C}'$ the subset of colors $r$ for which these walks exist, for any $1\leq i\leq m$ and any $y\in B_{i}$ we have:
\begin{equation*}
 d^{-}_{Y}(i) = \sum_{r\in\mathcal{C}'}|\{x\in C_{1}: c(x,y)=r\}| = \sum_{r\in\mathcal{C}'}d^{-}_{r}(y)
\end{equation*}
where $d^{-}_{r}(y)$ is the indegree of $y$ in the graph $\mathcal{G}_{r}$. These numbers are independent from $y$ by Exercise~\ref{ex:biparhom}\ref{it:bar1}, so all vertices in $\{1,\dotsc,m\}$ have the same degree in $Y$.

Now fix $x\in C_1$ and $y\in B_i$ such that $(x,y)$ is green: as shown in the proof of Exercise~\ref{ex:samesize}\ref{it:bibr2}, we can define the subset $\mathcal{C}'\subset\mathcal{C}$ of edge colors that sit inside a connected component of $\mathcal{G}_{\text{red}}$, i.e. $y'\in B_{i}$ if and only if $c(y',y)\in\mathcal{C}'$; the number $q$ of vertices $z\in B_{i}$ such that $(x,z)$ is green is then:
\begin{eqnarray*}
 q & = & |\{z\in B_{i}: (x,z)=\text{green}\}| \\
 & = & |\{z\in\Gamma: (z,y)\in\mathcal{C}',(x,z)=\text{green}\}| \\
 & = & \sum_{r\in\mathcal{C}'}|\{z\in\Gamma: (z,y)=r,(x,z)=\text{green}\}| \\ 
 &= & \sum_{r\in\mathcal{C}'}\gamma(r,\text{green},\text{green})
\end{eqnarray*}
The second equality is a consequence of the fact that we can pass from $z\in C_{2}$ to $z\in\Gamma$ without any problem, since in $z\in\Gamma\setminus C_{2}$ neither the condition $(x,z)=\text{green}$ nor $(z,y)\in\mathcal{C}'$ are ever realized.

The constants $\gamma$ do not depend on $x$, $y$ or $i$, so the same is true for $q$; on the other hand, for any $x\in C_{1}$ we have the same number $d^{+}_{\text{green}}$ of outgoing green edges by Exercise~\ref{ex:biparhom}\ref{it:bar1}. Combining these two things together, the number $d^{+}_{Y}(x)$ of outgoing edges in $Y$ is also independent from $x$, being simply $d^{+}_{\text{green}}/q$; thus, $Y$ is semiregular.
\bigskip

\noindent{\sc Solution to Ex.~\ref{ex:gita}} --- Either one between $A$ and $V_{1}\times V_{2}\setminus A$ has size $\leq |V_{1}||V_{2}|/2$: the graph $(V_1,V_2;V_{1}\times V_{2}\setminus A)$ moreover is still a non-trivial semiregular bipartite graph and $v_{1},v_{2}$ are twins in it if and only if they are twins in $(V_1,V_2;A)$, so we can suppose that we are in the situation where $|A|\leq |V_{1}||V_{2}|/2$ by passing to the complement if necessary.

Let $d_{2}$ be the degree of the vertices of $V_{2}$, let $S$ be a twin class in $V_{1}$ and let $v_{1}$ be any element of $S$: since the graph is non-trivial and semiregular, $v_{1}$ must have an edge going to a certain $v_{2}\in V_{2}$, therefore all elements of $S$ have $v_{2}$ as neighbor; this means that $d_{2}\geq |S|$, and then $|A|\geq |S||V_{2}|$, which, combined with our previous bound, gives $|S|\leq |V_{1}|/2$.
\bigskip

\noindent{\sc Solution to Ex.~\ref{ex:reduc2}} --- Suppose that we have $V_2 = C_1\cup C_2$, $C_1\cap C_2=\emptyset$ and that there are classes $S_{1},S_{2}\subset V_{1}$ of twins with respect to the graphs $(V_1,C_1;A\cap (V_1\times C_1))$ and $(V_1,C_2;A\cap (V_1\times C_2))$ (respectively). Suppose also that $|S_{1}|,|S_{2}|\geq|V_{1}|/2+1$: this implies that $|S_{1}\cap S_{2}|\geq 2$, so there exist two vertices $v,w\in V_{1}$ that are twins with respect to both graphs; but then, for every vertex $z\in V_{2}$, we have that $(v,z)\in A \Leftrightarrow (w,z)\in A$ whether we have $z\in C_{1}$ or $z\in C_{2}$: in other words, $v$ and $w$ are twins in the whole $(V_1,V_2;A)$.
\bigskip

\noindent{\sc Solution to Ex.~\ref{ex:twindet}} -- Two vertices $x,y\in C_{1}$ are not twins if and only if we have $c(x,z)=\text{brown}$ and $c(y,z)\neq\text{brown}$ (or vice versa) for at least one $z\in C_{2}$; in other words, $x,y$ are not twins if and only if at least one of the constants $\gamma(\text{brown}^{-1},r,c(x,y))$ or $\gamma(r^{-1},\text{brown},c(x,y))$ (for all $r\neq\text{brown}$) is positive. These constants are known to the color $c(x,y)$, so it knows whether $x$ and $y$ are twins or not.
\bigskip

\begin{center}
***
\end{center}

\bigskip
\begin{exer}[\S~\ref{subs:Schreier-Sims}]\label{hidden0}
Show that being able to describe the setwise stabilizer $G_{\{x_{1},...,x_{k}\}}$ of any permutation group $G$ acting on a finite set $\Omega$ for any choice of $x_{1},...,x_{k}\in\Omega$ implies being able to solve the string isomorphism problem.
\end{exer}
\bigskip

\begin{exer}[\S~\ref{subs:orbl}]\label{hidden1}
Let $G$ be a permutation group acting on a finite set $\Omega$. Then the orbits \mbox{$\{x^g: g\in G\}$} of the action of $G$ on $\Omega$ can be determined in polynomial time in $|\Omega|$.
\end{exer}
\bigskip

\begin{exer}[\S~\ref{subs:luks}]\label{hidden1b}
In this exercise we prove that the graph isomorphism problem for graphs of bounded degree reduces in polynomial time to the string isomorphism problem for a group $G$ whose composition factors are all bounded. We follow \cite{Lu}; we suppose all graphs to be undirected simple connected graphs.
\begin{enumerate}
\item\label{hidden1b-Aute} Let $\Gamma_{1},\Gamma_{2}$ be two graphs of size $n$ whose vertices have degree $\leq k$ ($k\geq 3$): show that we can reduce the problem of finding $\text{\normalfont Iso}(\Gamma_{1},\Gamma_{2})$ to $\leq kn$ problems of finding $\text{\normalfont Aut}_{e}(\Gamma)$ for some graph $\Gamma$ of size $2n+2$ and maximum degree $k$, where $\text{\normalfont Aut}_{e}(\Gamma)$ is the group of automorphisms of $\Gamma$ fixing a certain edge $e$.
\item\label{hidden1b-ker} Call $\Gamma_{r}$ the subgraph of $\Gamma$ given by all the vertices and edges appearing in walks of length $\leq r$ through $e$; call $\pi_{r}:\text{\normalfont Aut}_{e}(\Gamma_{r+1})\rightarrow\text{\normalfont Aut}_{e}(\Gamma_{r})$ the map that sends $\sigma$ to its restriction to $\Gamma_{r}$. Find a set of generators for $\ker\pi_{r}$.
\item\label{hidden1b-im} Let $K$ be the set of all subsets of size $\leq k$ of vertices of $\Gamma_{r}$. Show that finding a set of generators for the image of $\pi_{r}$ reduces to describing $\text{\normalfont Aut}_{G}(\mathbf{x})$ for some string $\mathbf{x}$ of length $|K|$ and some appropriate $G<\Sym(K)$.
\item\label{hidden1b-subsym} Let $\mathfrak{G}_{s}$ be the class of all groups $G$ such that all composition factors $G_{i}/G_{i+1}$ of a composition series $\{e\}=G_{m}\lhd G_{m-1}\lhd\ldots\lhd G_{1}\lhd G_{0}=G$ are subgroups of $\Sym_{s}$. Show that every subgroup of $\Sym_{s}$ is in $\mathfrak{G}_{s}$.
\item\label{hidden1b-subg} If $N\lhd G$, by the Jordan-H\"older Theorem $G\in\mathfrak{G}_{s}$ if and only if $N,G/N\in\mathfrak{G}_{s}$. Using this and (\ref{hidden1b-subsym}), show that if $G$ is in $\mathfrak{G}_{s}$ then every subgroup of $G$ is in $\mathfrak{G}_{s}$.
\item\label{hidden1b-AutG} Using Jordan-H\"older and (\ref{hidden1b-ker})-(\ref{hidden1b-subg}), show that $\text{\normalfont Aut}_{e}(\Gamma_{r})\in\mathfrak{G}_{k}$ for every $r$.
\item\label{hidden1b-concl} Conclude the proof.
\end{enumerate}
\end{exer}
\bigskip

\begin{exer}[\S~\ref{subs:secf1}]\label{hidden1d}
Let $\mathfrak{X}=(\Gamma,c)$ be a partition structure and let $\Gamma'\subset\Gamma$.
\begin{enumerate}
 \item\label{hidden1d-conf} If $\mathfrak{X}$ is a configuration, $\mathfrak{X}[\Gamma']$ is a configuration.
 \item\label{hidden1d-coh} If $\mathfrak{X}$ is coherent, $\mathfrak{X}[\Gamma']$ is not necessarily coherent.
\end{enumerate}
\end{exer}
\bigskip

\begin{exer}[Def.~\ref{def:uniprim}]\label{hidden2}
Let $(\Gamma,c)$ be a homogeneous classical coherent configuration. For any edge color $r$, let $\mathscr{G}_r$ be the graph $\{(x,y): x,y\in \Gamma, c(x,y)=r\}$.
\begin{enumerate}
\item Show that the outdegree $d^+_r(x)=|\{y\in \Gamma: (x,y)\in \mathscr{G}_r\}|$ is independent from $x$.
\item Show that any weakly connected component of $\mathscr{G}_r$ is connected.
\end{enumerate}
\end{exer}
\bigskip

\begin{exer}[Lemma~\ref{lem:mustaf}]\label{hidden3}
Let $G<\Sym(\Omega)$ be primitive. Let $\phi:G\to \Alt_k$ be an epimorphism with $k>\max(8,2+\log_2 |\Omega|)$. Then $\phi$ is an isomorphism.
\end{exer}
\bigskip

\begin{exer}[Prop.~\ref{prop:bicoup}]\label{hidden4}
Let $(\Gamma,c)$ be a uniprimitive classical coherent configuration. Then, for any edge color $r$, the graph $\overline{\mathscr{G}_r}=\{(x,y): x,y\in \Gamma, c(x,y)\neq r\}$ is of diameter $2$.
\end{exer}
\bigskip

\begin{exer}\label{WLcoh}
Let $\mathfrak{X}$ be a configuration; suppose that we apply the algorithm of Weisfeiler-Leman to $\mathfrak{X}$.
\begin{enumerate}
\item\label{WLcoh-coh} Prove that we obtain a coherent configuration.
\item\label{WLcoh-coarse} Prove that we obtain the coarsest coherent configuration that is a refinement of $\mathfrak{X}$.
\end{enumerate}
\end{exer}
\bigskip

\begin{exer}\label{canonord}
Let $\mathscr{C}$ be a finite set of colors; for any $k$-ary relational structure ${\mathfrak X}=(\Gamma,(R_{i})_{i\in\mathscr{C}})$, the functors $F_{1},F_{2},F_{3}$ (\S\S~\ref{subs:secf1}-\ref{subs:wl}) refine the coloring of ${\mathfrak X}$ giving a new $\mathscr{C}'$ (in general much bigger than $\mathscr{C}$) that eventually makes it a coherent configuration. Suppose that $\mathscr{C}$ is an ordered set (we consider all orderings to be well-orderings): show that the definitions of the three functors allow us to give canonically an ordering of $\mathscr{C}'$ from the ordering of $\mathscr{C}$.
\end{exer}

\noindent{\sc Comment} --- This is an important detail. A canonical ordering of the colors allows us to have no multiplication cost associated with the refinement: if the coloring is canonical but not canonically ordered, the two graphs $\Gamma_{1},\Gamma_{2}$ are divided into $m$ color classes each (say) but all possible correspondences between classes have to be considered, so that $m!$ combinations have to be examined by the algorithm. If the coloring is canonically ordered, the ``first'' color of $\Gamma_{1}$ necessarily corresponds to the ``first'' color of $\Gamma_{2}$, etc...; the algorithm has to run only for that one combination.
\bigskip

\begin{exer}\label{WLaut}
Let $\mathfrak{X}=(\Gamma,c)$ be a $k$-ary partition structure and let $\mathfrak{X}'=(\Gamma,c')$ be the coherent configuration obtained after applying functors $F_{2},F_{3}$ to $\mathfrak{X}$. Prove that for any $\sigma\in\mbox{\normalfont Aut}(\mathfrak{X})$ and any $\vec{x}\in\Gamma^{k}$ we have $c'(\sigma(\vec{x}))=c'(\vec{x})$.
\end{exer}

\noindent{\sc Comment} --- For $\sigma\in\mbox{Aut}(\mathfrak{X})$ we have $c(\sigma(\vec{x}))=c(\vec{x})$ by definition, but $c'$ is more refined than $c$. This is what makes the functors $F_{2},F_{3}$ interesting, especially Weisfeiler-Leman: thanks to it, the coloring of $\Gamma^{k}$ becomes more refined without losing any isomorphism that $\mathfrak{X}$ has in origin.
\bigskip

\begin{exer}\label{graphconfig}
Let $\Gamma$ be a non-trivial undirected simple graph. Define the coloring $c:\Gamma^{2}\rightarrow\mathscr{C}$ with $\mathscr{C}=\{\mbox{``vertex''},\mbox{``edge''},\mbox{``not edge''}\}$ in the obvious way; prove that $(\Gamma,c)$ is a configuration.

Prove also that, with such coloring:
\begin{enumerate}
 \item\label{graphconfig-Pn} $P_{n}$ (the path graph on $n$ vertices) is not coherent for $n\geq 3$;
 \item\label{graphconfig-C4} $C_{4}$ (the cycle graph on $4$ vertices) is coherent;
 \item\label{graphconfig-C6} $C_{6}$ is not coherent.
\end{enumerate}
\end{exer}
\bigskip

\begin{exer}\label{stronglyreg}
A graph $\Gamma$ is {\em strongly regular} if there are three constants $k,\lambda,\mu$ such that:
\begin{enumerate}
 \item every vertex has $k$ neighbors;
 \item every two adjacent vertices have $\lambda$ common neighbors;
 \item every two non-adjacent vertices have $\mu$ common neighbors.
\end{enumerate}
Show that a graph $\Gamma$ with the coloring defined in Exercise~\ref{graphconfig} is a coherent configuration if and only if $\Gamma$ is strongly regular.
\end{exer}
\bigskip

\begin{exer}\label{symmBIBD}
A BIBD (balanced incomplete block design, see \S~\ref{subs:hypdes}) of parameters $(v,u,\lambda)$ and number of blocks\footnote{In the paper we used ``edge'' instead of ``block'', since we were defining a BIBD starting from a hypergraph; in this exercise we use the term ``block'' (which is the usual term in literature, as proven by the acronym BIBD itself) in order to avoid confusion with the use of the term ``edge'' in the context of configurations.} $|\mathscr{A}|=b$ is {\em symmetric} when Fisher's inequality is in fact an equality, i.e. $v=b$; in a symmetric BIBD, the intersection of any two distinct blocks is of size $\lambda$. Call $\Gamma=V\cup\mathscr{A}$, and define a coloring $c$ of $\Gamma^{2}$ as follows:
\begin{enumerate}
 \item $c(x,x)=\mbox{``vertex''}$, for any $x\in V$ or $x\in\mathscr{A}$;
 \item $c(x,y)=\mbox{``white''}$, for any $(x,y)\in V\times V$ or $(x,y)\in\mathscr{A}\times\mathscr{A}$;
 \item $c(x,y)=c(y,x)=\mbox{``belongs''}$, for any $(x,y)\in V\times\mathscr{A}$ such that $x\in y$;
 \item $c(x,y)=c(y,x)=\mbox{``doesn't belong''}$, for any $(x,y)\in V\times\mathscr{A}$ such that $x\not\in y$.
\end{enumerate}
Show that $(\Gamma,c)$ is a coherent configuration.

Show that, if the BIBD is not symmetric, $(\Gamma,c)$ is a configuration but the Weisfeiler-Leman algorithm would give to vertices in $V$ and in $\mathscr{A}$ two different colors; therefore, $(\Gamma,c)$ is not coherent.

For a BIBD that is not symmetric, not all intersections of two distinct blocks have the same size (\cite[Thm.~3]{RChW}). Using this result, show the following: even if we refine the colors in the previous construction by distinguishing between ``vertex in $V$'' and ``vertex in $\mathscr{A}$'', ``white in $V$'' and ``white in $\mathscr{A}$'', ``belongs'' and ``contains'', ``doesn't belong'' and ``doesn't contain'', the configuration $(\Gamma,c)$ for a BIBD that is not symmetric is still not coherent.
\end{exer}

\noindent{\sc Comment} --- Notice that by the first part of the exercise we are not able to distinguish between vertices and blocks in a symmetric BIBD: Weisfeiler-Leman would not be able to give two different colors for $V$ and $\mathscr{A}$, for any $(x,y)$ of color ``belongs'' it would not be able to say whether $x\in y$ or $y\in x$. This is a powerful duality, reflecting the results $v=b$, $r=u$, $\lambda=\lambda$ that we know to be true for a symmetric BIBD and false for any other BIBD.
\bigskip

\begin{exer}\label{tdesign}
In this exercise we prove Proposition~\ref{prop:RW}. We follow \cite{RChW}.
\begin{enumerate}
 \item\label{tdesign-b} Let $(V,\mathscr{A})$ be a $t-(v,u,\lambda)$ design (see \S~\ref{subs:hypdes}) with $b=|\mathscr{A}|$. Prove that we have $b=\lambda\binom{v}{t}\binom{u}{t}^{-1}$.
 \item\label{tdesign-t'} Prove that any $t$-design is a $t'$-design for any $t'\leq t$, with new $\lambda_{t'}=\lambda\binom{v-t'}{t-t'}\binom{u-t'}{t-t'}^{-1}$. {\em Hint: show that a $t$-design is a $(t-1)$-design by fixing a $(t-1)$-set $X$ and counting pairs $(A,x)$ where $A$ is an edge containing $X\cup\{x\}$.}
 \item\label{tdesign-ij} Let $i,j\geq 0$ be two integers with $i+j\leq t$ and let $I,J$ be two disjoint subsets of $V$ with size $i,j$ respectively: prove that the number of edges that contain $I$ and are disjoint from $J$ is $\lambda_{i,j}=\lambda\binom{v-i-j}{u-i}\binom{v-t}{u-t}^{-1}$. {\em Hint: use (\ref{tdesign-t'}) to show that $\lambda_{i,j}\lambda^{-1}$ depends only on the parameters $v,u,t,i,j$ and not on the specific design, then consider the design having as edges all possible $u$-subsets of $V$.}
 \item\label{tdesign-FasE} From now on we suppose $t$ even, $t=2s$, and $v\geq s+u$. For any edge $A$, define $\hat{A}$ as the formal sum of all $s$-subsets of $V$ contained in $A$. Fix an $s$-subset $S_{0}\subseteq V$: define $E_{i}$ as the formal sum of all $s$-subsets $S$ with $|S\cap S_{0}|=s-i$; define $F_{j}$ as the sum of all $\hat{A}$ with $|A\cap S_{0}|=s-j$. Using (\ref{tdesign-ij}), express every $F_{j}$ as sum of the $E_{i}$.
 \item\label{tdesign-EasF} Show that it is possible to solve for the $E_{i}$ as linear combinations of the $F_{j}$.
 \item\label{tdesign-concl} Consider the free vector space generated by all $s$-subsets of $V$ (over ${\mathbb Q}$, say): it has dimension $\binom{v}{s}$. Using (\ref{tdesign-EasF}), conclude the proof under our hypotheses ($t=2s$, $v\geq s+u$).
 \item\label{tdesign-other} Cover the remaining cases, namely $t$ odd and $v<s+u$.
\end{enumerate}
\end{exer}
\bigskip

\begin{exer}\label{johnson}
Let $(\Gamma,c)$ be the following homogeneous classical coherent configuration. $\Gamma=\{v_{1},v_{2},...,v_{15}\}$, where $15=\binom{6}{2}$, and the coloring $c:\Gamma^{2}\rightarrow\mathscr{C}=\{x,y,z\}$ is defined by the following relations:
\bigskip

\begin{tabular}{|c|c|c|c|c|c|c|c|c|c|c|c|c|c|c|c|}
\hline
$c(\cdot,\cdot)$ & $v_{1}$ & $v_{2}$ & $v_{3}$ & $v_{4}$ & $v_{5}$ & $v_{6}$ & $v_{7}$ & $v_{8}$ & $v_{9}$ & $v_{10}$ & $v_{11}$ & $v_{12}$ & $v_{13}$ & $v_{14}$ & $v_{15}$ \\
\hline
$v_{1}$ & x & z & y & y & z & y & y & z & z & z & y & y & z & z & z \\
\hline
$v_{2}$ & z & x & y & z & z & y & y & y & z & y & z & z & z & z & y \\
\hline
$v_{3}$ & y & y & x & y & y & z & z & z & y & z & z & z & z & z & y \\
\hline
$v_{4}$ & y & z & y & x & z & z & z & y & z & y & z & z & y & y & z \\
\hline
$v_{5}$ & z & z & y & z & x & y & y & z & y & z & z & z & y & y & z \\
\hline
$v_{6}$ & y & y & z & z & y & x & z & z & z & y & z & y & y & z & z \\
\hline
$v_{7}$ & y & y & z & z & y & z & x & y & z & z & y & z & z & y & z \\
\hline
$v_{8}$ & z & y & z & y & z & z & y & x & y & z & z & y & y & z & z \\
\hline
$v_{9}$ & z & z & y & z & y & z & z & y & x & y & y & y & z & z & z \\
\hline
$v_{10}$ & z & y & z & y & z & y & z & z & y & x & y & z & z & y & z \\
\hline
$v_{11}$ & y & z & z & z & z & z & y & z & y & y & x & z & y & z & y \\
\hline
$v_{12}$ & y & z & z & z & z & y & z & y & y & z & z & x & z & y & y \\
\hline
$v_{13}$ & z & z & z & y & y & y & z & y & z & z & y & z & x & z & y \\
\hline
$v_{14}$ & z & z & z & y & y & z & y & z & z & y & z & y & z & x & y \\
\hline
$v_{15}$ & z & y & y & z & z & z & z & z & z & z & y & y & y & y & x \\
\hline
\end{tabular}
\bigskip

Following the procedure given in \S~\ref{subs:idgroup}, find two bijections $\iota:\Gamma\rightarrow\mathscr{S}_{2}(\Lambda)$ (with $|\Lambda|=6$) and $\iota':\mathscr{C}\rightarrow\{0,1,2\}$ that give an isomorphism between $(\Gamma,c)$ and the Johnson scheme defined by $\mathscr{S}_{2}(\Lambda)$; in other words, find $\iota,\iota'$ such that $|\iota(v_{i})\cap\iota(v_{j})|=\iota'(c(v_{i},v_{j}))$ for all $1\leq i,j\leq 15$.

Warning! Since the condition $m>(k+1)^{2}-2$ is not satisfied in this example with $(m,k)=(6,2)$, here we swap the definitions and impose that $\Upsilon$ is the biggest orbital and $\Delta$ is the smallest.
\end{exer}
\bigskip

\begin{exer}\label{notprimcoh}
Let $(\Gamma,c)$ be a homogeneous classical coherent configuration that is not primitive, and let $r$ be any edge color such that the graph $\mathscr{G}_{r}$ is disconnected. Call $B_{1},B_{2},...,B_{m}$ the connected components of $\mathscr{G}_{r}$; then, for any other edge color $r'$:
\begin{enumerate}
 \item either $\mathscr{G}_{r'}$ is also disconnected and the partition into connected components is equal or finer than the one determined by $B_{1},B_{2},...,B_{m}$,
 \item or $\mathscr{G}_{r'}$ is a subgraph of the complete $m$-partite graph defined by the sets $B_{1},B_{2},...,B_{m}$.
\end{enumerate}
\end{exer}
\bigskip

\begin{exer}\label{squaredom}
Let $(\Gamma,c)$ be a uniprimitive classical coherent configuration with unique vertex color $s$; for any edge color $r$, define $\gamma_{r}$ to be the number of pairs of color $r$ with the same starting point (so that $\gamma_{r}=\gamma(r^{-1},r,s)$, see Def.~\ref{def:coh}). Then, for any edge color $r$, there exists an edge color $r'$ such that $\gamma_{r'}\geq\sqrt{\gamma_{r}}$. {\em Hint: Exercise~\ref{hidden4} is a good start.}
\end{exer}

\noindent{\sc Comment} --- In particular, when we call Bipartite Split-or-Johnson for the first time inside Split-or-Johnson, we can suppose not only the bound $|V_{2}|<\beta|V_{1}|$ but also $|V_{2}|\geq \sqrt{|V_{1}|}$.
\bigskip

\begin{exer}\label{colpart}
Let $(\Gamma,c)$ be a uniprimitive classical coherent configuration with unique vertex color blue; suppose that there exists an edge color, say red, such that from any vertex $v$ there are $\gamma=\gamma(\mbox{red}^{-1},\mbox{red},\mbox{blue})$ red pairs $(v,w)$. Show that, giving to an arbitrarily chosen $v$ a new color (say green) and applying the functors $F_{2},F_{3}$, we obtain a $(1-(\gamma+1)^{-1})$-coloring of $\Gamma$.
\end{exer}

\noindent{\sc Comment} --- This is a process of {\em individualization} of one vertex, as described in \cite{SW}. It leads to a loss of canonicity that determines a multiplicative cost of $|\Gamma|$ in exchange for a reduction to problems of size $\leq (1-(\gamma+1)^{-1})|\Gamma|$; if $\gamma$ is bounded, this process is compatible with a bound of $n^{O(\log n)}$ time for the algorithm.
\bigskip

\begin{center}
***
\end{center}

\bigskip
\noindent{\sc Solution to Ex.~\ref{hidden0}} --- Suppose that the time to determine the setwise stabilizer of any subgroup of $\Sym(\Omega)$ for any subset of $\Omega$ is bounded by $T$. We prove that for every string $\mathbf{x}:\Omega\rightarrow\Sigma$ we can determine $\Aut_{G}(\mathbf{x})$ in $O(T|\Omega|)$.

We write $\Omega$ as the disjoint union of the preimages of each element $\sigma\in\Sigma$, say $\Omega=S_{1}\sqcup...\sqcup S_{m}$ with $S_ {i}=\mathbf{x}^{-1}(\sigma_{i})$ for $\Sigma=\{\sigma_{1},...,\sigma_{m}\}$; we define $G_{0}=G$ and $G_{i}=(G_{i-1})_{S_{i}}$ for every $1\leq i\leq m$. $G_{m}$ is the group of elements of $G$ that stabilize each preimage $\mathbf{x}^{-1}(\sigma_{i})$, so we have exactly $G_{m}=\Aut_{G}(\mathbf{x})$; the execution time is $O(Tm)\leq O(T|\Omega|)$.
\bigskip

\noindent{\sc Solution to Ex.~\ref{hidden1}} --- Let $A$ be a set of generators of $G$. First, we determine all the sets $A_{x}=\{x^{a}:a\in A\}$ for every $x\in\Omega$: this takes time $O(|A||\Omega|)$. After that, we follow this procedure: we start with any fixed $x_{0}\in\Omega$ and set $S_{x_{0}}=\{x_{0}\}\cup A_{x_{0}}$; at this stage, the only element of $S_{x_{0}}$ that we have examined is $x_{0}$. Then at every step we take a still unexamined $x\in S_{x_{0}}$ and we update $S_{x_{0}}$ by adding the elements of $A_{x}$ to it: the newly added elements are marked as unexamined, while $x$ now is examined; the procedure stops when $S_{x_{0}}$ becomes the orbit $\{x_{0}^g: g\in G\}$. If there is an element $x_{1}$ that has not yet been considered, we define $S_{x_{1}}=\{x_{1}\}\cup A_{x_{1}}$ and go through the whole procedure again, until we have considered all the elements of $\Omega$: the final sets $S_{x_{0}},S_{x_{1}},...,S_{x_{m}}$ are the orbits of the action of $G$ on $\Omega$; this part takes time $O(|\Omega|)$.

The execution time of this algorithm is $O(|A||\Omega|+|\Omega|)$. Thanks to the Schreier-Sims algorithm (Algorithm~1) we can suppose $|A|=O(|\Omega|^{2})$, so the runtime is $O(|\Omega|^{3})$.
\bigskip

\noindent{\sc Solution to Ex.~\ref{hidden1b}} --- (\ref{hidden1b-Aute}) Since we have $\leq kn$ edges for each graph, $\text{Iso}(\Gamma_{1},\Gamma_{2})$ is the union of $\leq kn$ sets $\text{Iso}_{e_{1}\mapsto e_{2}}(\Gamma_{1},\Gamma_{2})$ where $e_{1}$ is a fixed edge of $\Gamma_{1}$, $e_{2}$ runs through all edges of $\Gamma_{2}$ and the isomorphisms inside each of these sets send $e_{1}$ to $e_{2}$; in order to describe $\text{Iso}(\Gamma_{1},\Gamma_{2})$, we need to give a set $A$ of generators of $\text{Aut}(\Gamma_{1})$ (or $\text{Aut}(\Gamma_{2})$) and at least one element $\sigma$ of at least one of the $\text{Iso}_{e_{1}\mapsto e_{2}}(\Gamma_{1},\Gamma_{2})$: the coset $\sigma\langle A\rangle$ (or $\langle A\rangle\sigma$) will be exactly $\text{Iso}(\Gamma_{1},\Gamma_{2})$.

Call $e_{1}=\{x_{1},y_{1}\}$; for each choice of $e_{2}=\{x_{2},y_{2}\}$ we build the following $\Gamma$:
\begin{eqnarray*}
V(\Gamma) & = & V(\Gamma_{1})\cup V(\Gamma_{2})\cup\{p,q\} \\
E(\Gamma) & = & (E(\Gamma_{1})\setminus\{e_{1}\})\cup(E(\Gamma_{2})\setminus\{e_{2}\})\cup\{\{x_{1},p\},\{p,y_{1}\},\{x_{2},q\},\{q,y_{2}\},\{p,q\}\}
\end{eqnarray*}
Figuratively speaking, we have connected the ``middle points'' $p,q$ of $e_{1},e_{2}$ by a new edge $\{p,q\}$ (which we call $e$). $\Gamma$ is of size $2n+2$, its maximum degree is still $k$ ($p,q$ have degree $3$ by construction) and the automorphisms of $\Gamma$ that fix $e$ are exactly those that are either pairs of automorphisms of $\Gamma_{1},\Gamma_{2}$ (when $p,q$ are fixed) or isomorphisms sending $\Gamma_{1}$ to $\Gamma_{2}$ and $e_{1}$ to $e_{2}$ (when $p,q$ are transposed). For each set $A_{e_{2}}$ of generators of $\text{Aut}_{e}(\Gamma)$ found for every choice of $e_{2}$, we write $A_{e_{2}}=A_{e_{2}}^{(1)}\sqcup A_{e_{2}}^{(2)}$ for the subsets of elements with $p,q$ fixed and $p,q$ transposed, respectively: if $A_{e_{2}}^{(2)}=\emptyset$ for every choice of $e_{2}$, then we have that there are no isomorphisms sending $\Gamma_{1}$ to $\Gamma_{2}$, because otherwise for at least one $e_{2}$ at least one generator of $\text{Aut}_{e}(\Gamma)$ must be transposing $p,q$; if we have at least one element $\sigma$ inside $A_{e_{2}}^{(2)}$, then we have found the element of $\text{Iso}_{e_{1}\mapsto e_{2}}(\Gamma_{1},\Gamma_{2})$ that we were looking for, and to determine generators for the subgroup of automorphisms of $\Gamma$ that fix $p,q$ (which once restricted correspond to generators of $\text{Aut}(\Gamma_{1})$ or $\text{Aut}(\Gamma_{2})$) we simply need to use  the procedure given in Exercise~\ref{ex:fhl}\ref{it:richt}: this subgroup in fact has index $2$ and its polynomial-time test is simply ``is $p$ fixed?'', so this procedure takes polynomial time and we are done.

(\ref{hidden1b-ker}) Partition the vertices in $\Gamma_{r+1}\setminus\Gamma_{r}$ according to whether they have the same neighborhood in $\Gamma_{r+1}$, i.e. they are connected by an edge to the same vertices of $\Gamma_{r}$: an automorphism $\sigma\in\text{Aut}_{e}(\Gamma_{r+1})$ is in the kernel of $\pi_{r}$ if and only if it preserves all vertices of $\Gamma_{r}$, which means that in $\Gamma_{r+1}\setminus\Gamma_{r}$ it preserves the partition that we defined; on the other hand, every permutation of $\Gamma_{r+1}\setminus\Gamma_{r}$ preserving this partition gives obviously an element of $\pi_{r}$, just by extending it to the whole $\Gamma_{r+1}$ imposing that its action on $\Gamma_{r}$ is the identity. So $\ker\pi_{r}$ is isomorphic to the direct product of $\Sym(n_{r,i})$, where $n_{r,i}$ is the size of each part of our partition: each $\Sym(n_{r,i})$ is generated by two elements, for example $(1\ 2)$ and $(1\ 2\ \cdots\ n_{r,i})$, therefore we have found at most $|\Gamma_{r+1}\setminus\Gamma_{r}|$ generators of $\ker\pi_{r}$ of the form $(v_{1}\ v_{2})$ or $(v_{1}\ v_{2}\ \cdots\ v_{n_{r,i}})$.

(\ref{hidden1b-im}) In order to understand what the image of $\pi_{r}$ is, we have to consider the new structures that appear when we pass from $\Gamma_{r}$ to $\Gamma_{r+1}$: in considering the larger subgraph, we now have new vertices in $\Gamma_{r+1}\setminus\Gamma_{r}$ and new edges between vertices of $\Gamma_{r}$; any $\sigma\in\text{Aut}_{e}(\Gamma_{r})$ that extends to $\Gamma_{r+1}$ (i.e. it is in the image of $\pi_{r}$) must respect the restrictions given by these new vertices and edges.

Call $K_{i}$ the set of subsets of vertices of $\Gamma_{r}$ that are neighbors of the same $i$ vertices of $\Gamma_{r+1}\setminus\Gamma_{r}$: if $\sigma\in\text{im}\,\pi_{r}$, each subset that belongs to $K_{i}$ must go to another subset in $K_{i}$, since $i$ vertices of $\Gamma_{r+1}$ with the same neighbors in $\Gamma_{r}$ will necessarily go to other $i$ vertices with the same property; call $K'$ the set of $2$-subsets $\{x,y\}\subset\Gamma_{r}$ such that $\{x,y\}$ is a new edge of $\Gamma_{r+1}$: again, if $\sigma\in\text{im}\,\pi_{r}$ then new edges are sent to new edges by $\pi_{r}^{-1}(\sigma)$, therefore $K'$ is preserved by $\sigma$. On the other hand, we can show that preserving $K'$ and all $K_{i}$ is a sufficient condition for $\sigma$ to be a member of $\text{im}\,\pi_{r}$: given $\sigma\in\text{Aut}_{e}(\Gamma_{r})$, if $\sigma$ preserves each $K_{i}$ we define an extension $\sigma_{r+1}$ by sending each set $\{w_{1},...,w_{i}\}$ of neighbors of $\{v_{1},...,v_{j}\}$ to the set $\{z_{1},...,z_{i}\}$ of neighbors of $\{\sigma(v_{1}),...,\sigma(v_{j})\}$ (any possible way to send each $w_{i'}$ to some $z_{i''}$ is acceptable); $\sigma_{r+1}$ is a permutation of the vertices of $\Gamma_{r+1}$ that respect the old edges (because it agrees with $\sigma$ on $\Gamma_{r}$), the edges in $\Gamma_{r+1}\times\Gamma_{r}$ (by construction) and the new edges between vertices of $\Gamma_{r}$ (using the condition that $\sigma$ respects $K'$). Thus, $\sigma$ is in the image of $\pi_{r}$ if and only if it stabilizes $K'$ and every $K_{i}$.

Let $K$ be the set of all subsets of size $\leq k$ of vertices of $\Gamma_{r}$: in particular $K\supset K',K_{i}$ for every $i$, because any vertex in $\Gamma_{r+1}$ has at most $k$ neighbors in $\Gamma_{r}$; also, for the same reason every $K_{i}$ with $i\geq k$ is empty (every vertex in $\Gamma_{r}$ has at most $k-1$ neighbors in $\Gamma_{r+1}$). We build $\mathbf{x}$ that assigns to each element of $K$ a symbol in an alphabet $\Sigma$ of cardinality $\leq 2k$, according to whether this element belongs to $K'$ and/or to one of the $K_{i}$ (the $K_{i}$ are all disjoint from each other but they can intersect $K'$); there is a natural injection $\iota:\text{Aut}_{e}(\Gamma_{r})\hookrightarrow\Sym(K)$: by the discussion above, $\sigma\in\text{im}\,\pi_{r}$ if and only if $\iota(\sigma)\in\text{Aut}_{G}(\mathbf{x})$, where $G=\iota(\text{Aut}_{e}(\Gamma_{r}))$. After finding generators for $\text{Aut}_{G}(\mathbf{x})$, we can take their preimage and we will get generators of $\text{im}\,\pi_{r}$.

(\ref{hidden1b-subsym}) Let $G$ be a subgroup of $\Sym_{s}$ and let $N\lhd G$ be such that $G/N$ is a simple group: we want to prove that $G/N$ is a subgroup of $\Sym_{s}$ too. We have a chain of subgroups $N=G_{s-1}N<G_{s-2}N<...<G_{1}N<G_{0}N=G$; taking the smallest $j$ such that $G_{j+1}N\subsetneq G_{j}N$ we have $G=G_{j}N$: $G$ then acts transitively on $C=G_{j}N/G_{j+1}N=G/G_{j+1}N$, while $N$ acts trivially on $C$ because it is normal in $G$. Therefore we have an action of $G/N$ on $C$: the action is non-trivial by transitivity and faithful by the fact that $G/N$ is simple, which means that we can identify $G/N$ with a subgroup of $\Sym(C)$; finally, observing that $|C|\leq |G_{j}/G_{j+1}|\leq s$, we conclude that $G/N$ is a subgroup of $\Sym_{s}$ and that $G\in\mathfrak{G}_{s}$.

(\ref{hidden1b-subg}) Let $H<G\in\mathfrak{G}_{s}$: from a composition series $\{e\}=G_{m}\lhd G_{m-1}\lhd...\lhd G_{1}\lhd G_{0}=G$ we obtain a chain $\{e\}=G_{m}\cap H\lhd G_{m-1}\cap H\lhd...\lhd G_{1}\cap H\lhd G_{0}\cap H=H$. For every $i$, $G_{i}\cap H/G_{i+1}\cap H<G_{i}/G_{i+1}<\Sym_{s}$ since $G\in\mathfrak{G}_{s}$, thus giving $G_{i}\cap H/G_{i+1}\cap H\in\mathfrak{G}_{s}$ by (\ref{hidden1b-subsym}); then, using Jordan-H\"older and climbing up the chain we prove that each $G_{i}\cap H$ is in $\mathfrak{G}_{s}$, and in the end we get the same for $H=G_{0}\cap H$.

(\ref{hidden1b-AutG}) By (\ref{hidden1b-ker}) each of the $\ker\pi_{r}$ is a product of symmetric groups of rank $\leq k$, since the $n_{r,i}$ are sizes of sets of vertices all in the neighborhood of a same point, whose degree is $\leq k$; so $\ker\pi_{r}\in\mathfrak{G}_{k}$. Then we reason by induction on $r$: for $r=0$ we have $\text{Aut}_{e}(\Gamma_{1})=\ker\pi_{0}$, thus being in $\mathfrak{G}_{k}$; if $\text{Aut}_{e}(\Gamma_{r})\in\mathfrak{G}_{k}$ then its subgroup $\text{im}\,\pi_{r}$ is also in $\mathfrak{G}_{k}$ by (\ref{hidden1b-subg}), therefore using (\ref{hidden1b-ker}) and Jordan-H\"older on the natural isomorphism $\text{Aut}_{e}(\Gamma_{r+1})/\ker\pi_{r}\simeq\text{im}\,\pi_{r}$ we obtain $\text{Aut}_{e}(\Gamma_{r+1})\in\mathfrak{G}_{k}$. Every $\text{Aut}_{e}(\Gamma_{r})$ is therefore a member of $\mathfrak{G}_{k}$.

(\ref{hidden1b-concl}) Let us first examine steps (\ref{hidden1b-Aute})-(\ref{hidden1b-ker})-(\ref{hidden1b-im}): starting with $\Gamma_{1},\Gamma_{2}$ we have $\leq kn$ problems with $\text{Aut}_{e}(\Gamma)$, after whose solution we have polynomial-time tests on the generators to see whether they fix each vertex of $e$ and (possibly) a polynomial-time procedure as in Exercise~\ref{ex:fhl}\ref{it:richt} to determine generators of $\text{Aut}(\Gamma_{1})$ or $\text{Aut}(\Gamma_{2})$.

Now, for each problem involving $\text{Aut}_{e}(\Gamma)$, the construction of generators of the $\ker\pi_{r}$ involves only the determination of the partition described in (\ref{hidden1b-ker}): checking neighbors of all vertices of $\Gamma$ (of size $O(n)$ and maximum degree $k$) takes time $O(kn)$, and comparisons to determine whether two vertices belong to the same part of the partition takes time $O(kn^{2})$.

To construct the whole $\text{Aut}_{e}(\Gamma)$, we need to inductively build each $\text{Aut}_{e}(\Gamma_{r})$ and $\text{im}\,\pi_{r}$. $\text{Aut}_{e}(\Gamma_{1})=\ker\pi_{0}$, so we already have it. From $\text{Aut}_{e}(\Gamma_{r})$, we need first to construct the string $\mathbf{x}$: determining $K',K_{i}$ is again a matter of checking neighbors, which takes time $O(kn^{2})$, and the actual coloring of the string takes time $O(|K|)\leq n^{O(k)}$; after we have our string, we determine $\text{im}\,\pi_{r}=\text{Aut}_{G}(\mathbf{x})$ where $G=\iota(\text{Aut}_{e}(\Gamma_{r}))$ in a certain unknown time $T$, then from generators of $\ker\pi_{r}$ and $\text{im}\,\pi_{r}$ we build generators of $\text{Aut}_{e}(\Gamma_{r+1})$ in the following way: we extend $\sigma\in\text{im}\,\pi_{r}$ to the whole $\Gamma_{r+1}$ (arbitrarily: for example, if $\{x_{n_{i'}}\}_{i'\leq i}$ neighbors of $\{v_{n_{j'}}\}_{j'\leq j}$ have to go to $\{x_{n_{i''}}\}_{i''\leq i}$ neighbors of $\{\sigma(v_{n_{j'}})\}_{j'\leq j}$, we send the $x_{n_{i'}}$ with the lowest index to the $x_{n_{i''}}$ with the lowest index, etc...), we extend $\tau\in\ker\pi_{r}$ to the whole $\Gamma_{r+1}$ by setting $\tau=\text{Id}$ on $\Gamma_{r}$, and the generators of $\text{Aut}_{e}(\Gamma_{r+1})$ will be simply the collection of all these extended $\tau,\sigma$ (this process takes time $O(n^{2})$, as this is the bound on the number of generators of $\text{im}\,\pi_{r}$ that we can require, thanks to Schreier-Sims).

The recursion above constructs the whole $\text{Aut}_{e}(\Gamma)$ after $O(n)$ steps of induction. Therefore the entire procedure takes time at most:
\begin{equation*}
nk\left(n^{O(1)}+O(kn)+O(kn^{2})+O(n)\left(O(kn^{2})+n^{O(k)}+T+O(n^{2})\right)\right)\leq O(kn^{2}T+kn^{O(k)})
\end{equation*}
So for $k$ bounded we have reduced the graph isomorphism problem for graphs of maximum degree $k$ to the string isomorphism problem for $\text{Aut}_{G}(\mathbf{x})$, in polynomial time and with polynomial multiplicative cost; the string is of length $|K|\leq n^{k}$, i.e. polynomial length, and $G=\iota(\text{Aut}_{e}(\Gamma_{r}))\simeq\text{Aut}_{e}(\Gamma_{r})\in\mathfrak{G}_{k}$ by (\ref{hidden1b-AutG}), so we are done.
\bigskip

\noindent{\sc Solution to Ex.~\ref{hidden1d}} --- (\ref{hidden1d-conf}) Clearly, for any map $\tau:\{1,...,k\}\rightarrow\{1,...,k\}$, if the components of $\vec{x}$ are all inside $\Gamma'$ then the same is true for $\tau(\vec{x})$; then for $\mathfrak{X}[\Gamma']$ the maps $\rho,\eta$ are simply the restrictions of the maps $\rho,\eta$ of $\mathfrak{X}$. If we see it as a commutative diagram, for instance we have for $\eta$:
$$
\begin{array}{ccccc}
\Gamma'^{k} & \hookrightarrow & \Gamma^{k} & \stackrel{c}{\longrightarrow} & \mathscr{C} \\
{\scriptstyle \tau} \downarrow\ \ \ & \circlearrowleft & {\scriptstyle \tau} \downarrow\ \ \ & \circlearrowleft & {\scriptstyle \tau^{\eta}} \downarrow\ \ \ \\
\Gamma'^{k} & \hookrightarrow & \Gamma^{k} & \stackrel{c}{\longrightarrow} & \mathscr{C} \\
\end{array}
$$

(\ref{hidden1d-coh}) For $\Gamma'$ a color class, $\mathfrak{X}[\Gamma']$ is coherent (see Exercise~\ref{ex:rescoh}), but this is not true in general. Consider $\mathfrak{X}=(\Gamma,c)$ with $\Gamma=\mathbb{Z}/3\mathbb{Z}$ and coloring defined as $c(x,y)=x-y$: $\mathfrak{X}$ has two non-trivial automorphisms $\tau,\tau^{-1}$ where $\tau(x)=x+1$, and for every two pairs $(x,y),(x',y')$ with the same color we have either $\tau$ or $\tau^{-1}$ sending one to the other; therefore by Exercise~\ref{WLaut} there will be no refinement if we apply Weisfeiler-Leman, which means that $\mathfrak{X}$ is coherent. On the other hand, if $\Gamma'=\{0,1\}$ we have that $\gamma(1,-1,0)$ is $0$ for $(0,0)$ and $1$ for $(1,1)$, i.e. it is not independent anymore from the pair that we choose: $\mathfrak{X}[\Gamma']$ therefore is not coherent.
\bigskip

\noindent{\sc Solution to Ex.~\ref{hidden2}} --- (a) By the definition of configuration, for any pair of vertices $x,y\in\Gamma$, we have $c(x,y)=r$ if and only if $c(y,x)=r^{-1}$. Call $b$ the unique vertex color of our homogeneous configuration; then by coherence $\gamma(r^{-1},r,b)$ is a constant that does not depend on the choice of $x$ such that $c(x,x)=b$ (which is every possible vertex of $\Gamma$). So:
\begin{eqnarray*}
d^+_r(x) & = & |\{y\in \Gamma: (x,y)\in \mathscr{G}_r\}| \\ 
&= & |\{y\in \Gamma: c(x,y)=r\}| \\
 & = & |\{y\in \Gamma: c(y,x)=r^{-1},c(x,y)=r\}| \\ 
 &= & \gamma(r^{-1},r,b)
\end{eqnarray*}

(b) Suppose that there is a weakly connected component $B$ of $\mathscr{G}_r$ that is not strongly connected, i.e. $B=B_{1}\sqcup...\sqcup B_{m}$ where the $B_{i}$ are strongly connected components of $\mathscr{G}_r$ and $m>1$. For any two components $B_{i},B_{j}$, if there is an edge of color $r$ from some $x\in B_{i}$ to some $y\in B_{j}$ then there are no $r$-colored edges from $B_{j}$ to $B_{i}$: then we can define an oriented graph $\Gamma'$ with set of vertices $\{1,...,m\}$ and set of edges $\{(i,j)\in \{1,...,m\}^{2}: \exists \; r\mbox{-edge from }B_{i}\mbox{ to }B_{j}\}$; notice that in this new graph there are no (oriented) cycles, otherwise the components corresponding to the points in the cycle would be all strongly connected to each other.

Fix any vertex in $\Gamma'$ and start walking randomly through the edges of the graph, respecting their orientation: since $\Gamma'$ is a finite graph with no oriented cycle, our walk will necessarily end when we reach a vertex $i$ with no outgoing edges (in particular, this procedure shows that such a vertex exists). For every vertex $x$ inside the component $B_{i}$ and for every $(x,y)$ of color $r$, we must have $y\in B_{i}$, otherwise $i\in \Gamma'$ would have an outgoing edge to the component containing $y$; this means that for every $x\in B_{i}$ there are $d^+_r(x)$ edges of color $r$ starting from $x$ and entirely contained in $B_{i}$, and then $|B_{i}|d^+_r(x)$ $r$-colored edges inside $B_{i}$ in total.

Since $\Gamma$ is a configuration we have $d^+_{r^{-1}}(x)=d^+_r(x)$ for every vertex $x$; for the same reason we have $|B_{i}|d^+_r(x)$ edges of color $r^{-1}$ inside $B_{i}$, which means that {\em all} the edges of color $r^{-1}$ starting from vertices of $B_{i}$ are entirely contained in $B_{i}$: but this implies that $B_{i}$ is disconnected from the other components both through $r$-edges and through $r^{-1}$ edges, i.e. $B_{i}$ is not weakly connected to any other $B_{j}$. This contradicts our assumption that $B$ was weakly connected.
\bigskip

\noindent{\sc Solution to Ex.~\ref{hidden3}} --- We follow \cite{Ba}; many of the group-theoretic results used here can be found in \cite[Chapter 4]{DM}. Call $\text{Soc}(G)$ the socle of $G$; first we prove that $\text{Soc}(G)$ is not abelian. Supposing that it is, then $\text{Soc}(G)$ is the product of $s$ copies of $\mathbb{Z}/p\mathbb{Z}$ for $p$ such that $p^{s}=|\Omega|$, while $G/\text{Soc}(G)$ is a subgroup of $GL_{s}(\mathbb{Z}/p\mathbb{Z})$; also, using the fact that $\Alt_{k}$ is in this case a section of $GL_{s}(\mathbb{Z}/p\mathbb{Z})$, by \cite[Prop.~1.22a]{BaPS} we have $s\geq k-2$ ($k>8$ by hypothesis). But then we have $s\geq k-2>\log_2|\Omega|\geq\log_{p}|\Omega|=s$, which gives the contradiction that we were looking for. Thus, $\text{Soc}(G)$ is not abelian.

A consequence of the O'Nan-Scott Theorem is that since our $G$ has non-abelian $\text{Soc}(G)\simeq R^{s}$ ($\text{Soc}(G)$ can always be written in this way when $G<\Sym(\Omega)$ is primitive) we have $|\Omega|\geq 5^{s}$. Define $\psi:G\rightarrow\Sym_{s}$ by sending $g$ to the permutation of the copies of $R$ inside $\text{Soc}(G)$ given by conjugation by $g$ (we call these copies $R_{i}$): then $\ker\psi<\text{Aut}(R_{1})\times...\times\text{Aut}(R_{s})$ and $\ker\psi/\text{Soc}(G)<\text{Out}(R_{1})\times...\times\text{Out}(R_{s})$; by Schreier's conjecture each $\text{Out}(R_{i})$ is solvable, therefore $\ker\psi/\text{Soc}(G)$ is solvable too. As $\ker\psi/\text{Soc}(G)$ is solvable, it does not have $\Alt_{k}$ as a section; from $|\Omega|\geq 5^{s}$ we have $s<k$, so that $G/\ker\psi<\Sym_{s}$ cannot have $\Alt_{k}$ as a section either: thus, $G/\text{Soc}(G)$ itself has no section isomorphic to $\Alt_{k}$, which implies that $\text{Soc}(G)\not<\ker\phi$.

Consider any minimal $N\lhd G$. If $N\neq\text{Soc}(G)$, there is another (unique) minimal subgroup $N'$ which is the centralizer of $N$: $N'\simeq N$ by the structure of the socle, so that $N$ is regular and $|N|=|\Omega|$; we also have $s$ even and $|N|=|\Alt_{k}|^{s/2}$, therefore $|\Omega|=|\Alt_{k}|^{s/2}\geq|\Alt_{k}|>2^{k}>|\Omega|$, which is a contradiction: so $N=\text{Soc}(G)$, i.e. $\text{Soc}(G)$ is the unique minimal subgroup of $G$. $\text{Soc}(G)\not<\ker\phi$ then implies $\ker\phi=\{e\}$, and $\phi$ is an isomorphism.
\bigskip

\noindent{\sc Solution to Ex.~\ref{hidden4}} --- Call $b$ the unique vertex color: since there are edges of color $r$, the constant $\gamma(r^{-1},r,b)$ is nonzero, which means that every vertex is the starting point of some edge of color $r$. Fix $x\in\Gamma$: by our observation, not all vertices $y\neq x$ are neighbors of $x$ in $\overline{\mathscr{G}_r}$; given any edge color $r'\neq r$, by uniprimitivity $\mathscr{G}_{r'}$ is connected, which implies that $\overline{\mathscr{G}_r}$ is also connected: thus, there exists a vertex at distance $>1$ from $x$, and in particular there exists $y$ with $d(x,y)=2$. For this particular pair $(x,y)$ we have $c(x,y)=r$ and at least one vertex $z$ with $c(x,z),c(z,y)\neq r$: by coherence all $r$-colored edges will have a vertex $z'$ like this one; but then $x$ is at distance $1$ from any $w$ with $c(x,w)\neq r$ and at distance $2$ from any $w'$ with $c(x,w')=r$, and repeating the same reasoning for any $x$ we get that the diameter of $\overline{\mathscr{G}_r}$ is $2$.
\bigskip

\noindent{\sc Solution to Ex.~\ref{WLcoh}} --- (\ref{WLcoh-coh}) Suppose that we have reached the last iteration of Weisfeiler-Leman, i.e. at this point if $\vec{x},\vec{y}$ have the same color then the two vectors in the algorithm that refine $c(\vec{x})$ and $c(\vec{y})$ must also be the same; however, that is exactly what the definition of coherent configuration is.

(\ref{WLcoh-coarse}) Suppose that $\mathfrak{X}'=(\Gamma,c')$ is a coherent configuration whose coloring $c'$ is a refinement of the coloring of $\mathfrak{X}=(\Gamma,c)$: if we prove that every iteration of the Weisfeiler-Leman algorithm gives a coloring that is still less fine than $c'$, then the coherent configuration resulting after the end of the algorithm will also be coarser than $\mathfrak{X}'$, thus making it the coarsest of its kind. This is easy to see, though; if $c'(\vec{x})=c'(\vec{y})$ then we have that $c(\vec{x})=c(\vec{y})$, since $c'$ is a refinement of $c$, and for any $\vec{r'}\in\mathscr{C'}^{k}$:
\begin{equation*}
 \left|\left\{z\in \Gamma: c'(\vec{x}^{j}(z)) = r'_j\;\;\;\forall\; 1\leq j\leq k\right\}\right| =
 \left|\left\{z\in \Gamma: c'(\vec{y}^{j}(z)) = r'_j\;\;\;\forall\; 1\leq j\leq k\right\}\right|
\end{equation*}
by coherence. The colors $r'_{j}$ are refinements of the colors given by $c$, which means that the same equality holds for the appropriate coarser $c,\vec{r}$: so after the iteration of the algorithm, the new coloring is still the same for $\vec{x},\vec{y}$, which implies that $c'$ is finer than this new coloring too.
\bigskip

\noindent{\sc Solution to Ex.~\ref{canonord}} --- Lexicographical order solves all issues. Suppose that we have a $k$-ary relational structure with an ordered coloring $\mathscr{C}$; a partition structure is given canonically by $F_{1}$, where the set of colors is now the power set of $\mathscr{C}$: then we can associate to each new color a string of digits $0,1$ that correspond to old colors being included or not included in the new color we are considering (the ordering of $\mathscr{C}$ ensures that we actually have a string, rather than a jumble of digits). We can order these strings very naturally, for example lexicographically from largest to smallest according to the natural number they define.

Suppose now that we have a $k$-ary partition structure with an ordered coloring $\mathscr{C}$; a configuration is given canonically by $F_{2}$, where colors are given as strings $(\rho,c_{1},c_{2},...,c_{k^{k}})$ coming from $\left(\rho(\vec{x}),(c(\tau(\vec{x})))_{\tau}\right)$ (the $\tau$ are all the possible functions from $\{1,...,k\}$ to itself): the $\rho$ are equivalence relations on $\{1,...,k\}$, so there are $B_{k}$ of them ($B_{k}$ is the $k$-th Bell number) and they are naturally orderable, while the $c_{i}$ are inside $\mathscr{C}$. Again, we can order these strings lexicographically using the orderings on $\mathscr{C}$ and on the set of $\rho$.

Suppose finally that we have a $k$-ary configuration with an ordered coloring $\mathscr{C}$; a coherent configuration is given canonically through the algorithm of Weisfeiler-Leman (i.e. $F_{3}$): at every iteration, the new colors are given as strings $(c_{0},\gamma_{1},...,\gamma_{K})$ with $K=|\mathscr{C}|^{k}$. The numbers $\gamma_{i}$ are produced when considering elements of $\mathscr{C}^{k}$, so by lexicographically ordering this set we have a way to establish their order inside the string; the colors $c_{0}$ are ordered, being elements of $\mathscr{C}$, and the $\gamma_{i}$ are also ordered, being natural numbers: as before, we can order these strings lexicographically. This happens at every iteration of the algorithm, so in the end we have an ordering of the final coherent configuration.
\bigskip

\noindent{\sc Solution to Ex.~\ref{WLaut}} --- For any map $\tau:\{1,...,k\}\rightarrow\{1,...,k\}$, the components of $\tau(\vec{x})$ are a subset of the components of $\vec{x}$, so that it is clear that $\sigma(\tau(\vec{x}))=\tau(\sigma(\vec{x}))$; thus, using the fact that $\sigma\in\mbox{Aut}(\mathfrak{X})$ implies $c(\sigma(\vec{x}))=c(\vec{x})$, we get:
\begin{equation*}
c(\sigma(\tau(\vec{x})))=\tau^{\eta}(c(\sigma(\vec{x})))=\tau^{\eta}(c(\vec{x}))=c(\tau(\vec{x}))
\end{equation*}
which means that after $F_{2}$ the colors of $\sigma(\vec{x})$ and $\vec{x}$ are still the same.

Now suppose that $\mathfrak{X}=(\Gamma,c)$ is a $k$-ary configuration. For any $\vec{x}$ and any $z\in\Gamma$, the vector $\vec{x}^{j}(z)$ is sent to $\sigma(\vec{x})^{j}(\sigma(z))$, so they must have the same color because $\sigma$ is an automorphism of $\mathfrak{X}$: this means that for every $(k+1)$-tuple of vectors given by $\vec{x},z$ considered in the algorithm that has a certain coloring $(\vec{r},j)$, the images $\sigma(\vec{x}),\sigma(z)$ give a $(k+1)$-tuple of vectors with the same coloring. In particular,
\begin{equation*}\begin{aligned}
&\left|\left\{z\in \Gamma: c(\vec{x}^{j}(z)) = r_j\;\;\;\forall\; 1\leq j\leq k\right\}\right|\\ &=
\left|\left\{\sigma(z)\in \Gamma: c(\sigma(\vec{x})^{j}(\sigma(z))) = r_j\;\;\;\forall\; 1\leq j\leq k\right\}\right|
\end{aligned}\end{equation*}
so at every step $\vec{x},\sigma(\vec{x})$ have the same color; then, at the end of $F_{3}$ we will still have $c'(\sigma(\vec{x}))=c'(\vec{x})$.
\bigskip

\noindent{\sc Solution to Ex.~\ref{graphconfig}} --- A pair $(x,y)$ is colored ``vertex'' if and only if $x=y$, so by our construction ``vertex'' is a vertex color and ``edge'', ``not edge'' are edge colors. Since the graph is undirected, $(x,y)$ is an edge if and only if $(y,x)$ is an edge, so for the four functions:
\begin{equation*}
(1,2)\stackrel{\tau_{1}}{\mapsto}(1,2) \ \ \ \ \ (1,2)\stackrel{\tau_{2}}{\mapsto}(2,1) \ \ \ \ \ (1,2)\stackrel{\tau_{3}}{\mapsto}(1,1) \ \ \ \ \ (1,2)\stackrel{\tau_{4}}{\mapsto}(2,2)
\end{equation*}
we have:
\begin{eqnarray*}
(\text{vertex},\text{edge},\text{not edge}) & \stackrel{\tau_{1}^{\eta},\tau_{2}^{\eta}}{\mapsto} & (\text{vertex},\text{edge},\text{not edge}) \\
(\text{vertex},\text{edge},\text{not edge}) & \stackrel{\tau_{3}^{\eta},\tau_{4}^{\eta}}{\mapsto} & (\text{vertex},\text{vertex},\text{vertex})
\end{eqnarray*}
Therefore $(\Gamma,c)$ is a configuration.

(\ref{graphconfig-Pn}) Say that we have vertices $x_{1},x_{2},...,x_{n}$ with edges $\{x_{i},x_{i+1}\}$ for all $1\leq i<n$: then for $n\geq 3$ the vertex $x_{1}$ has only one edge while $x_{2}$ has two, so that the constant $\gamma(\text{edge},\text{edge},\text{vertex})$ is $1$ for $(x_{1},x_{1})$ and $2$ for $(x_{2},x_{2})$; thus $P_{n}$ is not coherent for $n\geq 3$.

(\ref{graphconfig-C4}) Say that we have vertices $x_{1},x_{2},x_{3},x_{4}$ with edges $\{x_{i},x_{i+1}\}$ for $1\leq i<4$ and $\{x_{4},x_{1}\}$. The cycle $C_{4}$ has 4 pairs colored ``vertex'', 8 pairs colored ``edge'' and 4 pairs colored ``not edge'': if we can show that for any two pairs of the same color there is an automorphism of $C_{4}$ that sends one to the other, by Exercise~\ref{WLaut} we know that the Weisfeiler-Leman algorithm will not refine the coloring at all, which means that $C_{4}$ is already coherent.

Define $\sigma,\tau:C_{4}\rightarrow C_{4}$ as:
\begin{eqnarray*}
(x_{1},x_{2},x_{3},x_{4}) & \stackrel{\sigma}{\mapsto} & (x_{2},x_{3},x_{4},x_{1}) \\
(x_{1},x_{2},x_{3},x_{4}) & \stackrel{\tau}{\mapsto} & (x_{1},x_{4},x_{3},x_{2})
\end{eqnarray*}
Basically $\sigma$ is a rotation and $\tau$ a reflection, obviously automorphisms of $C_{4}$; moreover, it is easy to see that the collection of seven automorphisms $\tau,\tau^{2},\tau^{3},\sigma,\tau\sigma,\tau^{2}\sigma,\tau^{3}\sigma$ satisfies the requirements that we need to conclude the proof. Thus $C_{4}$ is coherent.

(\ref{graphconfig-C6}) Say that we have vertices $x_{1},x_{2},...,x_{6}$ with edges $\{x_{i},x_{i+1}\}$ for $1\leq i<6$ and $\{x_{6},x_{1}\}$: then $d(x_{1},x_{3})=2$ and $d(x_{1},x_{4})>2$, so that the constant $\gamma(\text{edge},\text{edge},\text{not edge})$ is $1$ for $(x_{1},x_{3})$ and $0$ for $(x_{1},x_{4})$; thus $C_{6}$ is not coherent.
\bigskip

\noindent{\sc Solution to Ex.~\ref{stronglyreg}} --- Let us be given two vertices $x,y\in\Gamma$ and call $d(x),d(y),N(x,y)$ respectively the number of neighbors of $x$, of $y$ and of both $x$ and $y$; then:
\begin{eqnarray*}
|\{z\in\Gamma: c(z,y)=\text{edge},c(x,z)=\text{edge}\}| & = & N(x,y) \\
|\{z\in\Gamma: c(z,y)=\text{edge},c(x,z)=\text{not edge}\}| & = & d(y)-N(x,y)-\varepsilon \\
|\{z\in\Gamma: c(z,y)=\text{not edge},c(x,z)=\text{edge}\}| & = & d(x)-N(x,y)-\varepsilon \\
|\{z\in\Gamma: c(z,y)=\text{not edge},c(x,z)=\text{not edge}\}| & = & |\Gamma|-2(1-\varepsilon)\\ & &-d(x)-d(y)+N(x,y)
\end{eqnarray*}
where $\varepsilon=1$ if $(x,y)$ is an edge and $\varepsilon=0$ otherwise.

Suppose that $\Gamma$ is strongly regular. First we calculate the easy constants coming from impossible constructions or obvious identities:
\begin{eqnarray*}
\gamma(\text{vertex},\text{edge},\text{vertex})=0 & \ \ \ & \gamma(\text{edge},\text{vertex},\text{vertex})=0 \\
\gamma(\text{vertex},\text{not edge},\text{vertex})=0 & \ \ \ & \gamma(\text{not edge},\text{vertex},\text{vertex})=0 \\
\gamma(\text{edge},\text{not edge},\text{vertex})=0 & \ \ \ & \gamma(\text{not edge},\text{edge},\text{vertex})=0 \\
\gamma(\text{vertex},\text{edge},\text{edge})=1 & \ \ \ & \gamma(\text{edge},\text{vertex},\text{edge})=1 \\
\gamma(\text{vertex},\text{not edge},\text{not edge})=1 & \ \ \ & \gamma(\text{not edge},\text{vertex},\text{not edge})=1 \\
\gamma(\text{vertex},\text{not edge},\text{edge})=0 & \ \ \ & \gamma(\text{not edge},\text{vertex},\text{edge})=0 \\
\gamma(\text{vertex},\text{edge},\text{not edge})=0 & \ \ \ & \gamma(\text{edge},\text{vertex},\text{not edge})=0 \\
\gamma(\text{vertex},\text{vertex},\text{edge})=0 & \ \ \ & \gamma(\text{vertex},\text{vertex},\text{not edge})=0 \\
\gamma(\text{vertex},\text{vertex},\text{vertex})=1 & \ \ \ & \
\end{eqnarray*}
Then, since the graph is $k$-regular:
\begin{eqnarray*}
\gamma(\text{edge},\text{edge},\text{vertex})=k & \ \ \ & \gamma(\text{not edge},\text{not edge},\text{vertex})=|\Gamma|-1-k
\end{eqnarray*}
Finally, from the relations we established before:
\begin{eqnarray*}
\gamma(\text{edge},\text{edge},\text{edge})=\lambda & \ \ \ & \gamma(\text{edge},\text{not edge},\text{edge})=k-\lambda-1 \\
\gamma(\text{not edge},\text{edge},\text{edge})=k-\lambda-1 & \ \ \ & \gamma(\text{not edge},\text{not edge},\text{edge})=|\Gamma|-2k+\lambda \\
\gamma(\text{edge},\text{edge},\text{not edge})=\mu & \ \ \ & \gamma(\text{edge},\text{not edge},\text{not edge})=k-\mu \\
\gamma(\text{not edge},\text{edge},\text{not edge})=k-\mu & \ \ \ & \gamma(\text{not edge},\text{not edge},\text{not edge})\\ & &=|\Gamma|-2-2k+\mu
\end{eqnarray*}
So $(\Gamma,c)$ is a coherent configuration.

Suppose now that $(\Gamma,c)$ is a coherent configuration. Then it is easy to find $k,\lambda,\mu$:
\begin{eqnarray*}
k & = & \gamma(\text{edge},\text{edge},\text{vertex}) \\
\lambda & = & \gamma(\text{edge},\text{edge},\text{edge}) \\
\mu & = & \gamma(\text{edge},\text{edge},\text{not edge})
\end{eqnarray*}
and their independence from $(x,y)$ is guaranteed by the independence of the constants $\gamma$. So $\Gamma$ is strongly regular.
\bigskip

\noindent{\sc Solution to Ex.~\ref{symmBIBD}} --- First we show that $(\Gamma,c)$ is a coherent configuration; the fact that it is a configuration descends easily from the uniqueness of the vertex color ``vertex'' and from the fact that by construction every color is its own inverse. The verification of the existence of constants $\gamma=0,1$ for impossible constructions or obvious identities (as in Exercise~\ref{stronglyreg}) is trivial and we skip it here; $48$ of the $64$ combinations are of this type. Also, from now on we call the colors $V,W,B,D$ for the sake of brevity.

Since $v=b$ we have for both $V$ and $\mathscr{A}$:
\begin{eqnarray*}
\gamma(W,W,V)=v-1 & \ \ \ & \gamma(W,W,W)=v-2
\end{eqnarray*}
Every block has $u$ elements; any vertex in a BIBD belongs by definition to the same number of blocks $r$, which implies $rv=bu$ and then for a symmetric BIBD $r=u$, so for both $V$ and $\mathscr{A}$:
\begin{eqnarray*}
\gamma(B,B,V)=u & \ \ \ & \gamma(D,D,V)=v-u
\end{eqnarray*}
For the same reason, for any choice of $(x,y)\in V\times\mathscr{A}$, if $x\in y$ then we have exactly $u-1$ other $x'$ belonging to $y$ and exactly $u-1$ other $y'$ to which $x$ belongs, while if $x\not\in y$ the numbers are both $u$; therefore:
\begin{eqnarray*}
\gamma(W,B,B)=\gamma(B,W,B)=u-1 & \ \ \ & \gamma(W,D,B)=\gamma(D,W,B)=v-u \\
\gamma(W,B,D)=\gamma(B,W,D)=u & \ \ \ & \gamma(W,D,D)=\gamma(D,W,D)=v-u-1
\end{eqnarray*}
Since the BIBD is symmetric, any two vertices belong to exactly $\lambda$ blocks at the same time, and the intersection of any two blocks has size $\lambda$; this, and the fact that $v=b$ and $r=u$, gives:
\begin{eqnarray*}
\gamma(B,B,W)=\lambda & \ \ \ & \gamma(D,D,W)=v-2u+\lambda \\
\gamma(B,D,W)=\gamma(D,B,W)=u-\lambda & \ \ \ & \
\end{eqnarray*}
Thus $(\Gamma,c)$ is a coherent configuration when the BIBD is symmetric.

Suppose that the BIBD is not symmetric anymore. Then we have $b>v$, and $\gamma(W,W,V)$ would be $v-1$ for $x\in V$ and $b-1$ for $x\in\mathscr{A}$: so the configuration is not coherent.

Suppose that we have refined the coloring of $\Gamma$ as described, so that now we have $8$ colors: $V_{V},V_{\mathscr{A}},W_{V},W_{\mathscr{A}},B,B^{-1},D,D^{-1}$ (in particular, this repairs the problem of $v-1$ being different from $b-1$, since now we were referring to $\gamma(W_{V},W_{V},V_{V})$ and $\gamma(W_{\mathscr{A}},W_{\mathscr{A}},V_{\mathscr{A}})$). Nevertheless, there exist two pairs of blocks $(y_{1},y_{2}),(y'_{1},y'_{2})$ with two different intersection sizes: this means that they give different values for $\gamma(B,B^{-1},W_{\mathscr{A}})$ despite being of the same color $W_{\mathscr{A}}$; therefore the configuration is still not coherent.
\bigskip

\noindent{\sc Solution to Ex.~\ref{tdesign}} --- (\ref{tdesign-b}) Every $t$-subset of $V$ is contained in $\lambda$ edges and every edge contains $\binom{u}{t}$ such subsets: so counting pairs $(T,A)$ in two different ways we obtain $\binom{v}{t}\lambda=b\binom{u}{t}$, from which $b=\lambda\binom{v}{t}\binom{u}{t}^{-1}$.

(\ref{tdesign-t'}) Suppose that $(V,\mathscr{A})$ is a $t$-design, fix a $(t-1)$-subset $X$ and call $\lambda_{t-1}$ the number of edges containing $X$: we are going to count in two different ways the number of pairs $(A,x)\in\mathscr{A}\times V$ such that $x\not\in X$ and $A\supset X\cup\{x\}$. On one side, fixing $A$ we have $u-t+1$ possible choices of $x$; on the other side, fixing $x$ we have $\lambda$ choices of $A$: therefore $\lambda_{t-1}(u-t+1)=(v-t+1)\lambda$, from which $\lambda_{t-1}=\lambda\frac{v-t+1}{u-t+1}$ and since it is independent from our initial choice of $X$ the $t$-design is also a $(t-1)$-design.

Iterating the procedure, $(V,\mathscr{A})$ is a $t'$-design for any $t'\leq t$ with:
\begin{equation*}
\lambda_{t'}=\lambda\frac{v-t+1}{u-t+1}\frac{v-t+2}{u-t+2}\ldots\frac{v-t'}{u-t'}=\lambda\frac{(v-t')!}{(v-t)!}\frac{(u-t)!}{(u-t')!}=\lambda\binom{v-t'}{t-t'}\binom{u-t'}{t-t'}^{-1}
\end{equation*}

(\ref{tdesign-ij}) By (\ref{tdesign-t'}) our design is also an $(i+j')$-design for any $0\leq j'\leq j$. We reason by inclusion-exclusion: we count all edges containing $I$, we subtract all those containing $I\cup\{j_{1}\}$ for $j_{1}\in J$, we re-add those containing $I\cup\{j_{1},j_{2}\}$, etc...; doing this, we obtain:
\begin{equation*}
\lambda_{i,j}=\sum_{j'=0}^{j}(-1)^{j'}\binom{j}{j'}\lambda_{i+j'}=\lambda\sum_{j'=0}^{j}(-1)^{j'}\binom{j}{j'}\binom{v-i-j'}{t-i-j'}\binom{u-i-j'}{t-i-j'}^{-1}
\end{equation*}
Instead of solving this expression, we notice that $\lambda_{i,j}\lambda^{-1}$ does not depend on the particular design that we are working with, but depends only on its parameters: any other design with the same $v,u,t,i,j$ would give the same $\lambda_{i,j}$; defining $\mathscr{A}'$ to be the set of all $u$-subsets of $V$ we obtain trivially a $u$-design $(V,\mathscr{A}')$ (and then also a $t$-design) for which $\lambda'$ is equal to the number of all $U$-subsets containing a certain $t$-subset and $\lambda'_{i,j}$ is equal to the number of all $u$-subsets containing $I$ and disjoint from $J$, which implies:
\begin{equation*}
\lambda_{i,j}=\lambda\lambda'_{i,j}\lambda'^{-1}=\lambda\binom{v-i-j}{u-i}\binom{v-t}{u-t}^{-1}
\end{equation*}

(\ref{tdesign-FasE}) Any $S_{1}\in E_{i}$ appears in $F_{j}$ as many times as the number of $A$ such that $A\supset S_{1}$ and $|A\cap S_{0}|=s-j$; since $|S_{1}\cap S_{0}|=s-i$, when $j>i$ this number is $0$, while when $j\leq i$ it is $\binom{i}{j}\lambda_{s+i-j,j}$: so we have $F_{j}=\sum_{i=j}^{s}\binom{i}{j}\lambda_{s+i-j,j}E_{i}$.

(\ref{tdesign-EasF}) We have expressed the $F_{j}$ as linear combinations of the $E_{i}$: the linear system is upper triangular and each of the coefficients for $i=j$ is $\lambda_{s,i}\neq 0$ because $v\geq s+u\geq i+u$ implies $v-s-i\geq u-s$. Thus the system is invertible and we can express the $E_{i}$ as linear combinations of the $F_{j}$.

(\ref{tdesign-concl}) By (\ref{tdesign-EasF}), we can express every $E_{i}$ as combinations of the $\hat{A}$; $E_{0}=S_{0}$, so doing it for any $S_{0}$ we get that the $\hat{A}$ span the whole vector space generated by the $s$-subsets: in particular, the number of different $\hat{A}$ must be at least the dimension of the vector space, so $b\geq\binom{v}{s}$.

(\ref{tdesign-other}) For $t$ odd, $t=2s+1$, the $t$-design is also a $2s$-design by (\ref{tdesign-t'}), so we still have $b\geq\binom{v}{s}$.

Suppose that $v<s+u$ (and then $v-u<t$, whether $t=2s$ or $t=2s+1$): for every $(v-u)$-subset of $V$ by (\ref{tdesign-ij}) we have the same number of edges that are disjoint from it (using $I=\emptyset$); this number by construction can only be $0$ or $1$. If it is $0$, then we have $0=\lambda\binom{v-t}{u-t}^{-1}$ and then $\lambda=0$, which means that the design has no edges at all (in this case $b=0$); if it is $1$, then for every $(v-u)$-subset its complement is in $\mathscr{A}$, i.e. our design is the one formed by all possible $u$-subsets of $V$ (in this case $b=\binom{v}{u}$).
\bigskip

\noindent{\sc Solution to Ex.~\ref{johnson}} --- The intersections on the diagonal are the trivial ones, so we say $\iota'(x)=2$; there are $90$ pairs of color $y$ and $120$ pairs of color $z$, so we say $\iota'(y)=0$ and $\iota'(z)=1$ (remember, we swapped the definitions in this particular case).

Let us start with $v_{1}$. We have:
\begin{eqnarray*}
B(v_{1},v_{2})=\{v_{8},v_{10},v_{15}\} & \ \ B(v_{1},v_{5})=\{v_{9},v_{13},v_{14}\} \ \ & B(v_{1},v_{8})=\{v_{2},v_{9},v_{13}\} \\
B(v_{1},v_{9})=\{v_{5},v_{8},v_{10}\} & \ \ B(v_{1},v_{10})=\{v_{2},v_{9},v_{14}\} \ \ & B(v_{1},v_{13})=\{v_{5},v_{8},v_{15}\} \\
B(v_{1},v_{14})=\{v_{5},v_{10},v_{15}\} & \ \ B(v_{1},v_{15})=\{v_{2},v_{13},v_{14}\} \ \ & \ 
\end{eqnarray*}
and then:
\begin{eqnarray*}
C(v_{1},v_{2})=C(v_{1},v_{9})=C(v_{1},v_{13})=C(v_{1},v_{14})=\{v_{1},v_{5},v_{8},v_{10},v_{15}\}=\lambda_{1} \\
C(v_{1},v_{5})=C(v_{1},v_{8})=C(v_{1},v_{10})=C(v_{1},v_{15})=\{v_{1},v_{2},v_{9},v_{13},v_{14}\}=\lambda_{2}
\end{eqnarray*}
We have $\iota(v_{1})=\{\lambda_{1},\lambda_{2}\}$ (obviously, since we started with $v_{1}$) and $\lambda_{1}$ is the element contained in the images of $v_{1},v_{5},v_{8},v_{10},v_{15}$ (and analogously for $\lambda_{2}$). Doing the same for the other vertices, in the end we get:
\begin{eqnarray*}
\lambda_{1}=\{v_{1},v_{5},v_{8},v_{10},v_{15}\} & \ \ \lambda_{2}=\{v_{1},v_{2},v_{9},v_{13},v_{14}\} \ \ & \lambda_{3}=\{v_{2},v_{4},v_{5},v_{11},v_{12}\} \\
\lambda_{4}=\{v_{3},v_{7},v_{10},v_{12},v_{13}\} & \ \ \lambda_{5}=\{v_{3},v_{6},v_{8},v_{11},v_{14}\} \ \ & \lambda_{6}=\{v_{4},v_{6},v_{7},v_{9},v_{15}\}
\end{eqnarray*}
\begin{eqnarray*}
 & & v_{1}\stackrel{\iota}{\mapsto}\{\lambda_{1},\lambda_{2}\} \ \ \ \ v_{2}\stackrel{\iota}{\mapsto}\{\lambda_{2},\lambda_{3}\} \ \ \ \ v_{3}\stackrel{\iota}{\mapsto}\{\lambda_{4},\lambda_{5}\} \ \ \ \ v_{4}\stackrel{\iota}{\mapsto}\{\lambda_{3},\lambda_{6}\} \ \ \ \ v_{5}\stackrel{\iota}{\mapsto}\{\lambda_{1},\lambda_{3}\} \\
 & & v_{6}\stackrel{\iota}{\mapsto}\{\lambda_{5},\lambda_{6}\} \ \ \ \ v_{7}\stackrel{\iota}{\mapsto}\{\lambda_{4},\lambda_{6}\} \ \ \ \ v_{8}\stackrel{\iota}{\mapsto}\{\lambda_{1},\lambda_{5}\} \ \ \ \ v_{9}\stackrel{\iota}{\mapsto}\{\lambda_{2},\lambda_{6}\} \ \ \ \ v_{10}\stackrel{\iota}{\mapsto}\{\lambda_{1},\lambda_{4}\} \\
 & & v_{11}\stackrel{\iota}{\mapsto}\{\lambda_{3},\lambda_{5}\} \ \ \ v_{12}\stackrel{\iota}{\mapsto}\{\lambda_{3},\lambda_{4}\} \ \ \ v_{13}\stackrel{\iota}{\mapsto}\{\lambda_{2},\lambda_{4}\} \ \ \ v_{14}\stackrel{\iota}{\mapsto}\{\lambda_{2},\lambda_{5}\} \ \ \ v_{15}\stackrel{\iota}{\mapsto}\{\lambda_{1},\lambda_{6}\} \\
\end{eqnarray*}
\bigskip

\noindent{\sc Solution to Ex.~\ref{notprimcoh}} --- Consider any edge color $r'$: by Exercise~\ref{ex:samesize}\ref{it:bibr1}, $r'$ knows whether there is a walk of some length between its vertices that consist of all $r$-edges. If such a walk exists, every two vertices connected by an $r'$-edge are also connected by an $r$-walk: therefore every connected component of $\mathscr{G}_{r'}$ is a subset of a connected component $B_{i}$ of $\mathscr{G}_{r}$; if there is no such walk, no two vertices of an $r'$-edge are in the same $B_{i}$, which means that the complete $m$-partite graph built among the $B_{i}$ contains all the $r'$-edges.
\bigskip

\noindent{\sc Solution to Ex.~\ref{squaredom}} --- By Exercise~\ref{hidden4}, for any edge of color $r$ there exist two edges of color $r',r''$ (possibly $r'=r''$, but $r',r''\neq r$) that form a walk between the two vertices of the $r$-edge: therefore, starting from any fixed $x\in\Gamma$, the endpoints of all such walks contain all endpoints of $r$-edges starting from $x$. There are $\gamma_{r'}$ $r'$-colored edges from $x$, and from any of their endpoints there are $\gamma_{r''}$ $r''$-colored edges: therefore the number of possible walks is $\leq\gamma_{r'}\gamma_{r''}$ (some of them could have the same endpoint); we have then $\gamma_{r}\leq\gamma_{r'}\gamma_{r''}$, and for one of the two $r',r''$ we must have $\gamma_{r'}\geq\sqrt{\gamma_{r}}$ or $\gamma_{r''}\geq\sqrt{\gamma_{r}}$.
\bigskip

\noindent{\sc Solution to Ex.~\ref{colpart}} --- By Exercise~\ref{ex:samesize}\ref{it:bibr1}, in a coherent configuration any edge color knows the length of the shortest red walk between its vertices; moreover, in a configuration any edge color knows the color of both its vertices, and in a coherent configuration every vertex color knows how many edges of a given color are going from or to a vertex of that color: we use these three facts to prove our statement.

We color $v$ green and we apply the functors $F_{2},F_{3}$: for any $(v,w)$, the color $c(v,w)$ now knows that the edge is departing from the (only) green vertex and it also knows the shortest red walk from this green vertex to the endpoint; therefore, the color of $(v,w)$ knows $d(v,w)$. As we said, there is only one green vertex, so for any $w$ there is only one edge going to $w$ that has one of these edge colors that we described: according to the $d(v,w)$ that is encoded inside this unique color, the vertex $w$ will get a different color through $F_{3}$, which means that $w$ knows $d(v,w)$ too; so $\Gamma$ has a coloring of the vertices that is at least as fine as the partition into $D_{i}=\{w\in\Gamma: d(v,w)=i\}$: now we only need to bound their size (notice that we do not have a $D_{\infty}$ since $\mathscr{G}_{\text{red}}$ is connected).

For $i=0$ we have $|D_{0}|=|\{v\}|=1\leq\frac{1}{2}|\Gamma|$ (and $\gamma\geq 1$ concludes the argument). For every $i\geq 0$ we also have the bound $|D_{i+1}|\leq\gamma|D_{i}|$, because every $w\in D_{i}$ gives at most $\gamma$ vertices $w'\in D_{i+1}$ by the definition of $\gamma$; thanks to this, we have:
\begin{equation*}
|D_{i+1}|\leq\gamma|D_{i}|\leq\gamma(|\Gamma|-|D_{i+1}|) \Rightarrow |D_{i+1}|\leq\frac{\gamma}{\gamma+1}|\Gamma|=\left(1-\frac{1}{\gamma+1}\right)|\Gamma|
\end{equation*}
Therefore the new coloring is a $(1-(\gamma+1)^{-1})$-coloring.


\end{document}